\documentclass[aos,preprint]{imsart}
\usepackage{pdfpages}
\usepackage{amssymb}
\usepackage{import}
\usepackage{epsfig}
\usepackage{standalone}
\usepackage{amsthm}
\usepackage{booktabs}
\usepackage{amsmath}
\usepackage{appendix}
\usepackage{amsfonts}
\usepackage{color}
\usepackage{multirow}
\usepackage{float}
\usepackage{xr}
\RequirePackage[OT1]{fontenc} \RequirePackage{amsthm,amsmath}
\RequirePackage[numbers]{natbib}
\RequirePackage[colorlinks,citecolor=blue,urlcolor=blue]{hyperref}
\RequirePackage{hypernat} \makeatletter
\newcommand{\field}[1]{\mathbb{#1}}
\def\X{{\mathcal{X}}}
\def\Y{{\mathcal{Y}}}
\def\W{{\hat{W}}}
\newcommand{\Z}{\mathcal{Z}}
\newcommand{\R}{\field{R}}
\newcommand{\E}{\field{E}}
\newcommand*{\I}{\imath}
\def\qed{\hfill$\diamondsuit$}

\theoremstyle{example} \theoremstyle{remark} \theoremstyle{lemma}
\theoremstyle{definition} \theoremstyle{corol}
\theoremstyle{proposition} \theoremstyle{condition}
\theoremstyle{assumption}
\newtheorem{assumption}{\n{Assumption}}[section]
\newtheorem{theorem}{\n{Theorem}}[section]
\newtheorem{example}{\n{Example}}[section]
\newtheorem{remark}{\n{Remark}}[section]

\newtheorem{definition}{\n{Definition}}[section]

\newtheorem{proposition}{\n{Proposition}}[section]

\font\n=cmcsc10
\def\cov{{\mbox{cov}}}
\def\var{{\mbox{var}}}

\makeatletter

\newcommand{\Rmnum}[1]{\expandafter\@slowromancap\romannumeral #1@}

\begin{document}
\begin{frontmatter}
\title{Conditional Mean and Quantile Dependence Testing in High Dimension} \runtitle{Conditional Dependence Testing in High Dimension}
\begin{aug}
\author{\fnms{Xianyang}
\snm{Zhang}\ead[label=e1]{zhangxiany@stat.tamu.edu}},
\author{\fnms{Shun}
\snm{Yao}\ead[label=e2]{shunyao2@illinois.edu}}, and
\author{\fnms{Xiaofeng}
\snm{Shao}\ead[label=e3]{xshao@illinois.edu}}\thanks{Zhang's research
is partially supported by NSF DMS-1607320. Shao's research
is partially supported by NSF-DMS1104545, NSF-DMS1407037 and DMS-1607489.} \runauthor{X. Zhang et al.} \affiliation{Texas A\&M University and University of Illinois at Urbana-Champaign}
\address{
X. Zhang\\
Department of Statistics\\
Texas A\&M University\\
College Station, TX 77843.\\
\printead{e1}\\
\phantom{E-mail:\ }}

\address{
S. Yao\\
Department of Statistics\\
University of Illinois at Urbana-Champaign\\
Champaign, IL 61820.\\
\printead{e2}\\
\phantom{E-mail:\ }}

\address{
X. Shao\\
Department of Statistics\\
University of Illinois at Urbana-Champaign\\
Champaign, IL 61820.\\
\printead{e3}\\ \phantom{E-mail:\ }}
\end{aug}

\begin{abstract}
Motivated by applications in biological science, we propose a novel
test to assess the conditional mean dependence of a response
variable on a large number of covariates. Our procedure is built on
the martingale difference divergence recently proposed in Shao and
Zhang (2014), and it is able to detect certain type of departure from the
null hypothesis of conditional mean independence without making any
specific model assumptions. Theoretically, we establish the
asymptotic normality of the proposed test statistic under suitable
assumption on the eigenvalues of a Hermitian operator, which is
constructed based on the characteristic function of the covariates.
These conditions can be simplified under banded dependence structure on the covariates or Gaussian design.
To account for heterogeneity within the data, we further develop a
testing procedure for conditional quantile independence at a given
quantile level and provide an asymptotic justification. Empirically,
our test of conditional mean independence delivers comparable
results to the competitor, which was constructed under the linear
model framework, when the underlying model is linear. It
significantly outperforms the competitor when the conditional mean
admits a nonlinear form.
%The proposed test is employed to analyze a
%microarray data set on Yorkshire Gilts to detect significant gene
%ontology terms which are significantly associated with the thyroid hormone.
\end{abstract}


\end{frontmatter}

\section{Introduction}
Estimation and inference for regression models is of central
importance in statistics. For a response variable $Y\in \R$ and a
set of covariates $X\in \R^p$, it would be desirable to  determine
whether $X$
 is useful in modeling certain aspect of the response, such as
 the mean or quantiles of $Y$, even before constructing a parametric/nonparametric model.
  The main thrust of this article is to introduce new tests for conditional mean and quantile dependence in high dimension, which allow
  the dimension $p$ to be much larger than sample size $n$.
 Our dependence testing problem in the high dimensional setting is well motivated by an important problem in biological science, which is to
  test for the significance of a gene set and  identify significant sets of genes that
 are associated with certain clinical outcomes; see Subramanian et al. (2005),
 Efron and Tibshirani (2007) and Newton et al. (2007), among others. Since the size of a gene set can range from a few to thousands,
 this calls for a simultaneous test that can accommodate high dimensional covariates and possible dependence within the covariates.

%In statistical theory, there has been
%many developments that adopts the so-called "large-$p$, small-$n$" paradigm in recent theoretical studies for
% statistical estimation and testing; see Kosorok and Ma (2007), Fan, Hall and Yao (2007), Chen and Qin (2010), among others.
% New insights and phenomenon that cannot be uncovered under the traditional "fixed-$p$, large-$n$" setting
% have  been found lately in contemporary statistics.

 Recently, Zhong and Chen (2011) proposed  a simultaneous test for coefficient in high dimensional linear regression. Specifically,
 for a univariate response $Y$ and high dimensional covariates $X\in \R^p$, they assume $\E(Y|X)=\alpha+X^T\beta$ and $\var(Y|X)=\sigma^2$,
 and test $H_0:\beta=\beta_0$ versus $H_a:\beta\not=\beta_0$ for a specific $\beta_0\in \R^p$ based on a random sample from the joint distribution
 of $(X,Y)$. A special case of primary interest is when $\beta_0=0$, which indicates insignificance of all covariates.
 However, the assumption of a high dimensional linear model seems quite strong and we are not aware of any procedure to validate this linear model
 assumption in the high dimensional setting. To assess the usefulness of the covariates $X$ in modeling the mean of $Y$, we effectively want to test
 \begin{align*}
H_0: \E(Y|X)=\E(Y)~\mbox{almost surely},~\mbox{versus}~H_a:
P(\E(Y|X)=\E(Y))<1,
 \end{align*}
 in other words, the conditional mean independence of $Y$ on $X$. Thus it seems desirable to develop a test that can accommodate the high dimensionality
and dependence in the covariates without assuming a linear (or
parametric) model. It turns out that the above testing problem in a model-free setting  is very challenging since the class of
alternative we target is huge  owing to the  growing dimension and nonlinear dependence.

To circumvent the difficulty, we
focus on testing the weak null hypothesis:
\[
H_0': \E[Y|x_j]=\E[Y]\quad  \text{almost surely},~\text{for all} ~1\leq j\leq p, \]
where $X=(x_1,\cdots,x_p)^T$, which is itself an interesting and meaningful testing problem; see Section \ref{sec:mean} for more discussions.
 Since $H_0$ implies $H_0'$, a rejection of $H_0'$ automatically rejects $H_0$.
Furthermore, $H_0$ and $H_0'$ are equivalent in the framework of high dimensional linear models under some mild assumptions; see Example \ref{linear}.
In this paper, we  propose a new test for $H_0'$, i.e., the conditional mean independence
for each $x_j$,  that can allow high dimensionality and dependence in
$X$ without any parametric or homoscedastic model assumption. Our
test statistic is built on the so-called martingale difference
divergence (MDD, hereafter) [Shao and Zhang (2014)], which fully
characterizes the conditional mean independence of a univariate
response on the covariates of arbitrary dimension. For any random variable $\Y$ and random vector
  $\X$, the MDD is defined as the weighted $L_2$ norm of $\cov(\Y,e^{\I<s,\X>})$ with $\I=\sqrt{-1}$ and the weighting function
  similar to the one used in distance covariance [Sz\'ekely, Rizzo and Bakirov (2007)], and thus inherits
  many desirable properties of distance covariance; see Shao and Zhang (2014) and Park et al. (2015) for more details.
   Theoretically, we
establish the asymptotic normality of our MDD-based test statistic
under suitable assumption on the eigenvalues of a Hermitian operator
constructed based on the characteristic functions of the covariates.
These conditions can be further simplified under banded dependence structure on the covariates or
Gaussian design. The theoretical results are of independent interest since they
provide some insights on the asymptotic behavior of $U$-statistic in
high dimension. From a practical viewpoint, our test does not
involve any tuning parameter and uses standard $z$ critical value,
so it can be conveniently implemented. To reduce the finite sample size distortion, we further propose a studentized wild
bootstrap approach and justify its consistency. Simulation results show that the wild bootstrap leads to more accurate size
for most cases, at the expense of additional computation.

For heterogeneous data where $\var(Y|X)$ is not a constant,
conditional quantiles of $Y$ given $X$ can provide additional
information  not captured by conditional mean in describing the
relationship between $Y$ and $X$. Quantile regression [Koenker and
Bassett (1978)] has evolved into a quite  mature area of research
over the past several decades. For high dimensional heterogeneous
data, quantile regression models continue to be useful; see Wang, Wu
and Li (2012) and He, Wang and Hong (2013) for recent contributions.
For the use of quantile based models and methods in genomic data
analysis, see Wang and He (2007, 2008), among others. Motivated by
the usefulness of conditional quantile modeling, we propose to
extend the test for conditional mean independence $H_0'$ and introduce a
test for  conditional quantile independence at a particular quantile
level $\tau\in (0,1)$. It is worth  mentioning that we are not aware
of any tests in the literature developed for conditional quantile
independence in the high dimensional setting.

There is a recent surge of interest in inference for regression
coefficients in high dimensional setting, as motivated by the
emergence of high-dimensional data in various scientific areas such
as medical imaging, atmospheric and climate sciences, and biological
science. The analysis of high-dimensional data in the regression
framework is challenging due to the high dimensionality of
covariates, which can greatly exceed the sample size, and
traditional statistical methods developed for low and
fixed-dimensional covariates may not be applicable to modern large
data sets. Here we briefly mention some related work on testing for
dependence or the significance of regression coefficients in high dimensional regression
setting. Yata and Aoshima (2013) proposed a simultaneous test for
zero correlation between $Y$ and each component of $X$ in a
nonparametric setting; Feng, Zou and Wang (2013) proposed a
rank-based test building on the one in Zhong and Chen (2011) under
the linear model assumption; Testing whether a subset of $X$ is
significant or not after controlling for the effect of the remaining
covariates is addressed in Wang and Cui (2013) and Lan, Wang and
Tsai (2014); Sz\'{e}kely  and Rizzo (2013) studied the limiting
behavior of distance correlation [Sz\'{e}kely, Rizzo and Bakirov
(2007)] based test statistic for independence test when the
dimensions of $Y$ and $X$ both go to infinity while holding sample
size fixed; Goeman et al. (2006) introduced a score test against a
high dimensional alternative in an empirical Bayesian model.

It should be noted that our testing problem is intrinsically different from the marginal feature
screening problem in terms of the formulation, goal and methodology. The marginal screening has been extensively studied in the literature
since Fan and Lv's (2008) seminal work; see Liu, Zhong and Li (2015) for an excellent recent review.
Generally speaking,    marginal screening aims to rank the covariates using sample marginal utility (e.g., sample MDD), and reduce the dimension of covariates to a
moderate scale, and then perform variable selection using LASSO etc. For example, if we use MDD as the marginal utility, as done in Shao and Zhang (2014),
we then screen out the variables that do not contribute to the mean of $Y$. The success of the feature screening procedure depends on whether the subset of kept variables  is able to contain the set of all the important variables, which is referred as the sure screening consistency in theoretical investigations.
For the  testing problem we address in this paper, we aim to detect conditional mean dependence of $Y$ on $X$, and the
formulation and the output of the inference is very different from that delivered by the marginal screening. In
practice, the testing can be conducted prior to model building, which includes
marginal screening as an initial step. In particular, if our null hypothesis holds, then there is no need to
do marginal screening using MDD, as the mean of the response does not depend on covariates. It is worth noting that our test is closely tied to the
adaptive resampling test of McKeague and Qian (2015), who test for the marginal uncorrelatedness of $Y$ with $x_j$, for $j=1,\cdots,p$ in a linear regression framework.
We shall leave a more detailed discussion of the connection and difference of two papers and numerical comparison to the supplement.

 %They pointed out two limitations of the classical $F$ test and proposed a new test statistic that has limiting normal distribution under the null  in high dimensional setting where $p$ can grow much faster than sample size $n$.

The rest of the paper is organized as follows. We propose the
MDD-based conditional mean independence test and study its
asymptotic properties in Section \ref{sec:mean}. We describe a
testing procedure for conditional quantile independence in Section
\ref{sec:quant}.
%An extension to factorial designs is made in Section \ref{sec:factor}.
Section \ref{sec:num} is devoted to numerical studies. %and
%empirical analysis of a microarray dataset on Yorkshire Gilts with factorial designs.
Section \ref{sec:con} concludes. The technical
details, extension to factorial designs, additional discussions and numerical results are gathered in the supplementary material.

\section{Conditional mean dependence testing}\label{sec:mean}
For a scalar response variable $Y$ and a vector predictor
$X=(x_1,\dots,x_{p})^T\in\mathbb{R}^p$, we are interested in testing
the null hypothesis
\begin{equation}\label{eq:h0}
H_0: \E[Y|x_1,\dots,x_p]=\E[Y]\quad \text{almost surely}.
\end{equation}
Notice that under $H_0$, we have $\E[Y|x_j]=\E[Y]$ for all $1\leq
j\leq p$. This observation motivates us to consider the following
hypothesis
\begin{equation}
H_0': \E[Y|x_j]=\E[Y]\quad  \text{almost surely},
\end{equation}
for all $1\leq j\leq p$ versus the alternative that
\begin{equation}
H_a': P(\E[Y|x_j]\neq \E[Y])>0
\end{equation}
for some $1\leq j\leq p.$ We do not work directly on the original hypothesis in (\ref{eq:h0}),
because a MDD-based testing procedure targeting for the global null
$H_0$ generally fails to capture the nonlinear conditional mean
dependence when $p$ is large; see Remark \ref{rm:h0} for more
detailed discussions. In addition, if we regard $\E[Y|x_j]-\E[Y]$ as
the main effect of $x_j$ contributing to the mean of $Y$, then it
seems reasonable to test for the nullity of main effects first,
before proceeding to the higher-order interactions, in a manner
reminiscent of ``Effect Hierarchy Principle" and ``Effect Heredity
Principle" in design of experiment [see pages 172-173 of Wu and
Hamada (2009)]. Testing $H_0'$ versus $H_a'$ is a problem of own
interest from the viewpoint of nonparametric regression modeling,
where additive and separate effects from each $x_j$ usually first
enter into the model before adding interaction terms.

In some cases, $\E[Y|x_j]=0$ almost surely for $1\leq j\leq p$ implies
that $\E[Y|x_1,\dots,x_p]=0$ almost surely.
\begin{example}\label{linear}
{\rm Consider a linear regression model
$$Y=\beta_1x_1+\beta_2x_2+\cdots+\beta_px_p+\epsilon,$$
where $X=(x_1,\dots,x_p)^T$ follows the elliptical distribution $EC_p(0,\Sigma_0,\phi_0)$ with $\phi_0$ being the so-called characteristic generator of $X$ and $\Sigma_0=(\sigma_{0,ij})^{p}_{i,j=1}$ [see e.g. Fang et al. (1990)],
and $\epsilon$ is independent of
$X$ with $\E\epsilon=0.$ Under the elliptical distribution assumption, we have
$$\E[Y|x_j]=\sum^{p}_{k=1}\beta_k\E[x_k|x_j]=\sum^{p}_{k=1}\beta_k\frac{\sigma_{0,jk}}{\sigma_{0,jj}}x_j.$$
Thus $\E[Y|x_j]=0$ almost surely implies that
$\sum^{p}_{k=1}\beta_k\sigma_{0,jk}=0$ for all $1\leq j\leq p.$
Writing this in matrix form, we have $\Sigma_0\beta=\mathbf{0}_{p\times 1}$ for
$\beta=(\beta_1,\dots,\beta_p)^T$, which suggests that $\beta=\mathbf{0}_{p\times 1}$
provided that $\Sigma_0$ is non-singular. In this case, testing $H_0$
is equivalent to testing its marginal version $H_0'$. In particular, the above result holds if $X$ follows a multivariate normal distribution
with mean zero and covariance $\Sigma_0.$}
\end{example}

To characterize the conditional mean (in)dependence, we consider the
MDD proposed in Shao and Zhang (2014). Throughout the following
discussions, we let $\Y$ be a generic univariate random variable and
$\X$ be a generic $q$-dimensional random vector (e.g. $\X=X$ and
$q=p$, or $\X=x_j$ for some $1\leq j\leq p$ and $q=1$). Let
$L(y,y')=(y-y')^2/2$ and $K(x,x')=|x-x'|_q$ with
$|\cdot|_q$ denoting the Euclidean norm of $\mathbb{R}^q$. Under
finite second moment assumptions on $\Y$ and $\X$, the MDD which
characterizes the conditional mean dependence of $\Y$ given $\X$ is
defined as
\begin{equation}
\begin{split}
MDD(\Y|\X)^2=&\E[K(\X,\X')L(\Y,\Y')]+\E[K(\X,\X')]\E[L(\Y,\Y')]
\\&-2\E[K(\X,\X')L(\Y,\Y'')]
\\=&-\E[(\Y-\E \Y)(\Y'-\E \Y')|\X-\X'|_q],
\end{split}
\end{equation}
where $(\X',\Y')$ and $(\X'',\Y'')$ denote independent copies of
$(\X,\Y)$; see Shao and Zhang (2014) and Park et al. (2015). The MDD
is an analogue of distance covariance introduced in Sz\'{e}kely,
Rizzo and Bakirov (2007) in that the same type of weighting function
is used in its equivalent characteristic function-based definition and it completely characterizes
the conditional mean dependence in the sense that $MDD(\Y|\X)=0$ if and only if
$E(\Y|\X)=E(\Y)$ almost surely; see Section 2 of Shao and Zhang (2014). Note that the distance
covariance is used to measure the dependence, whereas MDD is for
conditional mean dependence.

In this paper, we introduce a new testing procedure based on an
unbiased estimator of $MDD(\Y|\X)^2$ in the high dimensional regime.
Compared to the modified $F$-test proposed in Zhong and Chen (2011),
our procedure is model-free and is able to detect nonlinear
conditional mean dependence. Moreover, the proposed method is
quite useful for detecting ``dense'' alternatives, i.e., there
are a large number of (possibly weak) signals, and it can be easily
extended to detect conditional quantile (in)dependence at a given
quantile level (see Section \ref{sec:quant}). It is also worth
mentioning that our theoretical argument is in fact applicable to
testing procedures based on a broad class of dependence measures
including MDD as a special case [see e.g. Sz\'ekely et al. (2007)
and Sejdinovick et al. (2013)].

\subsection{Unbiased estimator}\label{subsec:unbias}
Given $n$ independent observations $(\X_i,\Y_i)^{n}_{i=1}$ from the
joint distribution of $(\X,\Y)$ with
$\X_i=(\X_{i1},\X_{i2},\dots,\X_{iq})^T\in\mathbb{R}^q$ and
$\Y_i\in\mathbb{R}$, an unbiased estimator for $MDD(\Y|\X)^2$ can
be constructed using the idea of $\mathcal{U}$-centering
[Sz\'{e}kely and Rizzo (2014); Park et al. (2015)]. Define
$A=(A_{ij})^{n}_{i,j=1}$ and $B=(B_{ij})^{n}_{i,j=1}$, where
$A_{ij}=|\X_i-\X_j|_q$ and $B_{ij}=|\Y_i-\Y_j|^2/2.$ Following Park et al. (2015), we define the $\mathcal{U}$-centered versions of
$A_{ij}$ and $B_{ij}$ as
\begin{align*}
&\widetilde{A}_{ij}=A_{ij}-\frac{1}{n-2}\sum^{n}_{l=1}A_{il}-\frac{1}{n-2}\sum^{n}_{k=1}A_{kj}+\frac{1}{(n-1)(n-2)}\sum_{k,l=1}^{n}A_{kl},\\
&\widetilde{B}_{ij}=B_{ij}-\frac{1}{n-2}\sum^{n}_{l=1}B_{il}-\frac{1}{n-2}\sum^{n}_{k=1}B_{kj}+\frac{1}{(n-1)(n-2)}\sum_{k,l=1}^{n}B_{kl}.
\end{align*}
A suitable estimator for $MDD(\Y|\X)^2$ is given by
$$MDD_n(\Y|\X)^2=(\widetilde{A}\cdot\widetilde{B}):=\frac{1}{n(n-3)}\sum_{i\neq j}\widetilde{A}_{ij}\widetilde{B}_{ij},$$
where $\widetilde{A}=(\widetilde{A}_{ij})_{i,j=1}^n$ and
$\widetilde{B}=(\widetilde{B}_{ij})_{i,j=1}^{n}$. In Section
1.1 of the supplement, we show that $MDD_n(\Y|\X)^2$ is indeed an
unbiased estimator for $MDD(\Y|\X)^2$ and it has the expression
\begin{align*}
MDD_n(\Y|\X)^2=\frac{1}{\binom{n}{4}}\sum_{i<j<q<r}h(\Z_i,\Z_j,\Z_q,\Z_r),
\end{align*}
where
\begin{align*}
h(\Z_i,\Z_j,\Z_q,\Z_r)=\frac{1}{6}\sum_{s<t,u<v}^{(i,j,q,r)}(A_{st}B_{uv}+A_{st}B_{st})-\frac{1}{12}\sum_{(s,t,u)}^{(i,j,q,r)}A_{st}B_{su}
\end{align*}
with $\Z_i=(\X_i,\Y_i)$, and $\underset{s<t,u<v}{\overset{(i,j,q,r)}{\sum}}$ denotes the
summation over the set $\{(s,t,u,v):s<t,u<v,\text{$(s,t,u,v)$ is a permutation of $(i,j,q,r)$}\}$.
Therefore $MDD_{n}(\Y|\X)^2$ is a $U$-statistic of
order four. This observation plays an important role in subsequent
derivations.

\subsection{The MDD-based test}\label{sec:testing}
Given $n$ observations $(X_i,Y_i)_{i=1}^{n}$ with
$X_i=(x_{i1},\dots,x_{ip})^T$ from the joint distribution of
$(X,Y)$, we consider the MDD-based test statistic defined as
\begin{align}
T_n=\sqrt{\binom{n}{2}}\frac{\sum^{p}_{j=1}MDD_n(Y|x_j)^2}{\hat{\mathcal{S}}},
\end{align}
where $\hat{\mathcal{S}}^2$ is a suitable variance estimator defined
in (\ref{eq:s}) below. To introduce the variance estimator, let $(X',Y')$
with $X'=(x_1',\dots,x_p')^T$ be an independent copy of $(X,Y)$.
Define $U_j(x,x')=\E [K(x,x_j')]+\E [K(x_j,x')]-K(x,x')-\E [K(x_j,x_j')]$
and $V(y,y')=(y-\mu)(y'-\mu)$ with $\mu=\E Y$ and
$x,x',y,y'\in\mathbb{R}$. We further define the infeasible test
statistic
\begin{align*}
\breve{T}_n=&\sqrt{\binom{n}{2}}\frac{\sum^{p}_{j=1}MDD_n(Y|x_j)^2}{\mathcal{S}},
\end{align*}
where $\mathcal{S}$ is given in (\ref{eq:s0}). Because
$MDD_n(Y|x_j)^2$ is a $U$-statistic, by the Hoeffding decomposition
(see Section 1.2 of the supplement), we have under
$H_0'$,
$$MDD_n(Y|x_j)^2=\frac{1}{\binom{n}{2}}\sum_{1\leq k<l\leq n}U_j(x_{kj},x_{lj})V(Y_k,Y_l)+R_{j,n},$$
where $R_{j,n}$ is the remainder term. It thus implies that
\begin{align*}
\breve{T}_n=&\sum^{p}_{j=1}\frac{1}{\sqrt{\binom{n}{2}}\mathcal{S}}\sum_{1\leq
k<l\leq
n}U_j(x_{kj},x_{lj})V(Y_k,Y_l)+\frac{\sqrt{\binom{n}{2}}}{\mathcal{S}}\sum^{p}_{j=1}R_{j,n}
%\\=&\frac{1}{\sqrt{\binom{n}{2}}\mathcal{S}}\sum_{1\leq k<l\leq n}V(Y_k,Y_l)\sum^{p}_{j=1}U_j(x_{kj},x_{lj})+\frac{\sqrt{\binom{n}{2}}}{\mathcal{S}}\sum^{p}_{j=1}R_{j,n}
\\=&\frac{1}{\mathcal{S}}(J_{n,1}+J_{n,2}),
\end{align*}
where
$J_{n,1}=\binom{n}{2}^{-1/2}\sum_{1\leq k<l\leq n}\sum^{p}_{j=1}V(Y_k,Y_l)U_j(x_{kj},x_{lj})$ is the leading term and
$J_{n,2}=\binom{n}{2}^{1/2}\sum^{p}_{j=1}R_{j,n}$ is the remainder term.
Notice that under $H_0'$,
\begin{align}\label{eq:s0}
\var(J_{n,1})=\sum^{p}_{j,j'=1}\E V(Y,Y')^2U_j(x_{j},x_{j}')U_{j'}(x_{j'},x_{j'}').
\end{align}
Moreover, because the contribution from $J_{n,2}$ is asymptotically
negligible (see Theorem \ref{thm1}), we may set
$\mathcal{S}^2=\var(J_{n,1}).$ Based on the discussions in Section
\ref{subsec:unbias}, we propose the following variance estimator for $\mathcal{S}^2$,
\begin{align}\label{eq:s}
\hat{\mathcal{S}}^2=\frac{2}{n(n-1)c_n}\sum_{1\leq k<l\leq n}\sum^{p}_{j,j'=1}\widetilde{A}_{kl}(j)\widetilde{A}_{kl}(j')\widetilde{B}_{kl}^2,
\end{align}
where $\widetilde{A}_{kl}(j)$ is the $\mathcal{U}$-centered version
of $A_{kl}(j)=|x_{kj}-x_{lj}|$, and
$$c_n=\frac{(n-3)^4}{(n-1)^4}+\frac{2(n-3)^4}{(n-1)^4(n-2)^3}+\frac{2(n-3)}{(n-1)^4(n-2)^3}\approx \frac{(n-3)^4}{(n-1)^4},$$
is a finite sample adjustment factor to reduce the bias of $\hat{\mathcal{S}}^2$. See more details in the supplement.

\begin{remark}
{\rm We remark that $\hat{\mathcal{S}}^2$ is a biased estimator for
$\mathcal{S}^2$. To construct an unbiased estimator for
$\mathcal{S}^2$, we let
$a(x_1,x_2,x_3,x_4)=|x_1-x_2|-|x_1-x_3|-|x_2-x_4|+|x_3-x_4|$ and
$b(y_1,y_2,y_3,y_4)=(y_1-y_3)(y_2-y_4)$ for $x_i,y_i\in\mathbb{R}$.
Define the estimator
\begin{align*}
\breve{\mathcal{S}}^2=&\frac{1}{\binom{n}{10}}\sum_{(j_1,\dots,j_{10})}\sum_{j,j'=1}^p
a(x_{j_1,j},x_{j_2,j},x_{j_3,j},x_{j_4,j})a(x_{j_1,j'},x_{j_2,j'},x_{j_5,j'},x_{j_6,j'})
\\&\times
b(Y_{j_1},Y_{j_2},Y_{j_7},Y_{j_8})b(Y_{j_1},Y_{j_2},Y_{j_9},Y_{j_{10}}),
\end{align*}
where the summation $\sum_{(j_1,\dots,j_{10})}$ is over all
combination $(j_1,\dots,j_{10})$ chosen without replacement from
$\{1,2,\dots,n\}$. It is straightforward to verify that
$\E[\breve{\mathcal{S}}^2]=\mathcal{S}^2$. However, the calculation of $\breve{\mathcal{S}}^2$ is much more expensive and it seems less convenient to use
in the high-dimensional case.
}
\end{remark}

\begin{remark}\label{rm:h0}
{\rm An alternative test statistic can be constructed based on the
unbiased estimator of $MDD(Y|X)^2$. Specifically, we consider
$$\widetilde{T}_n=\sqrt{\binom{n}{2}}\frac{MDD_n(Y|X)^2}{\hat{\eta}},$$
with $\hat{\eta}^2$ being a suitable estimator for $\E
\bar{U}^2(X,X')V^2(Y,Y')$, where $\bar{U}(x,x')=\E [K(x,X')]+\E[
K(X,x')]-K(x,x')-\E [K(X,X')]$. Note that $\widetilde{T}_n$ targets
directly for the global hypothesis in (\ref{eq:h0}). Using Taylor
expansion and similar arguments in Sz\'{e}kely  and Rizzo (2013), it
is expected that
$\bar{U}(X,X')=\frac{<X-\mu_X,\mu_X-X'>}{\sqrt{\tau}}+R^*$ as
$p\rightarrow +\infty$, where $\tau=\E|X-X'|^2_p$, $\mu_X=\E X$ and
$R^*$ is the remainder term. We have up to a smaller order term
\begin{align*}
MDD(Y|X)^2=&-E[V(Y,Y')\bar{U}(X,X')]
\\ \approx& \frac{1}{\sqrt{\tau}}E[<(Y-\mu)(X-\mu_X),(Y'-\mu)(X'-\mu_X)>]
\\=&\frac{1}{\sqrt{\tau}}<\E(Y-\mu)(X-\mu_X),\E(Y'-\mu)(X'-\mu_X)>
\\=&\frac{1}{\sqrt{\tau}}\sum^{p}_{j=1}\cov^2(Y,x_j)
\end{align*}
with $\mu=\E Y$ and $<a,b>=a^Tb$ for $a,b\in\mathbb{R}^p.$
Thus for large $p$, $\widetilde{T}_n$ can be viewed
as a combination of the pairwise covariances. Because $\cov(Y,x_j)$
only captures the linear dependence, the proposed test $T_n$ enjoys
the advantage over $\widetilde{T}_n$ in the sense of detecting
certain degree of component-wise nonlinear dependence. This idea is also related to additive modeling, where the
effects of the covariates cumulate in an additive way. }
\end{remark}

%In fact, the test proposed in Chen and Qin (2010) also belongs to this class of
%tests which aggregates the marginal/component-wise metrics through summation.

Below we study the asymptotic properties of the proposed MDD-based
test. Let $Z'=(X',Y')$, $Z''=(X'',Y'')$ and $Z'''=(X''',Y''')$ be
independent copies of $Z=(X,Y)$. Define
$\widetilde{U}(X,X')=\sum^p_{j=1}U_j(x_j,x_j')$ and
$H(Z,Z')=V(Y,Y')\widetilde{U}(X,X')$. Further define
$G(Z,Z')=E[H(Z,Z'')H(Z',Z'')|(Z,Z')]$. We present the following
result regarding the asymptotic behavior of $\breve{T}_n$ under
$H_0'$.
\begin{theorem}\label{thm1}
Under the assumption that
\begin{align}
&\frac{\E [G(Z,Z')^2]}{\{\E [H(Z,Z')^2]\}^2}\rightarrow 0, \quad
\frac{\E [H(Z,Z')^4]/n+\E [H(Z,Z'')^2H(Z',Z'')^2]}{n\{\E
[H(Z,Z')^2]\}^2}\rightarrow 0, \label{eq:con1}
\end{align}
and the null hypothesis $H_0'$, we have
\begin{align}
\frac{J_{n,1}}{\mathcal{S}}=\frac{1}{\sqrt{\binom{n}{2}}\mathcal{S}}\sum_{1\leq
i<j\leq n}H(Z_i,Z_j)\rightarrow^d N(0,1).
\end{align}
Moreover assuming that
\begin{align}
&\frac{\E [\widetilde{U}(X,X'')^2V(Y,Y')^2]}{\mathcal{S}^2}=o(n),
\label{eq-r1}
%\\&\frac{\E \widetilde{U}(X,X'')\widetilde{U}(X',X'')V(Y,Y')^2}{\mathcal{S}^2}=O(1), \label{eq-r2}
\\&\frac{\var(Y)^2\sum_{j,j'=1}^pdcov(x_j,x_{j'})^2}{\mathcal{S}^2}=o(n^2), \label{eq-r3}
\end{align}
where $$dcov(x_j,x_{j'})^2=\E U_j(x_j,x_j')U_{j'}(x_{j'},x_{j'}')$$
denotes the (squared) distance covariance between $x_j$ and $x_{j'}$ [see the
definition in Sz\'{e}kely et al. (2007)]. Then we have $J_{n,2}/\mathcal{S}=o_p(1)$ and hence
\begin{align}\label{eq:u-approx2}
\breve{T}_n\rightarrow^d N(0,1).
\end{align}
\end{theorem}

The following theorem shows that $\hat{\mathcal{S}}^2$ is ratio-consistent under the null.
\begin{theorem}\label{thm11}
Under the conditions in Theorem \ref{thm11}, we have $\hat{\mathcal{S}}^2/\mathcal{S}^2\rightarrow^p 1.$ As a consequence, $T_n\rightarrow^d N(0,1).$
\end{theorem}

Condition (\ref{eq:con1}) ensures the asymptotic normality of the
leading term $J_{n,1}$ as well as the ratio-consistency for $\hat{\mathcal{S}}^2$,
while Conditions (\ref{eq-r1})-(\ref{eq-r3})
guarantee that the remainder $J_{n,2}$ is asymptotically negligible.
To gain more insight on Condition (\ref{eq:con1}), we define the
operator
$$A(g)(z)=\E[H(Z,z)g(Z)],$$
where $g:\mathbb{R}^{p+1}\rightarrow \mathbb{R}$ is a
measurable function. Let $\{\lambda_j\}^{+\infty}_{j=1}$ be the set
of eigenvalues for $A$. To obtain a normal limit for the leading
term, a typical condition [see Hall (1984)] is given by
\begin{align}\label{eq-eigen}
\frac{(\sum^{+\infty}_{j=1}\lambda_j^t)^{2/t}}{\sum^{+\infty}_{j=1}\lambda_j^2}\rightarrow
0,
\end{align}
for some $t>2.$ Note that $\E
G(Z,Z')^2=\sum^{+\infty}_{j=1}\lambda_j^4$ and $\E
H(Z,Z')^2=\sum^{+\infty}_{j=1}\lambda_j^2.$ Thus the condition $\E
G(Z,Z')^2/ \{\E H(Z,Z')^2\}^2\rightarrow 0$ is equivalent to
(\ref{eq-eigen}) with $t=4$.

%\begin{remark}\label{remark-cq}
%{\rm \textbf{Maybe remove this remark.} We connect our assumptions with those in the literature. With
%some abuse of notation, let $X_1,\dots,X_n$ be $n$ random samples from a
%$p$-dimensional distribution. We consider the degenerate
%$U$-statistic $Q_n:=\sum_{1\leq i\neq j\leq n}X_i^TX_j$ for testing
%the null hypothesis $\E[X_i]=0.$ In this case, $H$ is defined as the
%kernel associated with the degenerate $U$-statistic and it has the
%form $H(x,x')=x^Tx'$ for $x,x'\in\mathbb{R}^p$. Thus the operator
%$A$ is given by $A(g)(x)=\E[x^TXg(X)].$ Suppose $h(x)$ is the
%eigenfunction of $A$. Then we have $A(h)(x)=x^T\E[Xh(X)]=\lambda
%h(x).$ Thus $h(x)$ must be linear in $x$, i.e., $h(x)=x^T\gamma$ for
%$\gamma\in\mathbb{R}^p$. Under the null that $\E[X_i]=0$, we have
%$$A(h)(x)=x^T\E[XX^T]\gamma=x^T\widetilde{\Sigma}\gamma=\lambda x^T\gamma,$$
%with $\widetilde{\Sigma}$ being the covariance matrix of $X$. Since
%$x^T\widetilde{\Sigma}\gamma=\lambda x^T\gamma$ holds for every $x$,
%$\gamma$ must be the eigenvector of $\widetilde{\Sigma}$ and
%$\lambda$ is the corresponding eigenvalue. In this case,
%$\{\lambda_1,\dots,\lambda_p\}$ is the set of eigenvalues for
%$\widetilde{\Sigma}.$ When $t=4$, (\ref{eq-eigen}) reduces to
%\begin{align*}
%\frac{(\sum^{p}_{j=1}\lambda_j^4)^{1/2}}{\sum^{p}_{j=1}\lambda_j^2}=\frac{Tr^{1/2}(\widetilde{\Sigma}^4)}{Tr(\widetilde{\Sigma}^2)}\rightarrow0,
%\end{align*}
%with $Tr(\cdot)$ denoting the trace operator, which is essentially
%Condition (3.7) in Chen and Qin (2010). }
%\end{remark}

\subsection{Further discussions on the conditions}\label{sec:further}
Because MDD has an interpretation based on characteristic function,
we provide more discussions about Condition (\ref{eq-eigen}) based
on this interpretation. Let $\I=\sqrt{-1}$ and $w(t)=1/(\pi t^2)$, and denote by $\bar{a}$ the conjugate of a complex number $a$.
Define $G_j(u_j,y;t)=(f_{x_j}(t)-e^{\I tu_j})(\mu-y)$ for $u_j, y, t
\in \mathbb{R}$, where $f_{x_j}$ denotes the characteristic function for $x_j$.
By Lemma 1 of Sz\'{e}kely et al. (2007), it can be
shown that
\begin{align*}
H(z,z')=\int_{\mathbb{R}}\sum^{p}_{j=1}G_j(u_j,y;t)\overline{G_j(u_j',y';t)}w(t)dt,
\end{align*}
where $z=(u_1,\dots,u_p,y)^T$ and $z'=(u_1',\dots,u_p',y')^T$. Let
$L^2(w)$ be the space of functions such that $\int_{\mathbb{R}}
f(t)^2w(t)dt<\infty$ for $f\in L^2(w)$. For $g\in L^2(w)$, let $\varphi_{ji}(g)(t)=\int_{\mathbb{R}}
g(t')\E[G_j(x_j,Y;t)\overline{G_i(x_i,Y;t')}]w(t')dt'$.
%When $X=(x_1,\dots,x_p)^T$ and $(Y-\mu)^2$ are independent, we have
%\begin{equation*}
%\E[G_j(x_j,Y;t)\overline{G_i(x_i,Y;t')}]=\var(Y)(f_{x_j,x_i}(t,-t')-f_{x_j}(t)f_{x_i}(-t')),
%\end{equation*}
%in which case $\varphi_{ji}(g)(t)=\var(Y)\int_{\mathbb{R}}
%g(t')(f_{x_j,x_i}(t,-t')-f_{x_j}(t)f_{x_i}(-t'))w(t')dt'.$
For $g=(g_1,\dots,g_p)^T$ with $g_j\in L^2(w)$, define the operator
$$\Phi(g)=\left(\sum^{p}_{i=1}\varphi_{1i}(g_i),\dots,\sum^{p}_{i=1}\varphi_{pi}(g_i)\right).$$
It is shown in Section 1.4 of the supplement that $\lambda_j$ is the
eigenvalue associated with $\Phi.$ If we define $Tr$ as the nuclear
norm for the Hermitian operator $\Phi.$ Then with $t=4$, the
condition required is that
\begin{align*}
\frac{Tr^{1/2}(\Phi^4)}{Tr(\Phi^2)}\rightarrow 0,
\end{align*}
which can be viewed as a generalization of the condition in the
linear case [see Chen and Qin (2010); Zhong and Chen (2011)] to the nonlinear situation.

%e.g., $\frac{Tr^{1/2}(\widetilde{\Sigma}^4)}{Tr(\widetilde{\Sigma}^2)}\rightarrow 0$ (see Remark \ref{remark-cq})

Below we show that the assumptions in Theorem \ref{thm1} can be made more explicit in the following two cases: (i) banded dependence structure; (ii) Gaussian design.
\begin{assumption}\label{ass3}
Assume that
$$0<c\leq \text{var}(Y|X)\leq
\E[(Y-\E[Y|X])^4|X]^{1/2}\leq C<+\infty,$$
almost surely for some constants $c$ and
$C.$
\end{assumption}
Suppose $Y=g(X)+\sigma(X)\epsilon$, where $E[\epsilon|X]=0$,
$\E[\epsilon^2|X]\geq c_1>0$ and $\E[\epsilon^4|X]<C_1<\infty$. If $0<c<\sigma(X)<C<\infty,$ it can be verified that Assumption \ref{ass3} holds.
Note that Assumption~\ref{ass3} can be relaxed with additional conditions but at the expense of lengthy details.

\begin{proposition}
Suppose for an unknown permutation $\pi:\{1,2,\dots,p\}\rightarrow \{1,2,\dots,p\}$:
$$x_{\pi(i)}=g_{\pi(i)}(\epsilon_{i},\epsilon_{i+1},\dots,\epsilon_{i+L}),$$
for $1\leq i\leq p$ and $L\geq 0$, where $\{\epsilon_i\}$ are independent
random variables, and $g_i$ is measurable function such that $x_{i}$ is well defined. Suppose
\begin{align}\label{eq-band}
\frac{p(L+1)^3\max\left\{\max_{1\leq j\leq
p}n^{-1}\text{var}(x_j)^2,\max_{1\leq j\leq p}(\E[|x_j-\mu_j|])^4\right\}
}{(\sum_{|j-k|\leq L}dcov(x_j,x_{k})^2)^2}\rightarrow 0,
\end{align}
and Assumption \ref{ass3} holds. Then the conditions in Theorem \ref{thm1} are satisfied.
In particular, if $L=o(p^{1/3})$ and $\max_{1\leq j\leq
p}\text{var}(x_j)/\min_{1\leq j\leq p}dcov(x_j,x_j)^2$ is bounded from above,
(\ref{eq-band}) is fulfilled.
\end{proposition}

\begin{proposition}
Suppose
$X\sim N(0,\Sigma)$ with $\Sigma=(\sigma_{jk})_{j,k=1}^p$ and
$\sigma_{ii}=1$. Assume that $\max_{j\neq k}|\sigma_{jk}|\leq c<1$ and
\begin{align}\label{ass:normal}
\frac{\sum_{j,k,l,m=1}^{p}|\sigma_{jk}\sigma_{kl}\sigma_{lm}\sigma_{mj}|}{\text{Tr}^2(\Sigma^2)}\rightarrow
0,\quad \frac{\sum_{j=1}^{p}(\sum_{k=1}^{p}\sigma_{jk}^2)^3}{n\text{Tr}^2(\Sigma^2)}\rightarrow 0,
\end{align}
Then under Assumption \ref{ass3}, the conditions in Theorem \ref{thm1} are satisfied.
\end{proposition}

From the above two propositions, we can see that the asymptotic normality of our test statistic under the null hypothesis
 can be derived without any explicit constraints on the growth rate of $p$ as a function of $n$.
 In certain cases (e.g., in the Gaussian design with $\Sigma=I_p$), $p$ can grow to infinity freely as $n\rightarrow\infty$.

\subsection{Wild bootstrap}\label{sec:boot}
To improve the finite sample performance when sample size is small, we further propose a wild bootstrap approach to approximate
the null distribution of the MDD-based test statistic as a useful alternative to normal approximation.
Recall from Section \ref{subsec:unbias} that $MDD_n(Y|x_j)^2\\=\frac{1}{n(n-3)}\sum_{k\neq l}\widetilde{A}_{kl}(j)\widetilde{B}_{kl}$, we
then need to perturb the sample $MDD$s using a symmetric distributed random variable with mean 0 and variance 1.
Specifically, generate $\{e_{i}\}_{i=1}^n$ from standard normal distribution.
Then the bootstrap $MDD_n$ is defined as
\[
MDD^{*}_n(Y|x_j)^2=\frac{1}{n(n-1)}\sum_{k\neq l}\widetilde{A}_{kl}(j)\widetilde{B}_{kl}e_{k}e_{l}.
\]
%\textbf{For technical reason, it is convenient to use $n-1$ instead of $n-3$ in the denominator.}
The studentized bootstrap statistic is given by
$$T_n^{*}=\sqrt{\binom{n}{2}}\frac{\sum^{p}_{j=1}MDD^{*}_n(Y|x_j)^2 }{\hat S^{*}},
$$
where $$\hat S^{*2}=\frac{1}{\binom{n}{2}}\sum^{p}_{j,j'=1}\sum_{1\leq k<l\leq
n}\widetilde{A}_{kl}(j)\widetilde{A}_{kl}(j')\widetilde{B}_{kl}^2e_{k}^2e_{l}^2.
$$
The form of $\hat S^{*2}$ can be motivated from the same argument
used in (\ref{eq:s0}) by deriving the leading term in  the
unconditional variance of $\sum^{p}_{j=1}MDD^{*}_n(Y|x_j)^2$. Notice
that conditional on the sample, $\hat S^{*2}$ is an unbiased
estimator for
$var^*(\sqrt{\binom{n}{2}}\sum^{p}_{j=1}MDD^{*}_n(Y|x_j)^2)$, where
$var^*$ stands for the variance conditional on the sample. Repeat
the above procedure $B$ times, and denote by $T_n^{*(b)}$ the values
of the bootstrap statistics for $1\leq b\leq B.$ The corresponding
bootstrap p-value is given by $\frac{1}{B}\sum_{b=1}^{B} {\bf
1}(T_n^{*(b)}\ge T_n)$. We reject the null hypothesis when the
p-value is less than the nominal level $\alpha$. To study the
theoretical property of the bootstrap procedure, we first introduce
a notion of bootstrap consistency, also see Li et al. (2003).
\begin{definition}
Let $\zeta_n$ denote a statistic that depends on the random samples $\{Z_i\}^{n}_{i=1}$. We say that $(\zeta_n|Z_1,Z_2,\dots)$ converges to $(\zeta|Z_1,Z_2,\dots)$
in distribution in probability if for any subsequence $\zeta_{n_k}$ there is a further subsequence $\zeta_{n_{k_j}}$ such that $(\zeta_{n_{k_j}}|Z_1,Z_2,\dots)$ converges to $(\zeta|Z_1,Z_2,\dots)$ in distribution for almost every sequence $\{Z_1,Z_2,\dots\}$.
\end{definition}

The following theorem establishes the bootstrap consistency.
\begin{theorem}\label{thm-boot}
Suppose the conditions in Theorem \ref{thm1} hold. Then we have
$$\sqrt{\binom{n}{2}}\frac{\sum^{p}_{j=1}MDD^{*}_n(Y|x_j)^2}{\hat{S}^*}\rightarrow^{\mathcal{D}^*}
N(0,1),$$
where $\rightarrow^{\mathcal{D}^*}$ denotes convergence in distribution in probability.
\end{theorem}

\subsection{Asymptotic analysis under local alternatives}
Below we turn to the asymptotic analysis under local alternatives.
Define $\widetilde{\mathcal{L}}(x,y)=\E \widetilde{U}(x,\X)V(y,\Y)$
and
$H^*(Z,Z')=H(Z,Z')-\widetilde{\mathcal{L}}(X,Y)-\widetilde{\mathcal{L}}(X',Y')+\E[\widetilde{U}(X,X')V(Y,Y')].$
Note that $\E [H^*(z,Z')]=\E [H^*(Z,z')]=0.$ Thus $H^*$ can be viewed as
a demeaned version of $H$ under alternatives. With some abuse of notation, denote
$$G^*(Z,Z')=E[H^*(Z,Z'')H^*(Z',Z'')|(Z,Z')]$$
and let $\mathcal{S}^2=\var(H^*(Z,Z'))$.
\begin{theorem}\label{thm2}
Under the assumptions that
$\var(\widetilde{\mathcal{L}}(X,Y))=o(n^{-1}\mathcal{S}^2)$ and
\begin{align*}
&\frac{\E [G^*(Z,Z')^2]}{\{\E [H^*(Z,Z')^2]\}^2}\rightarrow 0, \quad
\frac{\E [H^*(Z,Z')^4]/n+\E [H^*(Z,Z'')^2H^*(Z',Z'')^2]}{n\{\E
[H^*(Z,Z')^2]\}^2}\rightarrow 0,
\end{align*}
We have
$$\frac{1}{\sqrt{\binom{n}{2}}\mathcal{S}}\sum_{1\leq k<l\leq
n}\left(H(Z_k,Z_l)-\sum^{p}_{j=1}MDD(Y|x_j)^2\right)\rightarrow^d
N(0,1).$$
Moreover assuming
(\ref{eq-r1}), (\ref{eq-r3}), and that
$\var(\widetilde{\mathcal{L}}(X,Y'))=o(\mathcal{S}^2)$ with $\mathcal{S}^2=\var(H^*(Z,Z'))$, we have
$$\sqrt{\binom{n}{2}}\frac{\sum^{p}_{j=1}\{MDD_n(Y|x_j)^2-MDD(Y|x_j)^2\}}{\mathcal{S}}\rightarrow^d N(0,1).$$
\end{theorem}

\begin{remark}
{\rm Under $H_0'$, we have $\widetilde{\mathcal{L}}(x,y)=0$ for any
$x$ and $y$, and $H^*(\cdot,\cdot)=H(\cdot,\cdot)$. Hence Theorem
\ref{thm2} provides a natural generalization of the results in
Theorem \ref{thm1} to the case of local alternatives. }
\end{remark}

\begin{remark}
{\rm The conditions in Theorem \ref{thm2} also have the characteristic function interpretation. In particular, the results in Section 2.3 also hold for local alternatives when $G_j(u_j,y;t)$ is replaced by
its demeaned version $G_j(u_j,y;t)-\E G_j(x_j,Y;t).$
}
\end{remark}

\begin{remark}
{\rm In Section 1.7 of the supplement, we show that the two conditions
\begin{align*}
&\var(\widetilde{\mathcal{L}}(X,Y))=o(n^{-1}\mathcal{S}^2), \quad \var(\widetilde{\mathcal{L}}(X,Y'))=o(\mathcal{S}^2),
\end{align*}
have similar interpretation as equation (4.2) in Zhong and Chen (2011). More precisely, let $h(t)=(h_1(t),\dots,h_p(t))$ with $h_j(t)=\E(f_{x_j}(t)-e^{\I tx_j})(\mu-Y)$ for $1\leq j\leq p$.
The above two conditions essentially impose constraints on the distance between $h(t)$ and $\mathbf{0}_{p\times 1}$ through a metric in some Hilbert space. The characterization of the local alternative model is somewhat abstract but this is sensible because MDD targets at very broad alternatives, i.e., arbitrary type of conditional mean dependence. In the special case that $(X,Y)$ are jointly Gaussian, it is in fact possible to make the conditions in Theorem~\ref{thm2} more explicit, but the details are fairly tedious and complicated, so are omitted due to space limitation.
%somewhat abstract but this is expected because MDD targets at very broad alternatives (arbitrary type of conditional mean dependence).
}
\end{remark}

\begin{remark}\label{rm:power}
{\rm When $\hat{\mathcal{S}}$ is close to $\mathcal{S}$ under the
local alternative, Theorem \ref{thm2} suggests that the power
function of $T_n$ is given by
\begin{align*}
\varphi\left(-c_{\alpha}+\sqrt{\binom{n}{2}}\sum^{p}_{j=1}MDD(Y|x_j)^2/\mathcal{S}\right),
\end{align*}
where $\varphi$ is the distribution function of $N(0,1)$ and
$c_{\alpha}$ is the $100(1-\alpha)\%$ quantile of $N(0,1)$. }
\end{remark}

Below we compare the asymptotic power of the MDD-based test with that of
Zhong and Chen (2011)'s test. Consider a linear regression model
$$Y=\beta_1x_1+\beta_2x_2+\cdots+\beta_px_p+\epsilon,$$
where $X=(x_1,\dots,x_p)^T\sim N(0,\breve{\Sigma})$ with
$\breve{\Sigma}=(\breve{\sigma}_{ij})^{p}_{i,j=1}$, and $\epsilon$ is independent of
$X$ with $\E\epsilon=0$ and $\var\epsilon=\sigma^2.$ By Example \ref{linear}, testing $H_0'$ is
equivalent to testing the global null $H_0$ provided that $\breve{\Sigma}$ is invertible. The power function of
Zhong and Chen (2011)'s test is given by
\begin{align*}
\varphi\left(-c_{\alpha}+\frac{n\beta^T\breve{\Sigma}^2\beta}{\sqrt{2\text{tr}(\breve{\Sigma}^2)}\sigma^2}\right).
\end{align*}
By Remark \ref{rm:power}, the asymptotic relative efficiency (Pitman efficiency) of the
two tests is determined by the ratio,
\begin{align*}
\mathcal{R}=&\frac{\sum^{p}_{j=1}MDD(Y|x_j)^2}{\mathcal{S}}\frac{\sqrt{\text{tr}(\breve{\Sigma}^2)}\sigma^2}{\beta^T\breve{\Sigma}^2\beta}.
\end{align*}
Under the local alternative considered in Theorem \ref{thm2}, we can
show that
\begin{align*}
&\sum^{p}_{j=1}MDD(Y|x_j)^2=-\sum^{p}_{j=1}\beta^T\E[XX^{'T}|x_j-x_j'|]\beta,\\
&\mathcal{S}^2=\sigma^4\sum^{p}_{j,j'=1}dcov(x_j,x_{j'})^2,
\end{align*}
which implies that
\begin{align*}
\mathcal{R}=\frac{-\sum^{p}_{j=1}\beta^T\E[XX^{'T}|x_j-x_j'|]\beta}{\sqrt{\sum^{p}_{j,j'=1}dcov(x_j,x_{j'})^2}}\frac{\sqrt{\text{tr}(\breve{\Sigma}^2)}}{\beta^T\breve{\Sigma}^2\beta}.
\end{align*}
When $\breve{\Sigma}=I_p$, direct calculation shows that
\begin{align*}
\mathcal{R}=\frac{MDD(x_1|x_1)^2}{dcov(x_1,x_{1})}\approx
0.891,
\end{align*}
where we have used Theorem 1 in Shao and Zhang (2014). It therefore
suggests that the MDD-based test is asymptotically less powerful
than Zhong and Chen (2011)'s test but the power loss is fairly moderate.
This is consistent with our
statistical intuition in the sense that Zhong and Chen (2011)'s test
is developed under the linear model assumption while the MDD-based
test is model-free and thus can be less powerful in the linear case.

\begin{remark}
{\rm We mainly focus on the degenerate case in the above
discussions. When the test statistic is non-degenerate under the
alternative, we have the following decomposition [see Serfling
(1980)],
$$\sum^{p}_{k=1}\{MDD_n(Y|x_k)^2-MDD(Y|x_k)^2\}=\frac{2}{n}\sum_{i=1}^n(\widetilde{\mathcal{L}}(X_i,Y_i)-\E \widetilde{\mathcal{L}}(X_i,Y_i))+\widetilde{\mathcal{R}}_{n},$$
with $\widetilde{\mathcal{R}}_{n}$ being the remainder term. In this
case, the asymptotic normality is still attainable under suitable
assumptions. }
\end{remark}

\section{Conditional quantile dependence testing}\label{sec:quant}
In this section, we consider the problem of testing the conditional
quantile dependence in high dimension. Let $Q_{\tau}(Y)$ and
$Q_{\tau}(Y|\cdot)$ be the unconditional and conditional quantiles
of $Y$ at the $\tau$th quantile level respectively. For
$\tau\in(0,1)$, we are interested in testing
$$H_{0,\tau}: Q_{\tau}(Y|x_j)=Q_{\tau}(Y)\quad \text{almost surely},$$
for all $1\leq j\leq p$ versus the alternative that
$$H_{a,\tau}: P(Q_{\tau}(Y|x_j)\neq Q_{\tau}(Y))>0$$
for some $1\leq j\leq p.$ To this end, we define $W_j=\tau-\mathbf{1}\{Y_j\leq Q_{\tau}(Y)\}$ and
$\W_j=\tau-\mathbf{1}\{Y_j\leq \hat{Q}_{\tau}\}$, where
$\hat{Q}_{\tau}$ is the $\tau$th sample quantile of $Y$. Define the
estimator $MDD_n(\W|x_j)^2$ for $MDD(W|x_j)^2$ based on the sample
$(\W_i,X_i)^{n}_{i=1}$. Let
\begin{align*}
\hat{\sigma}_{Q,j,j'}=\frac{1}{\binom{n}{2}}\sum_{1\leq k<l\leq
n}\widetilde{A}_{kl}(j)\widetilde{A}_{kl}(j')(\hat{B}_{Q,kl}^*)^2,
\end{align*}
with $\hat{B}^*_{Q,kl}$ being the $\mathcal{U}$-centered version of
$\hat{B}_{Q,kl}=|\hat{W}_{k}-\hat{W}_{l}|^2/2$. Define
$\hat{\mathcal{S}}^2_Q=\sum^{p}_{j,j'=1}\hat{\sigma}_{Q,j,j'}$. We
consider the test statistic
\begin{align*}
T_{Q,n}=\sqrt{\binom{n}{2}}\frac{\sum^{p}_{j=1}MDD_n(\W|x_j)^2}{\hat{\mathcal{S}}_Q},
\end{align*}
and its infeasible version
\begin{align*}
\breve{T}_{Q,n}=\sqrt{\binom{n}{2}}\frac{\sum^{p}_{j=1}MDD_n(\W|x_j)^2}{\mathcal{S}_Q},
\end{align*}
where $\mathcal{S}_Q^2=\sum^{p}_{j,j'=1}\E
W^2(W')^2U_j(x_{j},x_{j}')U_{j'}(x_{j'},x_{j'}')$. To facilitate the
derivation, we make the following assumptions.

\begin{assumption}\label{ass1}
The cumulative distribution function of the continuous response
variable $Y$, $F_Y$ is continuously differentiable in a small
neighborhood of $Q_{\tau}(Y)$, say $[Q_{\tau}(Y)-\delta_0,
Q_{\tau}(Y)+\delta_0]$ with $\delta_0>0$. Let
$G_{1}(\delta_0)=\inf_{y\in [Q_{\tau}(Y)-\delta_0,
Q_{\tau}(Y)+\delta_0]}f_Y(y)$ and $G_{2}(\delta_0)=\sup_{y\in
[Q_{\tau}(Y)-\delta_0, Q_{\tau}(Y)+\delta_0]}f_Y(y)$, where $f_Y$ is
the density function of $Y$. Assume that $0<G_1(\delta_0)\leq
G_2(\delta_0)<\infty$.
\end{assumption}

\begin{assumption}\label{ass2}
There exists $\delta>0$ such that the collection of random variables $\{\mathbf{1}\{Y\leq Q_{\tau}(Y)+a\},-\delta\leq a\leq \delta\}$ is independent of $X$.
\end{assumption}
Assumption \ref{ass1} was first introduced in Shao and Zhang (2014) and is quite mild.
We note that the conditional quantile independence at $\tau$th quantile level implies that $\mathbf{1}\{Y<Q_{\tau}(Y)\}$ is independent of $X$ as $\mathbf{1}\{Y<Q_{\tau}(Y)\}$ is a Bernoulli random variable.
Thus Assumption \ref{ass2} is quite a bit stronger than the independence between $X$ and $\mathbf{1}\{Y<Q_{\tau}(Y)\}$ and it can be interpreted as a local quantile independence assumption.
We provide some examples in Section~\ref{sec:num} where the local quantile independence assumption may or may not hold at some quantile level; see Examples \ref{ex53} and \ref{ex54}. Next, we
define $\widetilde{T}_{Q,n}$
by replacing $\W$ with $W$ in the definition of $\breve{T}_{Q,n}$.

%The local quantile independence assumption reflects the price we pay for not knowing the true quantile level $Q_{\tau}(Y)$.
%\begin{example}
%{\rm Consider the model $Y=\mathcal{G}(X)\epsilon_1+\epsilon_2$, where $\mathcal{G}(x)\geq a_0$ and $\epsilon_2\geq -a_0$ with $a_0>0$, and $\epsilon_1$ is a Bernoulli random variable with success probability $\gamma_0$.
%Suppose $X$, $\epsilon_1$ and $\epsilon_2$ are mutually independent. Then for $x<0,$ we have
%$$P(Y\leq x)=(1-\gamma_0)P(\epsilon_2\leq x)+\gamma_0P(\mathcal{G}(X)+\epsilon_2\leq x)=(1-\gamma_0)P(\epsilon_2\leq x).$$
%Therefore $Y$ satisfies the local quantile independence assumption for $Q_{\tau}(Y)<0.$}
%\end{example}

\begin{proposition}\label{prop-quan1} Suppose
Assumptions \ref{ass1}-\ref{ass2} hold. We have
$$\left|\widetilde{T}_{Q,n}-\breve{T}_{Q,n}\right|\rightarrow^p 0.$$
\end{proposition}
The above result suggests that the effect by replacing $W$ with $\W$
is asymptotically negligible. The following main result
follows immediately from Proposition \ref{prop-quan1}.

%The key here is to derive an expansion for $$\W_k\W_l-W_kW_l,$$
%based on the Bahadur representation,
%$$\hat{Q}_{\tau}=Q_{\tau}+\frac{\tau-F_n(Q_{\tau})}{f(Q_{\tau})}+R_n.$$
\begin{theorem}
Suppose the assumptions in Theorem \ref{thm1} hold with $Y$ and
$\mathcal{S}$ being replaced by $W$ and $\mathcal{S}_Q$
respectively. Then we have $\widetilde{T}_{Q,n}\rightarrow^d
N(0,1)$. Moreover under Assumptions \ref{ass1}-\ref{ass2}, we have
$\breve{T}_{Q,n}\rightarrow^d N(0,1)$.
\end{theorem}

\section{Simulation studies}\label{sec:num}

In this section, we present simulation results for conditional mean
dependence tests in Section \ref{sec:sim}; conditional quantile
dependence tests in Section \ref{sec:simquan}. We report the empirical sizes of the wild bootstrap approach
in Section \ref{sec:simboot}.

\subsection{Conditional mean independence}\label{sec:sim}

In this subsection, we conduct several simulations to assess the
finite sample performance of the proposed MDD-based test for
conditional mean independence (mdd, hereafter). For comparison, we
also implement the modified $F$-test proposed by Zhong and Chen
(2011) (ZC test, hereafter). Note that the ZC test is designed for
testing the regression coefficients simultaneously in
high-dimensional linear models while our procedure does not rely on
any specific model assumptions. All the results are based on 1000
Monte Carlo replications.

% % % % % % % % % % % % % % % % % % % % % % % % % %
% % % % % % % % % % % % % % % % % % % % % % % % % %
% % % % % % % % %EXAMPLE 1 % %% % % % % % % %  % %
% % % % % % % % % % % % % % % % % % % % % % % % % %
% % % % % % % % % % % % % % % % % % % % % % % % % %

\begin{example}\label{ex51}
{\rm
Following Zhong and Chen (2011), consider the
linear regression model with the simple random design,
\[
Y_{i}=X_i^T\beta+\epsilon_{i},\quad i=1,2,\dots,n,
\]
where $ X_i=(x_{i1}, x_{i2},\dots,x_{ip})^T$ is generated from the
following moving average model
\begin{equation}
\label{sim:DGP} x_{ij}=\alpha_1 z_{ij} +\alpha_2
z_{i(j+1)}+\cdots+\alpha_T z_{i(j+T-1)}+\mu_j,
\end{equation}
for $j=1,\dots,p$, and $T=10,20$. Here
$Z_i=(z_{i1},\dots,z_{i(p+T-1)})^T\overset{i.i.d.}{\sim}
N(0,I_{p+T-1})$ and $\{\mu_j\}_{j=1}^p$ are fixed realizations from
Uniform(2,3), i.e., the uniform distribution on $(2,3)$. The coefficients $\{\alpha_k\}_{k=1}^T$ are generated
from Uniform(0,1) and are kept fixed once generated. We consider two
distributions for $\epsilon_i$ namely, $N(0,4)$ and centralized
gamma distribution with shape parameter 1 and scale parameter 0.5.
}
\end{example}

Under the null, $\beta=\mathbf{0}_{p\times 1}$. Similar to Zhong and Chen (2011), we consider two
configurations under alternatives. One is the non-sparse case, in
which we equally allocate the first half $\beta_j$ to be positive so
that they all have the same magnitude of signals. The other one is
the sparse case, where we choose the first five elements of $\beta$
to be positive and also equal. In both cases, we fix
$|\beta|_p:=(\sum^{p}_{j=1}\beta_j^2)^{1/2}$ at the level 0.06. In
addition, we consider $n=40,60,80$, and
$p=34,54,76,310,400,550$.

Table \ref{table:eg51} presents the empirical sizes and powers of the
proposed test and the ZC test. The empirical sizes of both tests are reasonably close to the
10\% nominal level for two different error distributions. At the
5\% significance level, however, we see both tests have slightly
inflated rejection probabilities under the null. This is due to the
skewness of the test statistics in finite samples, which was also
pointed out in Zhong and Chen (2011). We provide in Figure 1 the density plot of the proposed test statistic as well as the density of standard Normal for $T=20$. The plot for $T=10$ is very similar and therefore omitted here.  As we can see from the plot, when sample size and dimension are both small, the distribution of the test statistic is skewed to the right, which explains the inflated type I error; when we increase dimension and sample size, the sampling distribution is getting closer to the standard Normal. Overall, both tests perform reasonably well under the null.

\begin{figure}[ht!]
\label{fig:density}
\centering
\includegraphics[width=5in]{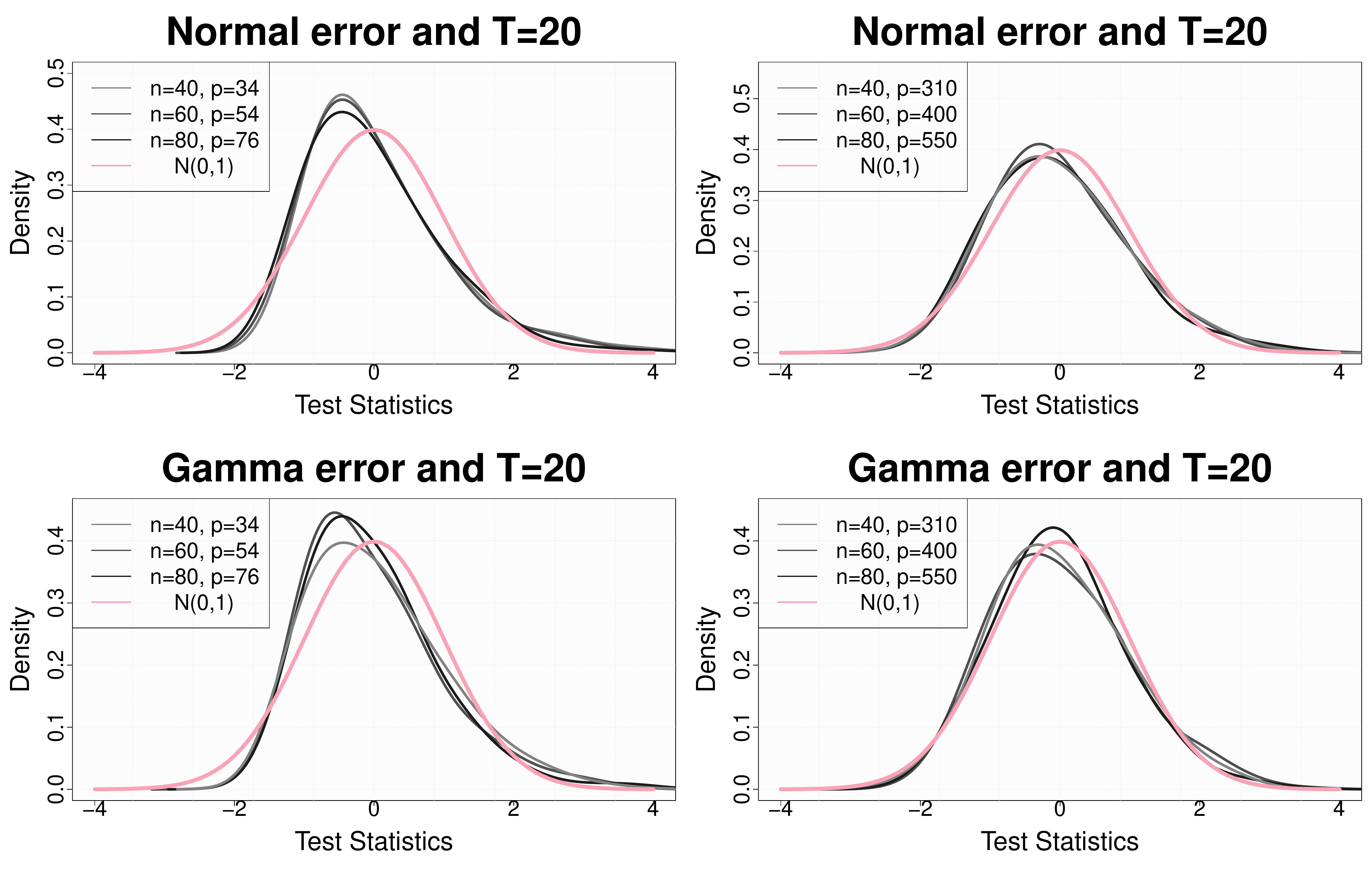}
\caption{Density curves of the sampling distribution of the proposed test statistic (green lines) and the standard normal distribution (red line). }
\end{figure}

For the empirical powers, since the underlying model is indeed a
simple linear model, we expect the ZC test to have more power under
correct model specification than our test. From Table
\ref{table:eg51}, we see that for the non-sparse alternative, the
powers are quite high and close to each other for the two tests; the
powers are comparable for two error distributions. Under the
sparse alternative, both tests have substantial power loss when the dimension is much larger than the sample size indicating that both the mdd and ZC tests target for dense
alternatives. Overall, our test is highly
comparable to the ZC test under the simple linear model.

% % % % % % % % % % % % % % % % % % % % % % % % % %
% % % % % % % % % % % % % % % % % % % % % % % % % %
% % % % % % % % %EXAMPLE 2 % %% % % % % % % % % % %
% % % % % % % % % % % % % % % % % % % % % % % % % %
% % % % % % % % % % % % % % % % % % % % % % % % % %

\begin{example}\label{ex52}
{\rm
The second experiment is designed to evaluate the
performance of the tests under non-linear model:
\[
Y_{i}=\sqrt{\sum^{p}_{j=1}\beta_j x_{ij}^2}+\epsilon_{i},\quad
i=1,2,\dots,n,
\]
where the covariates  $X_i=(x_{i1}, x_{i2},...,x_{ip})^T$ are
generated from the
following model:
\begin{equation}\label{eq:cov}
x_{ij}=(\varsigma_{ij}+\varsigma_{i0})/\sqrt{2},\quad j=1,2,...,p,
\end{equation}
where
$(\varsigma_{i0},\varsigma_{i1},\dots,\varsigma_{ip})^T\overset{i.i.d.}{\sim}
N(0,I_{p+1})$. Hence the covariates are strongly correlated with
the pairwise correlation equal to 0.5. We consider three
configurations for the error, $N(0,1)$, $t_3$, and $\chi^2-1$
respectively.
}
\end{example}

Under the null, we set $\beta=(\beta_1,\dots,\beta_p)^T=\mathbf{0}_{p\times 1}$. Under
the non-sparse alternative, we set $\beta_j=\mathbf{1}\{1\leq j\leq \lfloor n/2\rfloor \}$.
For the sparse alternative, $\beta_j=\mathbf{1}\{1\leq j\leq 5\}$. We consider $n=100$ and $p=50, 100, 200$.

 Table
\ref{table:eg52} reports the empirical size and power for the two tests. The size phenomenon is very similar to that of Example \ref{ex51}. Under the alternatives, the mdd test consistently outperforms the ZC test in all cases. We
observe that under the nonlinear model, the ZC test exhibits very
low powers. By contrast, our proposed test is much more powerful (especially in the nonsparse case)
regardless of $p$ and the error distribution.

\subsection{Conditional quantile independence} \label{sec:simquan}

In this subsection, we carry out additional simulations to evaluate
the performance of the proposed test for conditional quantile
independence.

% % % % % % % % % % % % % % % % % % % % % % % % % %
% % % % % % % % % % % % % % % % % % % % % % % % % %
% % % % % % % %  conditional quantile  % %% % % % %
% % % % % % % % % % % % % % % % % % % % % % % % % %
% % % % % % % % % % % % % % % % % % % % % % % % % %

\begin{example}\label{ex53}
{\rm
Consider the following mixture distribution:
\[
Y= (1+Z) \cdot  \epsilon - (1+X^T \beta) \cdot (1-\epsilon)
 \]
 where $\epsilon$ is a Bernoulli random variable with success probability 0.5; $Z$ is from Gamma(2,2) and   $X$ has $p$ components generated from i.i.d Gamma(6,1).
}
\end{example}
Here $\beta=\mathbf{0}_{p\times 1}$ under the null and it is positive for some components under the alternative. Hence $Y$ has a mixture distribution with half probability to be positive or negative.  We also consider non-sparse and sparse alternative and fix $|\beta|_p=0.06$ as in Example \ref{ex51}. We set $p=50,100,200$, $n=100$ and consider three different quantile levels, $\tau=0.25, 0.5, 0.75$.

It is not difficult to see that when $\beta=0$, $Y$ and $X$ are independent at all quantile levels; when $\beta\not=0$, $Y$ is conditionally quantile independent of $X$ only when $\tau\ge 0.5$. Furthermore, when $\beta\not=0$,  the local quantile independence described in Assumption \ref{ass2}  is satisfied when $\tau=0.75$, but does not hold at $\tau=0.5$.
As presented in Table \ref{table:eg53}, our proposed test  demonstrates nontrivial power only at quantile level 0.25 under both the sparse and non-sparse alternative, which is consistent with the theory.
%For $\tau=0.5, 0.75$, the rejection rate is reasonably close to the nominal level at 10\% and slightly inflated at 5\%.
% This is consistent with the fact that 0.5 and 0.75 level quantiles of $Y$ do not dependent on $X$.

\begin{example}\label{ex54}
{\rm
This example considers a simple model with
heteroscedasticity:
\[Y_{i}= \left(1+ X_i^T \beta \right)^2\epsilon_{i} \]
where $X_i=(x_{i1},x_{i2},\dots,x_{ip})^T$ is a $p$-dimensional
vector of covariates generated as in (\ref{eq:cov})  and $\epsilon_i$ is the error independent of $X_i$. All the
other configurations are the same as Example \ref{ex52}.
}
\end{example}

We examine the sizes and powers under both non-sparse
alternatives and sparse alternatives with four different error
distributions, i.e., $N(0,1)$, $t_3$, Cauchy(0,1) and $\chi^2_1-1$.
We are interested in testing the conditional quantile independence
at $\tau=0.25,0.5,0.75$. Consider $p=50,100,200$, and $n=100$. For
the purpose of comparison, we also include the results using our
MDD-based test for conditional mean independence. Note that when $\beta$ is nonzero, conditional quantile independence only holds at $\tau=0.5$ for the errors of
three types (Normal, Student-t and Cauchy), while conditional quantile dependence is present for other combinations of
quantile levels and error distributions. Further note that the local quantile independence assumption (Assumption \ref{ass2}) does not hold
when $\beta\not=0$ and $\tau=0.5$.

Table \ref{table:eg54} shows the empirical sizes and powers for various configurations. The sizes are generally precise at 10\%
level and slightly inflated at 5\% level under $H_0$, and they do not seem to
depend on the error distribution much. The empirical powers
apparently depend on the error distribution. For symmetric distributions, $N(0,1)$, $t_3$ and Cauchy(0,1) error, they have comparable powers at $\tau=0.25$ and $0.75$.
The powers suffer big reduction under sparse alternatives.
In contrast, power for $\chi^2_1-1$ is almost 1 under the non-sparse alternative at $\tau=0.25$ and remains high even under sparse alternative; at $\tau=0.75$, the proposed test has moderate power in the non-sparse case and suffers a great power loss in the sparse case. This phenomenon may be related to the fact that the chi-square distribution is skewed to the right.
At  $\tau=0.5$, we observe a significant upward size distortion for all the three symmetric distribution (i.e., Normal, student-$t $ and Cauchy), where the conditional median does not depend on the covariates (i.e., under the null). We speculate that this size distortion is related to
the violation of local quantile independence in Assumption \ref{ass2}, and indicates the necessity of Assumption \ref{ass2} in obtaining the
weak convergence of our test statistics to the standard normal distribution. For non-symmetric distribution (Chi-square), the power is satisfactory at the median.

In
comparison, the MDD-based test for conditional mean has  power close to the level for $N(0,1)$, $t_3$ and  $\chi^2_1-1$ as the covariates only contribute to the conditional variance of the response but not to the conditional mean. In addition, when the mean does not exist, for
example in the Cauchy(0,1) case, we observe that the MDD-based test
for conditional mean has very little power. Therefore, the test for
the conditional quantile dependence can be a good complement to the
conditional mean independence test, especially for  models that exhibit conditional heteroscedasticity and heavy tail.

% % % % % % % % % % % % % % % % % % % % % % % % % %
% % % % % % % % % % % % % % % % % % % % % % % % % %
% % % % % % % %  Bootstrap   % %% % % % %
% % % % % % % % % % % % % % % % % % % % % % % % % %
% % % % % % % % % % % % % % % % % % % % % % % % % %

\subsection{Wild bootstrap}\label{sec:simboot}
\begin{figure}[ht!]
\centering
\includegraphics[width=4in]{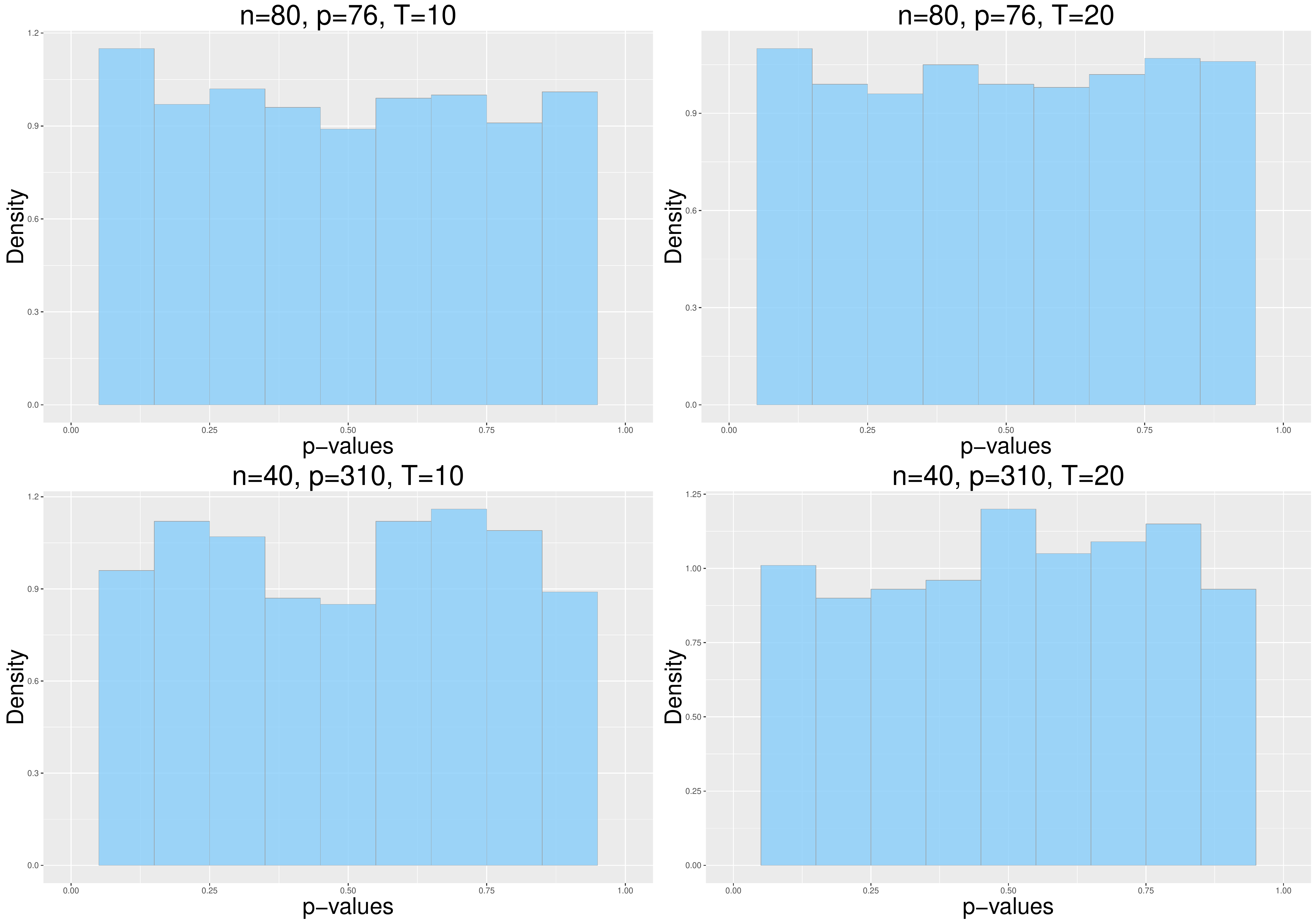}
\caption{Histograms of the p-values of the test statistics using the studentized bootstrap critical value. }
\label{fig:pvalue}
\end{figure}

%To alleviate the size distortion we observed in the simulation section, and improve the finite sample performance when sample size is small, we propose a wild bootstrap approach to approximate
%the null distribution of the MDD-based test statistic as a useful alternative to normal approximation.
In this section, we implement the Wild bootstrap proposed in Section \ref{sec:boot}. We again consider Example \ref{ex51} under the null of $\beta=\mathbf{0}_{p\times 1}$ in Section \ref{sec:sim}. The results are summarized in Table \ref{table:boot}, where we
compare the type I errors delivered by normal approximation (mdd) with the wild bootstrap (Boot) approach.
It is observed that the size is usually more accurate using wild bootstrap. This is not surprising as it is known in the literature that application of bootstrap method to pivotal statistic has better asymptotic properties, e.g.,
faster rate of convergence of the actual significance level to the nominal significance level, see the discussions in Westfall and Young (1993).
To understand whether the bootstrap distribution is mimicking the finite sample distribution of our test statistics under the null, we also provide in Figure \ref{fig:pvalue} the histograms of the p-values corresponding to the wild bootstrap. As we can see from these histograms, the p-values are close to being uniformly distributed between 0 and 1 indicating the bootstrap distribution is consistently approximating the null distribution of the test statistics.

%\begin{figure}[ht!]
%\centering
%\includegraphics*[width=3.5in]{hist.pdf}
%\caption{Histogram of the test statistics (blue) and the studentized bootstrap statistics (pink) from one simulation with $n=80, p=76, T=20$ for Normal error. }
%\label{fig:hist}
%\end{figure}

\section{Conclusion}\label{sec:con}

In this paper, we proposed new tests for testing conditional mean
and conditional quantile independence in the high dimensional
setting, where the number of covariates is allowed to grow rapidly
with the sample size. Our test statistics are built on the MDD, a recently proposed metric for quantifying conditional mean dependence, and do not
involve any tuning parameters.  A distinctive feature of our test is
that it is totally model free and is able to capture nonlinear
dependence which  cannot be revealed by using traditional  linear
dependence metric. Our test for conditional mean independence can be viewed as a nonparametric model-free
counterpart of Zhong and Chen's test, which was developed for high dimensional linear models with homoscedastic errors;  our test for conditional quantile independence  seems to be the first
 simultaneous test  proposed in the high dimensional setting that detects conditional quantile dependence.
 Besides the methodological advance, our
theoretical analysis reveals a key assumption generally needed in
establishing the asymptotic normality for a broad class of
$U$-statistic based tests in high dimension, and shed light on the
asymptotic theory for high dimensional $U$-statistic.

Our analysis of local
asymptotic power shows that our test is  less powerful than Zhong and Chen's test
when the high dimensional linear model holds, but the efficiency loss is fairly moderate.
As a tradeoff, our test for conditional mean dependence achieves great robustness against model mis-specification, under which
the power of Zhong and Chen's test can greatly deteriorate, as we demonstrate in our simulations. In particular,
 the proposed conditional mean dependence test can
significantly outperform the modified $F$-test in Zhong and Chen
(2011) when the underlying relationship between $Y$ and $X$ is
nonlinear. Since there seem no well-developed methods yet to check
the validity of a high dimensional linear model and it is often
unknown whether the response is linearly related to covariates, it
might be safe and more reliable to apply our model-free dependence
test in practice.

To conclude, we point out a few future research directions. Firstly,
given a set of variables $\mathcal{X}^*$ which are believed to have
significant impact on the response $Y$, it is of interest to test
whether the conditional mean or quantile of $Y$ depends on another
set of variables $\mathcal{X}$, adjusting for the effect of
$\mathcal{X}^*$. Based on the recently proposed concept of partial
MDD [Park et al. (2015)], we expect that the results could be
extended to test the conditional mean or quantile dependence
controlling for  a set of significant variables. Secondly, we only
handle a single quantile level in our test, and it would be
interesting to extend our test to testing conditional quantile
independence over a range of quantile levels, say $[\tau_1,\tau_2],$
where $0<\tau_1<\tau_2<1$.  Thirdly, the studentized wild bootstrap showed great promise in reducing the
size distortion associated with normal approximation in Section \ref{sec:simboot}, but there is yet
 any theory that states the second order
accuracy for the wild bootstrap approach. Fourthly, our MDD-based tests mainly target the marginal conditional mean dependence and may
fail to capture conditional mean dependence at higher order, which needs to be taken into account in developing new tests.
 Finally, our test statistics are of
$L^2$-type, and it may be of interest to develop a
$L^{\infty}$-version, which targets for sparse and strong
alternatives. We leave these topics for future research.

%This problem is also closely related with the
%estimation of graphical model.

\vskip 0.3cm

\textbf{Acknowledgments.} The authors would like to thank the Associate Editor
and the reviewers for their constructive comments and helpful suggestions, which substantially
improved the paper. We also thank Dr. Ping-Shou Zhong for providing the data set
and the R codes used in our simulation studies and data illustration. We are also grateful to
Dr. Min Qian for providing us with the matlab code for the tests in McKeague and Qian (2015).

\vskip 0.3cm

\centerline{SUPPLEMENTARY MATERIAL}
This supplement contains proofs of the
main results in the paper, extension to factorial designs, additional discussions and numerical results.

% % % % % % % % % % % % % % % % % % % % % % % % % %
% % % % % % % % % % % % % % % % % % % % % % % % % %
% % % % % % % % % % % EXAMPLE 1 % % % % % % % %
% % % % % % % % % % % % % % % % % % % % % % % % % %
% % % % % % % % % % % % % % % % % % % % % % % % % %

\begin{table}[htp!]
\footnotesize
    \centering
    \caption{Empirical sizes and powers of the MDD-based test for conditional mean independence and the ZC test at significance levels 5\% and 10\% for Example \ref{ex51}.}
    \label{table:eg51}

    \begin{tabular}{cccccccccccc}
        \toprule
        &&&&\multicolumn{4}{c}{Normal error}&\multicolumn{4}{c}{Gamma error}\\
        &&&&\multicolumn{2}{c}{mdd}&\multicolumn{2}{c}{ZC}&\multicolumn{2}{c}{mdd}&\multicolumn{2}{c}{ZC}\\
        $T$ & case & $n$ & $p$ & 5\% & 10\% & 5\% & 10\% & 5\% & 10\% & 5\% & 10\% \\
        \hline

        \multirow{18}{*}{\scriptsize $10$} & \multirow{6}{*}{\small $H_0$ }   &   40 &   34 & 0.069 & 0.125 & 0.073 & 0.116 & 0.055 & 0.095 & 0.055 & 0.094 \\
   &  &   60 &   54 & 0.075 & 0.115 & 0.068 & 0.114 & 0.072 & 0.118 & 0.067 & 0.113 \\
   &  &   80 &   76 & 0.069 & 0.114 & 0.075 & 0.107 & 0.063 & 0.101 & 0.056 & 0.100 \\
   &  &   40 &  310 & 0.043 & 0.088 & 0.048 & 0.092 & 0.049 & 0.093 & 0.053 & 0.088 \\
   &  &   60 &  400 & 0.060 & 0.104 & 0.058 & 0.106 & 0.052 & 0.098 & 0.054 & 0.097 \\
   &  &   80 &  550 & 0.060 & 0.106 & 0.065 & 0.107 & 0.047 & 0.090 & 0.047 & 0.082 \\
        \cline{3-12}
           & \multirow{6}{*}{\parbox{1.2cm}{\small \centering non-sparse $H_a$}} &         40 &   34 & 0.634 & 0.714 & 0.678 & 0.746 & 0.672 & 0.748 & 0.704 & 0.766 \\
           &  &   60 &   54 & 0.781 & 0.833 & 0.814 & 0.866 & 0.803 & 0.846 & 0.817 & 0.871 \\
           &  &   80 &   76 & 0.848 & 0.901 & 0.875 & 0.914 & 0.866 & 0.908 & 0.883 & 0.915 \\
           &  &   40 &  310 & 0.262 & 0.354 & 0.287 & 0.392 & 0.276 & 0.391 & 0.317 & 0.419 \\
           &  &   60 &  400 & 0.385 & 0.489 & 0.410 & 0.530 & 0.386 & 0.516 & 0.406 & 0.523 \\
           &  &   80 &  550 & 0.435 & 0.544 & 0.473 & 0.574 & 0.454 & 0.567 & 0.474 & 0.586 \\
        \cline{3-12}
        &\multirow{6}{*}{\parbox{1.2cm}{\small \centering sparse $H_a$}} &   40 &   34 & 0.310 & 0.396 & 0.337 & 0.411 & 0.331 & 0.410 & 0.357 & 0.435 \\
           &  &   60 &   54 & 0.313 & 0.411 & 0.344 & 0.435 & 0.368 & 0.455 & 0.375 & 0.459 \\
           &  &   80 &   76 & 0.350 & 0.442 & 0.354 & 0.436 & 0.343 & 0.427 & 0.368 & 0.449 \\
           &  &   40 &  310 & 0.067 & 0.134 & 0.086 & 0.145 & 0.081 & 0.157 & 0.088 & 0.157 \\
           &  &   60 &  400 & 0.110 & 0.176 & 0.108 & 0.182 & 0.098 & 0.161 & 0.094 & 0.174 \\
           &  &   80 &  550 & 0.102 & 0.171 & 0.110 & 0.171 & 0.093 & 0.157 & 0.096 & 0.167 \\
        \hline
        \multirow{18}{*}{\scriptsize $20$} & \multirow{6}{*}{\small $H_0$ } & 40 &   34 & 0.069 & 0.114 & 0.081 & 0.118 & 0.063 & 0.107 & 0.062 & 0.090 \\
           &  &   60 &   54 & 0.068 & 0.100 & 0.076 & 0.108 & 0.052 & 0.089 & 0.058 & 0.093 \\
           &  &   80 &   76 & 0.058 & 0.102 & 0.059 & 0.100 & 0.073 & 0.111 & 0.072 & 0.114 \\
           &  &   40 &  310 & 0.057 & 0.098 & 0.063 & 0.101 & 0.077 & 0.105 & 0.062 & 0.096 \\
           &  &   60 &  400 & 0.061 & 0.089 & 0.071 & 0.099 & 0.045 & 0.092 & 0.046 & 0.096 \\
           &  &   80 &  550 & 0.070 & 0.104 & 0.074 & 0.107 & 0.054 & 0.100 & 0.059 & 0.099 \\
        \cline{3-12}
        & \multirow{6}{*}{\parbox{1.2cm}{\small \centering non-sparse $H_a$}} &   40 &   34 & 0.976 & 0.984 & 0.985 & 0.990 & 0.971 & 0.979 & 0.977 & 0.985 \\
           &  &   60 &   54 & 1.000 & 1.000 & 1.000 & 1.000 & 0.997 & 0.998 & 0.998 & 0.999 \\
           &  &   80 &   76 & 1.000 & 1.000 & 1.000 & 1.000 & 1.000 & 1.000 & 1.000 & 1.000 \\
           &  &   40 &  310 & 0.848 & 0.903 & 0.877 & 0.918 & 0.858 & 0.893 & 0.877 & 0.917 \\
           &  &   60 &  400 & 0.947 & 0.967 & 0.961 & 0.981 & 0.946 & 0.973 & 0.955 & 0.976 \\
           &  &   80 &  550 & 0.990 & 0.994 & 0.994 & 0.996 & 0.981 & 0.990 & 0.983 & 0.990 \\
        \cline{3-12}
        &\multirow{6}{*}{\parbox{1.2cm}{\small \centering sparse $H_a$}} & 40 &   34 & 0.705 & 0.763 & 0.737 & 0.782 & 0.715 & 0.779 & 0.741 & 0.786 \\
           &  &   60 &   54 & 0.743 & 0.808 & 0.764 & 0.821 & 0.753 & 0.812 & 0.768 & 0.819 \\
           &  &   80 &   76 & 0.764 & 0.813 & 0.782 & 0.835 & 0.743 & 0.799 & 0.765 & 0.819 \\
           &  &   40 &  310 & 0.137 & 0.212 & 0.156 & 0.236 & 0.185 & 0.274 & 0.193 & 0.269 \\
           &  &   60 &  400 & 0.186 & 0.268 & 0.209 & 0.275 & 0.183 & 0.283 & 0.195 & 0.285 \\
           &  &   80 &  550 & 0.176 & 0.271 & 0.200 & 0.291 & 0.208 & 0.288 & 0.213 & 0.289 \\
        \hline
    \end{tabular}
\end{table}

% % % % % % % % % % % % % % % % % % % % % % % % % %
% % % % % % % % % % % % % % % % % % % % % % % % % %
% % % % % % % % % % % EXAMPLE 2 % % % % % % % % % %
% % % % % % % % % % % % % % % % % % % % % % % % % %
% % % % % % % % % % % % % % % % % % % % % % % % % %

\begin{table}[t]
\footnotesize
    \centering
    \caption{\small Empirical sizes and powers of the MDD-based test for conditional mean independence and the ZC test at significance levels 5\% and 10\% for Examples \ref{ex52}.}
    \label{table:eg52}
    \begin{tabular}{ccccccc}
        \toprule
        &&&\multicolumn{2}{c}{mdd}&\multicolumn{2}{c}{ZC}\\
        error & case & $p$ & 5\% & 10\% & 5\% & 10\%\\
        \hline
           \multirow{9}{*}{$N(0,1)$} & \multirow{3}{*}{$H_0$ } &  50 & 0.078 & 0.112 & 0.075 & 0.098 \\
              &  &  100 & 0.075 & 0.111 & 0.078 & 0.109 \\
              &  &  200 & 0.065 & 0.091 & 0.062 & 0.089 \\
           \cline{3-7}
 &  \multirow{3}{*}{\parbox{1.2cm}{\centering non-sparse $H_a$}} &   50 & 0.927 & 0.990 & 0.200 & 0.221 \\
    &  &  100 & 0.970 & 0.997 & 0.213 & 0.238 \\
    &  &  200 & 0.980 & 0.998 & 0.230 & 0.249 \\
        \cline{3-7}
        &\multirow{3}{*}{\parbox{1.2cm}{\centering sparse $H_a$}} &  50 & 0.428 & 0.583 & 0.143 & 0.176 \\
           &  &  100 & 0.370 & 0.519 & 0.147 & 0.172 \\
           &  &  200 & 0.331 & 0.477 & 0.129 & 0.161 \\
        \hline
         \multirow{9}{*}{$t_3$} & \multirow{3}{*}{$H_0$ } &   50 & 0.084 & 0.110 & 0.080 & 0.104 \\
            &  &  100 & 0.075 & 0.106 & 0.072 & 0.097 \\
            &  &  200 & 0.082 & 0.112 & 0.063 & 0.091 \\
                 \cline{3-7}
          & \multirow{3}{*}{\parbox{1.2cm}{\centering non-sparse $H_a$}} &   50 & 0.759 & 0.877 & 0.178 & 0.206 \\
             &  &  100 & 0.870 & 0.955 & 0.182 & 0.209 \\
             &  &  200 & 0.939 & 0.983 & 0.199 & 0.226 \\
        \cline{3-7}
        &\multirow{3}{*}{\parbox{1.2cm}{\centering sparse $H_a$}} &  50 & 0.246 & 0.355 & 0.108 & 0.141 \\
           &  &  100 & 0.208 & 0.302 & 0.105 & 0.136 \\
           &  &  200 & 0.209 & 0.311 & 0.108 & 0.131 \\
        \hline
        \multirow{9}{*}{$\chi^2_1-1$} & \multirow{3}{*}{$H_0$ }&
        50 & 0.083 & 0.114 & 0.083 & 0.110 \\
           &  &  100 & 0.061 & 0.096 & 0.058 & 0.086 \\
           &  &  200 & 0.058 & 0.084 & 0.047 & 0.078 \\
                \cline{3-7}
         & \multirow{3}{*}{\parbox{1.2cm}{\centering non-sparse $H_a$}} &  50 & 0.838 & 0.937 & 0.162 & 0.202 \\
            &  &  100 & 0.933 & 0.983 & 0.198 & 0.226 \\
            &  &  200 & 0.967 & 0.995 & 0.216 & 0.245 \\
        \cline{3-7}
        &\multirow{3}{*}{\parbox{1.2cm}{\centering sparse $H_a$}} &    50 & 0.269 & 0.384 & 0.113 & 0.150 \\
           &  &  100 & 0.247 & 0.364 & 0.124 & 0.155 \\
           &  &  200 & 0.233 & 0.323 & 0.111 & 0.141 \\
        \hline
    \end{tabular}
\end{table}

% % % % % % % % % % % % % % % % % % % % % % % % % %
% % % % % % % % % % % % % % % % % % % % % % % % % %
% % % % % % % %  conditional quantile  % %% % % % %
% % % % % % % % % % % % % % % % % % % % % % % % % %
% % % % % % % % % % % % % % % % % % % % % % % % % %

\begin{table}[b]
\footnotesize
    \centering
    \caption{Empirical sizes and powers of the MDD-based test for conditional quantile independence at significance levels 5\% and 10\% for Example \ref{ex53}.}
    \label{table:eg53}
\begin{tabular}{cccccccc}
  \hline

  & &   \multicolumn{2}{c}{$ H_0$}  &   \multicolumn{2}{c}{Non-sparse $H_0 / H_a$ }  &   \multicolumn{2}{c}{Sparse $H_0/H_a$} \\
  $\tau$ & $p$ & 5\% & 10\%  & 5\% & 10\%  & 5\% & 10\%\\
  \hline
 \multirow{3}{*}{\scriptsize 0.25}    &   50 & 0.054 & 0.104 & 0.522 & 0.644 & 0.554 & 0.663 \\
    &  100 & 0.048 & 0.094 & 0.322 & 0.464 & 0.352 & 0.473 \\
    &  200 & 0.044 & 0.100 & 0.231 & 0.368 & 0.230 & 0.349 \\
   \hline
 \multirow{3}{*}{\scriptsize 0.50}    &  50 & 0.058 & 0.106 & 0.071 & 0.129 & 0.069 & 0.130 \\
    &  100 & 0.047 & 0.093 & 0.062 & 0.104 & 0.059 & 0.104 \\
    &  200 & 0.056 & 0.093 & 0.065 & 0.104 & 0.058 & 0.107 \\
   \hline
 \multirow{3}{*}{\scriptsize 0.75}   & 50 & 0.050 & 0.108 & 0.050 & 0.108 & 0.050 & 0.108 \\
    &  100 & 0.058 & 0.103 & 0.058 & 0.103 & 0.058 & 0.103 \\
    &  200 & 0.067 & 0.120 & 0.067 & 0.120 & 0.067 & 0.120 \\
   \hline
\end{tabular}
\end{table}

\begin{table}[h]
\footnotesize
    \centering
    \caption{Empirical sizes and powers of the MDD-based test for conditional quantile independence at significance levels 5\% and 10\% for Example \ref{ex54}.}
    \label{table:eg54}
    \begin{tabular}{ccccccccccc}
        \toprule
        &&&\multicolumn{2}{c}{$N(0,1)$}&\multicolumn{2}{c}{$t_3$}&\multicolumn{2}{c}{Cauchy(0,1)} &\multicolumn{2}{c}{$\chi^2_1-1$}\\
        $\tau$ & case & $p$ & 5\% & 10\% & 5\% & 10\% & 5\% & 10\% & 5\% & 10\%  \\
        \hline
       \multirow{9}{*}{0.25} & \multirow{3}{*}{$H_0$ } &   50 & 0.084 & 0.112 & 0.058 & 0.077 & 0.073 & 0.094 & 0.074 & 0.099 \\
          &  &  100 & 0.073 & 0.104 & 0.067 & 0.094 & 0.069 & 0.101 & 0.071 & 0.099 \\
          &  &  200 & 0.064 & 0.085 & 0.058 & 0.095 & 0.062 & 0.082 & 0.069 & 0.094 \\
        \cline{3-11}
        & \multirow{3}{*}{\parbox{1.2cm}{ \centering non-sparse $H_a$}} &   50 & 0.573 & 0.640 & 0.494 & 0.554 & 0.368 & 0.424 & 0.999 & 0.999 \\
           &  &  100 & 0.840 & 0.864 & 0.727 & 0.778 & 0.578 & 0.632 & 1.000 & 1.000 \\
           &  &  200 & 0.936 & 0.956 & 0.925 & 0.941 & 0.842 & 0.880 & 1.000 & 1.000 \\
        \cline{3-11}
        &\multirow{3}{*}{\parbox{1.2cm}{\centering sparse $H_a$}} &    50 & 0.207 & 0.243 & 0.146 & 0.196 & 0.122 & 0.162 & 0.951 & 0.970 \\
           &  &  100 & 0.198 & 0.240 & 0.154 & 0.196 & 0.125 & 0.161 & 0.956 & 0.967 \\
           &  &  200 & 0.199 & 0.235 & 0.174 & 0.224 & 0.125 & 0.168 & 0.952 & 0.965 \\
        \hline
        \multirow{9}{*}{0.5} & \multirow{3}{*}{$H_0$ } & 50 & 0.083 & 0.106 & 0.076 & 0.097 & 0.068 & 0.099 & 0.064 & 0.089 \\
           &  &  100 & 0.072 & 0.100 & 0.063 & 0.087 & 0.057 & 0.092 & 0.062 & 0.091 \\
           &  &  200 & 0.073 & 0.100 & 0.053 & 0.073 & 0.071 & 0.095 & 0.063 & 0.088 \\
        \cline{3-11}
        & \multirow{3}{*}{\parbox{1.2cm}{ \centering non-sparse $H_0/ H_a$}} &  50 & 0.110 & 0.145 & 0.093 & 0.125 & 0.091 & 0.127 & 0.704 & 0.755 \\
           &  &  100 & 0.151 & 0.184 & 0.123 & 0.168 & 0.126 & 0.169 & 0.862 & 0.888 \\
           &  &  200 & 0.230 & 0.274 & 0.232 & 0.269 & 0.196 & 0.242 & 0.939 & 0.953 \\
        \cline{3-11}
        &\multirow{3}{*}{\parbox{1.2cm}{\centering sparse $H_0/ H_a$}} &  50 & 0.083 & 0.110 & 0.075 & 0.100 & 0.077 & 0.099 & 0.234 & 0.288 \\
           &  &  100 & 0.079 & 0.108 & 0.063 & 0.088 & 0.066 & 0.107 & 0.252 & 0.299 \\
           &  &  200 & 0.077 & 0.109 & 0.057 & 0.085 & 0.079 & 0.098 & 0.242 & 0.303 \\
        \hline
        \multirow{9}{*}{0.75} & \multirow{3}{*}{$H_0$ }&
        50 & 0.067 & 0.098 & 0.080 & 0.108 & 0.062 & 0.084 & 0.065 & 0.085 \\
           &  &  100 & 0.057 & 0.085 & 0.078 & 0.098 & 0.076 & 0.111 & 0.073 & 0.095 \\
           &  &  200 & 0.064 & 0.095 & 0.063 & 0.085 & 0.069 & 0.095 & 0.065 & 0.088 \\
        \cline{3-11}
        & \multirow{3}{*}{\parbox{1.2cm}{ \centering non-sparse $H_a$}}& 50 & 0.585 & 0.638 & 0.516 & 0.577 & 0.384 & 0.445 & 0.129 & 0.161 \\
           &  &  100 & 0.834 & 0.875 & 0.744 & 0.800 & 0.587 & 0.645 & 0.199 & 0.233 \\
           &  &  200 & 0.944 & 0.961 & 0.925 & 0.940 & 0.824 & 0.865 & 0.329 & 0.385 \\
        \cline{3-11}
        &\multirow{3}{*}{\parbox{1.2cm}{\centering sparse $H_a$}} &  50 & 0.187 & 0.235 & 0.191 & 0.232 & 0.131 & 0.168 & 0.076 & 0.096 \\
           &  &  100 & 0.193 & 0.237 & 0.197 & 0.236 & 0.131 & 0.172 & 0.092 & 0.116 \\
           &  &  200 & 0.195 & 0.244 & 0.167 & 0.213 & 0.135 & 0.170 & 0.078 & 0.109 \\
        \hline
        \multirow{9}{*}{mean} & \multirow{3}{*}{$H_0$ }&
        50 & 0.071 & 0.103 & 0.069 & 0.101 & 0.067 & 0.098 & 0.052 & 0.097 \\
           &  &  100 & 0.074 & 0.101 & 0.075 & 0.097 & 0.071 & 0.101 & 0.060 & 0.103 \\
           &  &  200 & 0.067 & 0.094 & 0.078 & 0.105 & 0.066 & 0.088 & 0.054 & 0.096 \\
        \cline{3-11}
        & \multirow{3}{*}{\parbox{1.2cm}{ \centering non-sparse $H_a$}}&  50 & 0.077 & 0.107 & 0.071 & 0.107 & 0.109 & 0.144 & 0.058 & 0.107 \\
           &  &  100 & 0.077 & 0.102 & 0.063 & 0.091 & 0.117 & 0.143 & 0.057 & 0.107 \\
           &  &  200 & 0.076 & 0.109 & 0.078 & 0.117 & 0.146 & 0.166 & 0.059 & 0.102 \\
        \cline{3-11}
        &\multirow{3}{*}{\parbox{1.2cm}{\centering sparse $H_a$}} &  50 & 0.072 & 0.103 & 0.066 & 0.102 & 0.080 & 0.110 & 0.057 & 0.094 \\
           &  &  100 & 0.067 & 0.100 & 0.059 & 0.084 & 0.080 & 0.105 & 0.060 & 0.112 \\
           &  &  200 & 0.068 & 0.102 & 0.075 & 0.108 & 0.091 & 0.116 & 0.053 & 0.095 \\
        \hline

    \end{tabular}
\end{table}

\begin{table}[h]
\small
    \centering
    \caption{Size comparison for the proposed test using the normal approximation (mdd) and the wild bootstrap approximation for Example \ref{ex51} }
        \label{table:boot}
\resizebox{\columnwidth}{!}{%
\begin{tabular}{ccccccc|cccc}
  \hline
        &&&\multicolumn{4}{c}{Normal error}&\multicolumn{4}{c}{Gamma error}\\
        &&  &\multicolumn{2}{c}{5\%}  &\multicolumn{2}{c}{10\%} &\multicolumn{2}{c}{5\%} &\multicolumn{2}{c}{10\%}  \\
 \cmidrule(lr){4-5} \cmidrule(lr){6-7}  \cmidrule(lr){8-9} \cmidrule(lr){10-11}
         $T$  & $n$ & $p$  & mdd & Boot  & mdd & Boot  & mdd & Boot  & mdd & Boot   \\
        \hline
 \multirow{6}{*}{\scriptsize $10$}  &   40 &   34 & 0.069 & 0.055 & 0.125 & 0.109 & 0.055 & 0.044 & 0.095 & 0.087 \\
   &   60 &   54 & 0.075 & 0.068 & 0.115 & 0.113 & 0.072 & 0.053 & 0.118 & 0.115 \\
   &   80 &   76 & 0.069 & 0.053 & 0.114 & 0.106 & 0.063 & 0.048 & 0.101 & 0.093 \\
   &   40 &  310 & 0.043 & 0.037 & 0.088 & 0.083 & 0.049 & 0.046 & 0.093 & 0.085 \\
   &   60 &  400 & 0.060 & 0.053 & 0.104 & 0.101 & 0.052 & 0.052 & 0.098 & 0.096 \\
   &   80 &  550 & 0.060 & 0.058 & 0.106 & 0.101 & 0.047 & 0.039 & 0.090 & 0.090 \\
       \hline
           \multirow{6}{*}{\scriptsize $20$}
   &   40 &   34 & 0.069 & 0.059 & 0.114 & 0.103 & 0.063 & 0.048 & 0.107 & 0.097 \\
   &   60 &   54 & 0.068 & 0.056 & 0.100 & 0.092 & 0.052 & 0.035 & 0.089 & 0.081 \\
   &   80 &   76 & 0.058 & 0.040 & 0.102 & 0.094 & 0.073 & 0.051 & 0.111 & 0.106 \\
   &   40 &  310 & 0.057 & 0.046 & 0.098 & 0.091 & 0.077 & 0.064 & 0.105 & 0.097 \\
   &   60 &  400 & 0.061 & 0.053 & 0.089 & 0.088 & 0.045 & 0.039 & 0.092 & 0.088 \\
   &   80 &  550 & 0.070 & 0.056 & 0.104 & 0.093 & 0.054 & 0.044 & 0.100 & 0.091 \\
   \hline
\end{tabular}
}
\end{table}


\begin{thebibliography}{99}
\bibitem{cq}
\par\noindent\hangindent2.3em\hangafter 1
\textsc{Chen, S. X.} and \textsc{Qin, Y.} (2010). A two sample test
for high dimensional data with applications to gene-set testing.
{\it Ann. Statist.} {\bf 38} 808-835.

\bibitem{et07}
\par\noindent\hangindent2.3em\hangafter 1
\textsc{Efron, B.} and \textsc{Tibshirani, R.} (2007). On testing
the significance of sets of genes. {\it Ann. Appl. Stat.} {\bf 1}
107-129.

\bibitem{fl08}
\par\noindent\hangindent2.3em\hangafter 1
\textsc{Fan, J.} and \textsc{Lv, J.} (2008) Sure independence screening for ultra-high dimensional feature space (with discussion)
{\it J. Roy. Stat. Soc. Ser. B.} {\bf 70}
849-911.

\bibitem{fkn90}
\par\noindent\hangindent2.3em\hangafter 1
\textsc{Fang, K., Kotz, S.} and \textsc{Ng, K.-W.} (1990). {\it Symmetric Multivariate and Related
Distributions}. Chapman \& Hall, London.

\bibitem{fzw13}
\par\noindent\hangindent2.3em\hangafter 1
\textsc{Feng, L., Zou, C.} and \textsc{Wang, Z.} (2013). Rank-based
score tests for high-dimensional regression coefficients. {\it
Electron. J. Stat.} {\bf 7} 2131-2149.

\bibitem{ggh}
\par\noindent\hangindent2.3em\hangafter 1
\textsc{Goeman, J., van de Geer, S. A.} and \textsc{van Houwelingen,
H. C.} (2006). Testing against a high dimensional alternative. {\it
J. R. Statist. Soc. B}. {\bf 68} 477-493.


\bibitem{h1984}
\par\noindent\hangindent2.3em\hangafter 1
\textsc{Hall, P.}  (1984). Central limit theorem for integrated
square error of multivariate nonparametric density estimators. {\it
J. Multivariate Anal.} {\bf 14} 1-16.


\bibitem{hh}
\par\noindent\hangindent2.3em\hangafter 1
\textsc{Hall, P.} and \textsc{Heyde, C. C.} (1980). {\it Martingale
Limit Theory and Its Application.} Academic Press, New York.


\bibitem{hwh13}
\par\noindent\hangindent2.3em\hangafter 1
 \textsc{He, X.}, \textsc {Wang, L.} and \textsc{Hong, H.} (2013). Quantile-adaptive model-free nonlinear feature screening for high-dimensional heterogeneous data.
{\it Ann. Statist.} {\bf 41} 342-369.


\bibitem{kb}
\par\noindent\hangindent2.3em\hangafter 1
\textsc{Koenker, R.} and \textsc{Bassett, G.} (1978). Regression
quantiles. {\it Econometrica}. {\bf 46} 33-50.

\bibitem{lwt}
\par\noindent\hangindent2.3em\hangafter 1
\textsc{Lan, W., Wang, H.} and \textsc{Tsai, C. L.}  (2014). Testing
covariates in high-dimensional regression. {\it Ann. Inst. Statist.
Math.} {\bf 66} 279-301.


\bibitem{lhz}
\par\noindent\hangindent2.3em\hangafter 1
\textsc{Li, Q., Hsiao, C.} and \textsc{Zinn, J.} (2003). Consistent specification tests for semiparametric/nonparametric models based on series estimation methods.
{\it Journal of Econometrics}. {\bf 112} 295-325.


\bibitem{lzl}
\par\noindent\hangindent2.3em\hangafter 1
\textsc{Liu, J., Zhong, W.} and \textsc{Li, R.} (2015) A selective overview of feature screening for ultrahigh-dimensional data. {\it Science China Mathematics.} {\bf 58}
 2033-2054.

\bibitem{lqc}
\par\noindent\hangindent2.3em\hangafter 1
\textsc{Lkhagvadorj, S., Qu, L., Cai, W., Couture, O., Barb, C.,
Hausman, G., Nettleton, D., Anderson, L., Dekkers, J.} and
\textsc{Tuggle, C.} (2009). Microarray gene expression profiles of
fasting induced changes in liver and adipose tissues of pigs
expressing the melanocortin-4 receptor D298N variant. {\it
Physiological Genomics}. {\bf 38} 98-111.

\bibitem{mq}
\par\noindent\hangindent2.3em\hangafter 1
\textsc{McKeague, I. W.} and \textsc{Qian, M.}  (2015) An adaptive resampling
test for detecting the presence of significant predictors. {\it J. Amer. Statist. Assoc.} {\bf 110} 1422-1433.

\bibitem{n07}
\par\noindent\hangindent2.3em\hangafter 1
\textsc{Newton, M., Quintana, F., Den Boon, J., Sengupta, S.} and
\textsc{Ahlquist, P.} (2007). Random-set methods identify distinct
aspects of the enrichment signal in gene-set analysis. {\it Ann.
Appl. Stat.} {\bf 1} 85-106.


\bibitem{psy}
\par\noindent\hangindent2.3em\hangafter 1
\textsc{Park, T.}, \textsc{Shao, X.} and \textsc{Yao, S.} (2015).
Partial martingale difference correlation. {\it Electron. J. Stat.}
{\bf 9} 1492-1517.


\bibitem{ssgf}
\par\noindent\hangindent2.3em\hangafter 1
\textsc{Sejdinovick, D., Sriperumbudur, B., Gretton, A.} and
\textsc{Fukumizu, K.} (2013). Equivalence of distance-based and
RKHS-based statistics in hypothesis testing. {\it Ann. Statist.}
\textbf{41} 2263-2291.


\bibitem{serf}
\par\noindent\hangindent2.3em\hangafter 1
\textsc{Serfling, R. J.} (1980).  {\it Approximation Theorems of
Mathematical Statistics.} Wiley, New York.

\bibitem{sz14}
\par\noindent\hangindent2.3em\hangafter 1
\textsc{Shao, X.} and \textsc{ Zhang, J.} (2014). Martingale
difference correlation and its use in high dimensional variable
screening. {\it J. Amer. Statist. Assoc.} {\bf 109} 1302-1318.


\bibitem{s05}
\par\noindent\hangindent2.3em\hangafter 1
\textsc{Subramanian, A., Mootha, V. K., Mukherjee, S., Ebert, B. L.,
Gillette, M. A., Paulovich, A., Pomeroy, S. L., Golub, T. R.,
Lander, E. S.} and \textsc{Mesirov, J. P.} (2005). Gene set
enrichment analysis: a knowledge-based approach for interpreting
genome-wide expression profiles. {\it Proc. Natl. Acad. Sci. USA}
{\bf 102} 15545-15550


\bibitem{sr13}
\par\noindent\hangindent2.3em\hangafter 1
\textsc{Sz\'ekely, G. J.} and \textsc {Rizzo, M. L.} (2013). The
distance correlation t-test of independence in high dimension. {\it
J. Multivariate Anal.} {\bf 117} 193-213.


\bibitem{sr}
\par\noindent\hangindent2.3em\hangafter 1
\textsc{Sz\'ekely, G. J.} and \textsc {Rizzo, M. L.} (2014). Partial
distance correlation with methods for dissimilarities. {\it Ann.
Statist.} {\bf 42} 2382-2412.


\bibitem{srb}
\par\noindent\hangindent2.3em\hangafter 1
\textsc{Sz\'ekely, G. J., Rizzo, M. L.} and \textsc{Bakirov, N. K}.
(2007). Measuring and testing independence by correlation of
distances. {\it Ann. Statist.} {\bf 35} 2769-2794.


\bibitem{wc13}
\par\noindent\hangindent2.3em\hangafter 1
\textsc {Wang, S.} and \textsc {Cui, H.} (2013). Generalized F test
for high dimensional linear regression coefficients. {\it J.
Multivariate Anal.} {\bf 117} 134-149.

\bibitem{wh07}
\par\noindent\hangindent2.3em\hangafter 1
\textsc{Wang, H.} and \textsc{He, X.} (2007). Detecting differential
expressions in geneChip microarray studies: A quantile approach.
{\it J. Amer. Statist. Assoc.} {\bf 102} 104-112.

\bibitem{wh08}
\par\noindent\hangindent2.3em\hangafter 1
\textsc{Wang, H.} and \textsc{He, X.} (2008). An enhanced quantile
approach for assessing differential gene expressions. {\it
Biometrics}. {\bf 64} 449-457.

\bibitem{wwl12}
\par\noindent\hangindent2.3em\hangafter 1
\textsc{Wang, L.}, \textsc{Wu, Y.} and \textsc{Li, R.} (2012).
Quantile regression of analyzing heterogeneity in ultra-high
dimension. {\it J. Amer. Statist. Assoc.} {\bf 107} 214-222.


\bibitem{wy93}
\par\noindent\hangindent2.3em\hangafter 1
\textsc{Westfall, P. H.} and \textsc{Young, S. S.} (1993). {\it Resampling-based Multiple Testing: Examples and Methods for p-value Adjustment}. John Wiley and Sons,
New York.

\bibitem{wh09}
\par\noindent\hangindent2.3em\hangafter 1
\textsc{Wu, C. F. J.} and \textsc{Hamada, M. S.} (2009). {\it
Experiments Planning, Analysis and Optimization.} Wiley, New York.

\bibitem{ya13}
\par\noindent\hangindent2.3em\hangafter 1
\textsc{Yata, K.} and \textsc{Aoshima, M.} (2013). Correlation tests
for high-dimensional data using extended cross-data-matrix
methodology. {\it J. Multivariate Anal.} {\bf 117} 313-331.


\bibitem{zc11}
\par\noindent\hangindent2.3em\hangafter 1
\textsc{Zhong, P. S.} and \textsc{Chen, S. X.} (2011). Testing for
high-dimensional regression coefficients with factorial designs.
{\it J. Amer. Statist. Assoc.} {\bf 106} 260-274.
\end{thebibliography}
\end{document}

\title{\Large
  Supplement to ``Conditional Mean and Quantile Dependence Testing in High Dimension''}
\author{\large Xianyang Zhang, Shun Yao, and Xiaofeng
Shao\thanks{Xianyang Zhang is Assistant Professor in the Department
of Statistics, Texas A\&M University, College Station, TX 77843,
USA. E-mail: zhangxiany@stat.tamu.edu. Shun Yao is a Ph.D. candidate
and Xiaofeng Shao is Professor in the Department of
Statistics, University of Illinois at Urbana-Champaign, Champaign,
IL 61820, USA. E-mails: shunyao2,xshao@illinois.edu. Zhang's research
is partially supported by NSF DMS-1607320. Shao's research
is partially supported by NSF-DMS1104545, NSF-DMS1407037 and DMS-1607489. We
would like to thank Dr. Ping-Shou Zhong for providing the data set
and the R codes used for the tests in the models with fractional design. We are also grateful to
Dr. Min Qian for providing the matlab code for the methods developed in McKeague and Qian  (2015).  }
\medskip\\
{\normalsize \it Texas A\&M University and University of
Illinois at Urbana-Champaign}}
\date{\normalsize This version: \today}
\maketitle
\sloppy%
\onehalfspacing
\small

\section{Technical Appendix}\label{sec:appendix}
Throughout the technical appendix, we let $c,c',c'',C,C',C'',c_i,C_i$ be generic constants which can be different from line to line.
\subsection{Unbiasedness of $MDD_n(\Y|\X)^2$}\label{sec:proof-unbias}
By Lemma 1 of Park et al. (2014), we have
$n(n-3)(\widetilde{A}\cdot\widetilde{B})=\sum_{i\neq
j}\widetilde{A}_{ij}\widetilde{B}_{ij}=\sum_{i\neq
j}A_{ij}\widetilde{\widetilde{B}}_{ij}.$ Using the fact that
$B_{ii}=0$, it can be verified that
$\widetilde{\widetilde{B}}_{ij}=\widetilde{B}_{ij}$. By the
definition of $\mathcal{U}$-centering, we have
\begin{align*}
\sum_{i\neq j}A_{ij}\widetilde{B}_{ij}=&\sum_{i\neq
j}A_{ij}\left(B_{ij}-\frac{1}{n-2}\sum^{n}_{l=1}B_{il}-\frac{1}{n-2}\sum^{n}_{k=1}B_{kj}+\frac{1}{(n-1)(n-2)}\sum_{k,l=1}^{n}B_{kl}\right)
\\=&\sum_{i\neq
j}A_{ij}B_{ij}-\frac{1}{n-2}\sum_{i\neq
j}A_{ij}\sum^{n}_{l=1}B_{il}-\frac{1}{n-2}\sum_{i\neq
j}A_{ij}\sum^{n}_{k=1}B_{kj}
\\&+\frac{1}{(n-1)(n-2)}\sum_{i,j=1}^nA_{ij}\sum_{k,l=1}^{n}B_{kl}
\\=&\text{tr}(AB)+\frac{\mathbf{1}^T_pA\mathbf{1}_p\mathbf{1}^T_pB\mathbf{1}_p}{(n-1)(n-2)}-\frac{2\mathbf{1}^T_pAB\mathbf{1}_p}{(n-2)},
\end{align*}
where $\mathbf{1}_p\in \mathbb{R}^p$ is the vector of all ones. Thus
we have
$$MDD_n(\Y|\X)^2=\frac{1}{n(n-3)}\left(\text{tr}(AB)+\frac{\mathbf{1}^T_pA\mathbf{1}_p\mathbf{1}^T_pB\mathbf{1}_p}{(n-1)(n-2)}-\frac{2\mathbf{1}^T_pAB\mathbf{1}_p}{(n-2)}\right).$$

Let $(n)_k=n!/(n-k)!$ and $I^{n}_k$ be the collections of $k$-tuples
of indices (chosen from $\{1,2,\dots,n\}$) such that each index
occurs exactly once. It can be shown that
\begin{align*}
&(n)_2^{-1}\E\left[\sum_{(i,j)\in I^n_2}A_{ij}B_{ij}\right]=(n)_2^{-1}\E[\text{tr}(AB)]=\E[K(\X,\X')L(\Y,\Y')],\\
&(n)_4^{-1}\E\left[\sum_{(i,j,q,r)\in I^n_4}A_{ij}B_{qr}\right]=(n)_4^{-1}\E[\mathbf{1}^T_pA\mathbf{1}_p\mathbf{1}^T_pB\mathbf{1}_p-4\mathbf{1}^T_pAB\mathbf{1}_p+2\text{tr}(AB)]=\E[K(\X,\X')]\E[L(\Y,\Y')],\\
&(n)_3^{-1}\E\left[\sum_{(i,j,r)\in
I^n_3}A_{ij}B_{ir}\right]=(n)_3^{-1}\E
[\mathbf{1}^T_pAB\mathbf{1}_p-\text{tr}(AB)]=\E[K(\X,\X')L(\Y,\Y'')],\\
&MDD_n(\Y|\X)^2=(n)_2^{-1}\sum_{(i,j)\in
I^n_2}A_{ij}B_{ij}+(n)_4^{-1}\sum_{(i,j,q,r)\in
I^n_4}A_{ij}B_{qr}-2(n)_3^{-1}\sum_{(i,j,r)\in I^n_3}A_{ij}B_{ir},
\end{align*}
which implies the unbiasedness of $MDD_n(\Y|\X)^2$. The above
derivation also indicates that
\begin{align*}
MDD_n(\Y|\X)^2=&(n)_4^{-1}\sum_{(i,j,q,r)\in
I^n_4}(A_{ij}B_{qr}+A_{ij}B_{ij}-2A_{ij}B_{ir})=\frac{1}{\binom{n}{4}}\sum_{i<j<q<r}h(\Z_i,\Z_j,\Z_q,\Z_r),
\end{align*}
where
\begin{equation}\label{eq:h}
\begin{split}
h(\Z_i,\Z_j,\Z_q,\Z_r)=&\frac{1}{4!}\sum_{(s,t,u,v)}^{(i,j,q,r)}(A_{st}B_{uv}+A_{st}B_{st}-2A_{st}B_{su})
\\=&\frac{1}{6}\sum_{s<t,u<v}^{(i,j,q,r)}(A_{st}B_{uv}+A_{st}B_{st})-\frac{1}{12}\sum_{(s,t,u)}^{(i,j,q,r)}A_{st}B_{su},
\end{split}
\end{equation}
with $\Z_i=(\X_i,\Y_i)$, and the summation is over all permutation
of the 4-tuples of indices $(i,j,q,r).$ Therefore $MDD_n(\Y|\X)^2$
is a $U$-statistic of order four. From the above arguments, we have
$$\sum^{p}_{k=1}MDD_n(Y|x_k)^2=\frac{1}{\binom{n}{4}}\sum_{i<j<q<r}\sum_{k=1}^p\mathfrak{h}_k(Z_{ik},Z_{jk},Z_{qk},Z_{rk}),$$
where $Z_{ik}=(x_{ik},Y_i)$ and $\mathfrak{h}_k$ is defined in a
similar way as $h$ by replacing $A_{st}$ with $A_{st}(k).$ Thus
$\sum^{p}_{k=1}MDD_n(Y|x_k)^2$ is also a fourth order $U$-statistic
with the kernel $\sum^{p}_{k=1}\mathfrak{h}_k.$

\subsection{Hoeffding decomposition}\label{sec:proof-hoef}
Recall the definition of $h$ in (\ref{eq:h}).
Define $h_c(w_1,\dots,w_c)=\E
h(w_1,\dots,w_c,\Z_{c+1},\cdots,\Z_4)$, where
$\Z_i=(\X_i,\Y_i)\overset{D}{=}(\X,\Y)$ and $c=1,2,3,4.$ The symbol ``$\overset{D}{=}$'' here means ``equal in distribution''.
Let $w=(x,y)$,
$w'=(x',y')$, $w''=(x'',y'')$ and $w'''=(x''',y''')$, where
$x,x',x'',x'''\in \mathbb{R}^q$ and $y,y',y'',y'''\in\mathbb{R}$. Further let $\Z'=(\X',\Y')$,
$\Z''=(\X'',\Y'')$ and $\Z'''=(\X''',\Y''')$ be independent copies
of $\Z=(\X,\Y)$. With some abuse of notation, we define
$U(x,x')=\E [K(x,\X')]+\E [K(\X,x')]-K(x,x')-\E [K(\X,\X')]$ and
$V(y,y')=(y-\mu_{\Y})(y'-\mu_{\Y})$ with $\mu_{\Y}=\E \Y$. Direct
calculation shows that
\begin{align*}
h_1(w)%=&\frac{1}{2}\bigg\{\E [K(x,\X)L(\Y',\Y'')]+\E [K(\X',\X'')L(y,\Y)]+\E [K(\X,\X')L(\Y,\Y')]+\E [K(x,\X)L(y,\Y)]
%\\&-\E [K(x,\X)L(y,\Y')]-\E [K(x,\X)L(\Y,\Y')]-\E [K(\X,\X')L(\Y,y)]-\E [K(\X,\X')L(\Y,\Y'')]\bigg\}
%\\=&\frac{1}{2}\bigg\{\E [K(\X,x)(L(\Y',\Y'')+L(y,\Y)-L(y,\Y')-L(\Y,\Y''))]
%\\&+\E [K(\X,\X')(L(y,\Y'')+L(\Y,\Y')-L(\Y,y)-L(\Y',\Y''))]\bigg\}
%\\=&-\frac{1}{2}\bigg\{\E [K(x,\X)(V(\Y',\Y'')+V(y,\Y)-V(y,\Y')-V(\Y,\Y''))]
%\\&+\E [K(\X,\X')(V(y,\Y'')+V(\Y,\Y')-V(\Y,y)-V(\Y',\Y''))]\bigg\}
%\\=&-\frac{1}{2}\bigg\{\E [K(x,\X)V(y,\Y)]+\E [K(\X,\X')V(\Y,\Y')]-\E [K(\X,\X')V(\Y,y)]\bigg\}
=&\frac{1}{2}\bigg\{\E [U(x,\X)V(y,\Y)]+MDD(\Y|\X)^2\bigg\},
\end{align*}
and
\begin{align*}
h_2(w,w')%=&\frac{1}{6}\bigg\{\E [K(x,x')L(\Y,\Y')]+\E [K(\X,\X')L(y,y')]+2\E [K(x,\X)L(y',\Y')]+2\E [K(x',\X)L(y,\Y')]
%\\&+K(x,x')L(y,y')+\E [K(\X,\X')L(\Y,\Y')]+2\E [K(x,\X)L(y,\Y)]+2\E [K(x',\X)L(y',\Y)]
%\\&-\E [K(x,x')L(y,\Y)]-\E [K(x',x)L(y',\Y)]-\E [K(x,\X)L(y,\Y')]-\E [K(x,\X)L(y,y')]
%\\&-\E [K(\X,x)L(\Y,\Y')]-\E [K(\X,x)L(\Y,y')]-\E [K(x',\X)L(y',\Y')]-\E [K(x',\X)L(y',y)]
%\\&-\E [K(\X,x')L(\Y,\Y')]-\E [K(\X,x')L(\Y,y)]-\E [K(\X,\X')L(\Y,y)]-\E
%[K(\X,\X')L(\Y,y')]\bigg\}
%\\=&-\frac{1}{6}\E\bigg\{K(x,x')V(y,y')+K(\X,\X')\left(V(y,y')+V(\Y,\Y')-V(y,\Y')-V(\Y,y')\right)
%\\&+K(x,\X)\left(2V(y,\Y)-V(y,y')-V(\Y,y')\right)+K(x',\X)\left(2V(y',\Y)-V(y,y')-V(\Y,y)\right) \bigg\}
=&\frac{1}{6}\bigg\{U(x,x')V(y,y')+MDD(\Y|\X)^2+\E
[U(x,\X)V(y,\Y)]+\E [U(x',\X)V(y',\Y)]
\\&+\E[(U(x,\X)-U(x',\X))(V(y,\Y)-V(y',\Y))]\bigg\}.
\end{align*}
Moreover, we have
%\begin{align*}
%&h_3(w,w',w'')
%\\=&\frac{1}{6}\E\bigg\{K(x,x')(L(y'',\Y)+L(y,y'))+K(x,x'')(L(y',\Y)+L(y,y''))+K(x,\X)(L(y',y'')+L(y,\Y))
%\\&+K(x',x'')(L(y,\Y)+L(y',y''))+K(x',\X)(L(y,y'')+L(y',\Y))+K(x'',\X)(L(y,y')+L(y'',\Y))\bigg\}
%\\&-\frac{1}{12}\E\bigg\{K(x,x')L(y,y'')+K(x,x'')L(y,y')+K(x,x')L(y,\Y)+K(x,\X)L(y,y')+K(x,x'')L(y,\Y)
%\\&+K(x,\X)L(y,y'')+K(x',x)L(y',y'')+K(x',x'')L(y',y)+K(x',x)L(y',\Y)+K(x',\X)L(y',y)
%\\&+K(x',x'')L(y',\Y)+K(x',\X)L(y',y'')+K(x'',x)L(y'',y')+K(x'',x')L(y'',y)+K(x'',x)L(y'',\Y)
%\\&+K(x'',\X)L(y'',y)+K(x'',x')L(y'',\Y)+K(x'',\X)L(y'',y')+K(\X,x)L(\Y,y')+K(\X,x')L(\Y,y)
%\\&+K(\X,x)L(\Y,y'')+K(\X,x'')L(\Y,y)+K(\X,x')L(\Y,y'')+K(\X,x'')L(\Y,y')\bigg\}
%\\=&\frac{1}{12}\E\bigg\{K(x,x')(V(y,y'')+V(y',y'')-2V(y,y'))+K(x,x'')(V(y,y')+V(y'',y')-2V(y,y''))
%\\&+K(x',x'')(V(y',y)+V(y'',y)-2V(y',y''))+K(x,\X)(V(y,y')+V(y,y'')-2V(y',y''))
%\\&+K(\X,x')(V(y',y)+V(y',y'')-2V(y,y''))+K(x'',\X)(V(y,y'')+V(y',y'')-2V(y,y'))\bigg\}
%\\=&\frac{1}{12}\E\bigg\{(K(x,x')+K(x'',\X))(V(y,y'')+V(y',y'')-2V(y,y'))
%\\&+(K(x,x'')+K(x',\X))(V(y,y')+V(y'',y')-2V(y,y''))
%\\&+(K(x',x'')+K(x,\X))(V(y',y)+V(y'',y)-2V(y',y''))
%\\&+K(x'',\X)(V(y,\Y)+V(y',\Y)-2V(y'',\Y))+K(x',\X)(V(y,\Y)+V(y'',\Y)-2V(y',\Y))
%\\&+K(x,\X)(V(y',\Y)+V(y'',\Y)-2V(y,\Y))\bigg\}.
%\end{align*}
%Rearranging the above terms, we deduce that
\begin{align*}
h_3(w,w',w'')=&\frac{1}{12}\bigg\{(2U(x,x')-U(x',x'')-U(x,x''))V(y,y')
\\ &+(2U(x,x'')-U(x,x')-U(x',x''))V(y,y'')
\\ &+(2U(x',x'')-U(x,x')-U(x,x''))V(y',y'')
\\ &+\E[(2U(x,\X)-U(x',\X)-U(x'',\X))V(y,\Y)]
\\ &+\E[(2U(x',\X)-U(x,\X)-U(x'',\X))V(y',\Y)]
\\ &+\E[(2U(x'',\X)-U(x,\X)-U(x',\X))V(y'',\Y)]\bigg\}.
\end{align*}
Using similar calculation as that for $h_3(w,w',w'')$, we obtain
\begin{align*}
&h_4(w,w',w'',w''')
%\\=&\frac{1}{12}\bigg\{(K(x,x')+K(x'',x'''))(V(y,y'')+V(y',y'')-2V(y,y'))
%\\&+(K(x,x'')+K(x',x'''))(V(y,y')+V(y'',y')-2V(y,y''))
%\\&+(K(x',x'')+K(x,x'''))(V(y',y)+V(y'',y)-2V(y',y''))
%\\&+(K(x,x')+K(x'',x'''))(V(y,y''')+V(y',y''')-2V(y'',y'''))
%\\&+(K(x,x'')+K(x',x'''))(V(y,y''')+V(y'',y''')-2V(y',y'''))
%\\&+(K(x',x'')+K(x,x'''))(V(y',y''')+V(y'',y''')-2V(y,y''')))\bigg\}
\\=&\frac{1}{12}\bigg\{(2U(x,x')+2U(x'',x''')-U(x,x'')-U(x,x''')-U(x',x'')-U(x',x'''))(V(y,y')+V(y'',y'''))
\\&+(2U(x,x'')+2U(x',x''')-U(x,x')-U(x,x''')-U(x'',x')-U(x'',x'''))(V(y,y'')+V(y',y'''))
\\&+(2U(x,x''')+2U(x'',x')-U(x,x'')-U(x,x')-U(x''',x'')-U(x''',x'))(V(y,y''')+V(y',y''))\bigg\}.
\end{align*}

\emph{Analysis under the null hypothesis} When $MDD(\Y|\X)^2=0$, it
can be verified that $h_1(w)=0$. In this case, we have
$h_2(w,w')=U(x,x')V(y,y')/6$ and
\begin{align*}
h_3(w,w',w'')=&\frac{1}{12}\bigg\{(2U(x,x')-U(x',x'')-U(x,x''))V(y,y')
\\ &+(2U(x,x'')-U(x,x')-U(x',x''))V(y,y'')
\\ &+(2U(x',x'')-U(x,x')-U(x,x''))V(y',y'')
\bigg\}.
\end{align*}
It is not hard to verify that under the null
\begin{align*}
\var(h_2(\Z,\Z'))=\frac{1}{36}\E [U^2(\X,\X')V^2(\Y,\Y')]=\frac{1}{36}\xi^2,
\end{align*}
and
\begin{align*}
\var(h_3(\Z,\Z',\Z''))=&\frac{3}{144}\var\{(2U(\X,\X')-U(\X',\X'')-U(\X,\X''))V(\Y,\Y')\}
\\=&\frac{3}{144}\bigg\{4\xi^2+2\E [U(\X,\X'')^2V(\Y,\Y')^2]+2\E
[U(\X,\X'')U(\X',\X'')V(\Y,\Y')^2]\bigg\},
\end{align*}
where $\xi^2=\var(U(X,X')V(Y,Y')).$ Furthermore, careful calculation
yields that
\begin{align*}
\var(h_4(\Z,\Z',\Z'',\Z'''))=&\frac{6}{144}\E\bigg\{
V(\Y,\Y')^2[U(\X,\X'')+U(\X',\X''')+U(\X',\X'')
\\&+U(\X,\X''')-2U(\X,\X')-2U(\X'',\X''')]^2\bigg\}
\\=&\frac{1}{6}\{\E [V(\Y,\Y')^2U(\X,\X'')U(\X',\X'')]+\E [V(\Y,\Y')^2U(\X,\X'')^2]
\\&+\E [V(\Y,\Y')^2]\E [U(\X,\X')^2]+\xi^2\}.
\end{align*}
By the Cauchy-Schwarz inequality, we have
\begin{align*}
\E [U(\X,\X'')U(\X',\X'')V(\Y,\Y')^2]\leq& \{\E [U(\X,\X'')^2V(\Y,\Y')^2]\}^{1/2} \{\E [U(\X',\X'')^2V(\Y,\Y')^2]\}^{1/2}
\\=&\E [U(\X,\X'')^2V(\Y,\Y')^2].
\end{align*}
Under the assumption that
\begin{align*}
&\frac{\E [U(\X,\X'')^2V(\Y,\Y')^2]}{\xi^2}=o(n),
\\&\frac{\E [V(\Y,\Y')^2]\E[U(\X,\X')^2]}{\xi^2}=o(n^2),
\end{align*}
we have
\begin{align}\label{eq:u-approx}
MDD_n(\Y|\X)^2=\frac{1}{\binom{n}{2}}\sum_{1\leq i<j\leq
n}U(\X_i,\X_j)V(\Y_i,\Y_j)+\mathcal{R}_n,
\end{align}
where $\mathcal{R}_n$ is the remainder term which is asymptotically
negligible [see Serfling (1980)].

\emph{Analysis under local alternatives} We consider the case where
$MDD(\Y|\X)^2$ is nonzero, i.e., the conditional mean of $\Y$
may depend on $\X$. Recall that $\mathcal{L}(x,y)=\E
[U(x,\X)V(y,\Y)]$. Under the assumption that
\begin{align}\label{eq-add}
\var(\mathcal{L}(\X,\Y))=o(n^{-1}\xi^2),\quad
\var(\mathcal{L}(\X,\Y'))=o(\xi^2),
\end{align}
we get
\begin{align*}
\var(h_1)=o(n^{-1}\xi^2),\quad \var(h_2)=\frac{\xi^2}{36}(1+o(1)).
\end{align*}
Moreover, we have
\begin{align*}
\var(h_3(\Z,\Z',\Z''))\leq C\bigg\{\xi^2+\E
[U(\X,\X'')^2V(\Y,\Y')^2]+\E [U(\X,\X'')U(\X',\X'')V(\Y,\Y')^2]\bigg\},
\end{align*}
and
\begin{align*}
\var(h_4(\Z,\Z',\Z'',\Z'''))\leq &C'\{\E
[V(\Y,\Y')^2U(\X,\X'')U(\X',\X'')]+\E [V(\Y,\Y')^2U(\X,\X'')^2]
\\&+\E [V(\Y,\Y')^2]\E [U(\X,\X')^2]+\xi^2\}.
\end{align*}
Thus (\ref{eq:u-approx}) still holds under Assumption (\ref{eq-add}). Applying the above arguments to
$\sum^{p}_{k=1}MDD_n(Y|x_k)^2$, we deduce that
\begin{align*}
\sum^{p}_{k=1}\{MDD_n(Y|x_k)^2-MDD(Y|x_k)^2\}=\frac{1}{\binom{n}{2}}\sum_{1\leq
i<j\leq n}\widetilde{U}(X_i,X_j)V(Y_i,Y_j)+\sum^{p}_{k=1}R_{k,n}.
\end{align*}
The remainder term $\sum^{p}_{k=1}R_{k,n}$ is
asymptotically negligible provided that
\begin{align*}
&\frac{\E [\widetilde{U}(X,X'')^2V(Y,Y')^2]}{\mathcal{S}^2}=o(n),
%\\&\frac{\E \widetilde{U}(X,X'')\widetilde{U}(X',X'')V(Y,Y')^2}{\mathcal{S}^2}=O(1),
\\&\frac{\E [\widetilde{U}(X,X')^2]\E
[V(Y,Y')^2]}{\mathcal{S}^2}=o(n^2),\\
&\var(\widetilde{\mathcal{L}}(X,Y))=o(n^{-1}\mathcal{S}^2),\quad
\var(\widetilde{\mathcal{L}}(X,Y'))=o(\mathcal{S}^2),
\end{align*}
where $\widetilde{\mathcal{L}}(x,y)=\E [\widetilde{U}(x,X)V(y,Y)]$.

\subsection{Asymptotic normality under the null and alternatives}\label{sec:normal}
We shall establish the asymptotic normality for $\breve{T}_n$ using
the Central Limit Theorem for martingale difference sequences. We
first restrict our analysis under $H_0'$. Let
$$S_{r}:=\sum_{j=2}^r\sum^{j-1}_{i=1}\widetilde{U}(X_i,X_j)V(Y_i,Y_j)=\sum_{j=2}^r\sum^{j-1}_{i=1}H(Z_i,Z_j).$$
Define the filtration $\mathcal{F}_{r}=\sigma\{Z_1,Z_2,\dots,Z_r\}$
with $Z_i=(X_i,Y_i)$. It is not hard to see that $S_{r}$ is adaptive
to $\mathcal{F}_r$ and $S_r$ is a mean-zero martingale sequence,
i.e., $\E[S_r]=0$ and
$$\E[S_{r'}|
\mathcal{F}_r]=S_r+\sum_{j=r+1}^{r'}\sum^{j-1}_{i=1}\E[\E[\widetilde{U}(X_i,X_j)V(Y_i,Y_j)|\mathcal{F}_r,X_i,X_j]|\mathcal{F}_r]=S_r$$
for $r'\geq r.$ Thus by verifying the following two conditions [see
e.g. Lemmas 2-3 of Chen and Qin (2010)], we can establish the
asymptotic normality for $\breve{T}_{n}$ by Corollary 3.1 of Hall
and Heyde (1980). Specifically, define
$\mathcal{W}_j=\sum^{j-1}_{i=1}H(Z_i,Z_j)$. We need to show that
\begin{align}\label{eq:B}
\sum^{n}_{j=1}B^{-2}\E[\mathcal{W}_j^2\mathbf{I}\{|\mathcal{W}_j|>\epsilon
B\}|\mathcal{F}_{j-1}]\rightarrow^p 0,
\end{align}
for $B$ such that
\begin{align}\label{eq:B2}
\sum^{n}_{j=1}\E[\mathcal{W}_j^2|\mathcal{F}_{j-1}]/B^2\rightarrow^p
C>0.
\end{align}
We shall first prove that
\begin{equation}\label{eq:B3}
\frac{2}{n(n-1)\mathcal{S}^2}\sum^{n}_{j=1}\E[\mathcal{W}_j^2|\mathcal{F}_{j-1}]\rightarrow^p 1,
\end{equation}
i.e., $B^2=n(n-1)\mathcal{S}^2/2$ and $C=1$ in (\ref{eq:B2}). Notice that
\begin{align*}
\E[\mathcal{W}_j^2|\mathcal{F}_{j-1}]=&\E\left[\sum^{j-1}_{i,k=1}H(Z_i,Z_j)H(Z_k,Z_j)\bigg|\mathcal{F}_{j-1}\right]=\sum^{j-1}_{i,k=1}G(Z_i,Z_k),
\end{align*}
and
\begin{align*}
\frac{2}{n(n-1)}\sum^{n}_{j=2}\E[\mathcal{W}_j^2]=&\frac{2}{n(n-1)}\sum^{n}_{j=2}\E\left[\sum^{j-1}_{i,k=1}(Y_i-\mu)(Y_k-\mu)(Y_j-\mu)^2\widetilde{U}(X_i,X_j)\widetilde{U}(X_k,X_j)\right]
\\=&\frac{2}{n(n-1)}\sum^{n}_{j=2}\E\left[\sum^{j-1}_{i=1}(Y_i-\mu)^2(Y_j-\mu)^2\widetilde{U}(X_i,X_j)^2\right]
\\=&\E\left[V(Y,Y')^2\widetilde{U}(X,X')^2\right]=E[H(Z,Z')^2]=\mathcal{S}^2.
\end{align*}
Define the following quantities
\begin{align*}
&\mathcal{D}_1=\E
[H(Z,Z'')^2H(Z',Z'')^2]-(\E[H(Z,Z')^2])^2=\var(G(Z,Z)),\\
&\mathcal{D}_2=\E[H(Z,Z')H(Z',Z'')H(Z'',Z''')H(Z''',Z)]=E[G(Z,Z')^2].
\end{align*}
We have for $j\geq j'$
\begin{align*}
\cov\left(\E[\mathcal{W}_j^2|\mathcal{F}_{j-1}],\E[\mathcal{W}_{j'}^2|\mathcal{F}_{j'-1}]\right)=&\sum^{j-1}_{i,k=1}\sum^{j'-1}_{i',k'=1}\cov(G(Z_i,Z_k),G(Z_{i'},Z_{k'}))
\\=&(j'-1)\mathcal{D}_1+2(j'-1)(j'-2)\mathcal{D}_2.
\end{align*}
Under the assumption that
$$\frac{\E [G(Z,Z')^2]}{\{\E
[H(Z,Z')^2]\}^2}\rightarrow 0, \quad \frac{\E
[H(Z,Z'')^2H(Z',Z'')^2]}{n\{\E [H(Z,Z')^2]\}^2}\rightarrow 0,$$ we
have
\begin{align*}
\frac{4}{n^2(n-1)^2}\sum^{n}_{j,j'=2}\cov\left(\E[\mathcal{W}_j^2|\mathcal{F}_{j-1}],\E[\mathcal{W}_{j'}^2|\mathcal{F}_{j'-1}]\right)=O(\mathcal{D}_1/n+\mathcal{D}_2)=o(\mathcal{S}^4),
\end{align*}
which ensures (\ref{eq:B3}). To show (\ref{eq:B}), we note
that
\begin{equation}\label{eq-1}
\sum^{n}_{j=1}B^{-2}\E[\mathcal{W}_j^2\mathbf{I}\{|\mathcal{W}_j|>\epsilon
B\}|\mathcal{F}_{j-1}]\leq
B^{-2-s}\epsilon^{-s}\sum^{n}_{j=1}\E[|\mathcal{W}_j|^{2+s}|\mathcal{F}_{j-1}],
\end{equation}
for some $s>0.$ With $s=2$, we prove that
$$B^{-4}\sum^{n}_{j=1}\E[|\mathcal{W}_j|^{4}|\mathcal{F}_{j-1}]\rightarrow^p 0,$$
where $B^2=n(n-1)\mathcal{S}^2/2$. To this end, it suffices to show
that
$$B^{-4}\sum^{n}_{j=1}\E[|\mathcal{W}_j|^{4}]\rightarrow 0.$$
Under the assumption
$$\frac{\E [H(Z,Z')^4]/n+\E [H(Z,Z'')^2H(Z',Z'')^2]}{n\{\E
[H(Z,Z')^2]\}^2}\rightarrow 0,$$ we have
\begin{align*}
\sum^{n}_{j=2}\E[|\mathcal{W}_j|^{4}]=&\sum^{n}_{j=2}\sum^{j-1}_{i_1,i_2,i_3,i_4=1}\E [H(Z_{i_1},Z_j)H(Z_{i_2},Z_j)H(Z_{i_3},Z_j)H(Z_{i_4},Z_j)]
\\=&\frac{n(n-1)}{2}\E [H(Z,Z')^4]+3\sum^{n}_{j=2}\sum_{i_1\neq i_2}\E [H(Z_{i_1},Z_j)^2H(Z_{i_2},Z_j)^2]
\\=&\frac{n(n-1)}{2}\E [H(Z,Z')^4]+3\sum^{n}_{j=2}(j-1)(j-2)\E [H(Z,Z'')^2H(Z',Z'')^2]=o(B^4),
\end{align*}
which implies (\ref{eq:B}).

We next extend the above arguments to local alternatives. Under the
assumption that
$\var(\widetilde{\mathcal{L}}(X,Y))=o(n^{-1}\mathcal{S}^2)$, we
have
\begin{align*}
&\frac{1}{\sqrt{\binom{n}{2}}\mathcal{S}}\sum_{1\leq i<j\leq
n}\left(\widetilde{U}(X_i,X_j)V(Y_i,Y_j)-\E[\widetilde{U}(X_i,X_j)V(Y_i,Y_j)]\right)
\\=&\frac{1}{\sqrt{\binom{n}{2}}\mathcal{S}}\sum_{1\leq i<j\leq
n}H^*(Z_i,Z_j)+\frac{1}{\sqrt{\binom{n}{2}}\mathcal{S}}\sum_{1\leq
i<j\leq
n}\left(\widetilde{\mathcal{L}}(X_i,Y_i)+\widetilde{\mathcal{L}}(X_j,Y_j)-2\E[\widetilde{U}(X_i,X_j)V(Y_i,Y_j)]\right)
\\=&\frac{1}{\sqrt{\binom{n}{2}}\mathcal{S}}\sum_{1\leq i<j\leq
n}H^*(Z_i,Z_j)+o_p(1).
\end{align*}
Using similar arguments by replacing $H$ with $H^*$, we can show
that
$$\frac{1}{\sqrt{\binom{n}{2}}\mathcal{S}}\sum_{1\leq i<j\leq
n}H^*(Z_i,Z_j)\rightarrow^d N(0,1),$$
which implies that
$$\frac{1}{\sqrt{\binom{n}{2}}\mathcal{S}}\sum_{1\leq i<j\leq
n}\left(\widetilde{U}(X_i,X_j)V(Y_i,Y_j)-\E[\widetilde{U}(X_i,X_j)V(Y_i,Y_j)]\right)\rightarrow^d
N(0,1).$$ The proof is thus completed.

\subsection{Interpretation of Condition (13) based on Fourier embedding}\label{sec:proof-eigen}
In this subsection, we provide more discussions on Condition
(13) in the paper. Consider a generic random variable
$\Z=(\X,\Y)$ and a random sample $\{\Z_i\}^{n}_{i=1}$ from the
distribution of $\Z$. Let $\I=\sqrt{-1}$ and
$c_q=\pi^{(1+q)/2}/\Gamma((1+q)/2)$. For $a_1,a_2\in\mathbb{R}^q$,
define the inner product $<a_1,a_2>=a_1^Ta_2$. Denote by
$\overline{a}$ the conjugate of a complex number $a$. First note
that an alternative expression for MDD is given by
\begin{align*}
MDD(\Y|\X)^2=\frac{1}{c_q}\int_{\mathbb{R}^q}\frac{|\E[\Y
e^{\I<t,\X>}]-\mu_{\Y}f_{\X}(t)|}{|t|_q^{1+q}}dt,
\end{align*}
where $f_{\X}(t)=\E[e^{\I<t,\X>}]$ and $\mu_{\Y}=\E[\Y]$. Recall that $U(x,x')=\E [K(x,\X')]+\E [K(\X,x')]-K(x,x')-\E [K(\X,\X')]$ and
$V(y,y')=(y-\mu_{\Y})(y'-\mu_{\Y})$ with $\mu_{\Y}=\E \Y$ and $\Z'=(\X',\Y')\overset{D}{=}(\X,\Y)$. By the
analysis in Sections \ref{sec:proof-unbias} and \ref{sec:proof-hoef},
we know $MDD(\Y|\X)^2$ is a fourth-order $U$-statistic, whose major
term is given by
\begin{equation}
\mathcal{V}_n=\frac{1}{n(n-1)}\sum_{1\leq i\neq j\leq
n}U(\X_i,\X_j)V(\Y_i,\Y_j).
\end{equation}
By Lemma 1 of Sz\'{e}kely et al. (2007), we have
\begin{align*}
U(x,x')=&\int_{\mathbb{R}^q}
\E(e^{\I<t,\X>}-e^{\I<t,x>})(\overline{e^{\I<t,\X'>}-e^{\I<t,x'>}})w_q(t)dt
\\=&\int_{\mathbb{R}^q}
(f_{\X}(t)-e^{\I<t,x>})(\overline{f_{\X}(t)}-e^{-\I<t,x'>})w_q(t)dt,
\end{align*}
where $w_q(t)=1/(c_q|t|_q^{1+q})$. Let $K(z,z')=U(x,x')V(y,y')$ with
$z=(x,y)$ and $z'=(x',y').$ Thus we have
\begin{align*}
K(z,z')=&\int_{\mathbb{R}^q}
(f_{\X}(t)-e^{\I<t,x>})(\mu_{\Y}-y)(\overline{f_{\X}(t)}-e^{-\I<t,x'>})(\mu_{\Y}-y')w_q(t)dt
\\=&\int_{\mathbb{R}^q} G(x,y;t)\overline{G(x',y';t)}w_q(t)dt,
\end{align*}
where $G(x,y;t)=(f_{\X}(t)-e^{\I<t,x>})(\mu_{\Y}-y).$ Note that
$\E[G(\X,\Y;t)]=\E[\Y e^{\I<t,\X>}]-\mu_{\Y}f_{\X}(t)$, and thus
$$\E [K(\Z,\Z')]=\int_{\mathbb{R}^q}|\E[\Y
e^{\I<t,\X>}]-\mu_{\Y}f_{\X}(t)|^2w_q(t)dt=MDD(\Y|\X)^2.$$ Define
the operator $A_Kg(z)=\E[K(\Z,z)g(\Z)]$ for $g:\mathbb{R}^{q+1}\rightarrow \mathbb{R}$. The eigenfunction $h$ of
$A_K$ satisfies that
\begin{align*}
\E[K(\Z,z)h(\Z)]=\int_{\mathbb{R}^q}
\E[G(\X,\Y;t)h(\X,\Y)]\overline{G(x,y;t)}w_q(t)dt=\lambda h(x,y),
\end{align*}
which implies that $h$ has the form of
$$h(x,y)=\int_{\mathbb{R}^q}
\eta(t)\overline{G(x,y;t)}w_q(t)dt.$$ Plugging back into the above
equation, we obtain,
\begin{align*}
&\int_{\mathbb{R}^q} \E[G(\X,\Y;t)h(\X,\Y)]\overline{G(x,y;t)}w_q(t)dt
\\=&\int_{\mathbb{R}^q} \left(\int_{\mathbb{R}^q}
\eta(t')\E[G(\X,\Y;t)\overline{G(\X,\Y;t')}]w_q(t')dt'\right)\overline{G(x,y;t)}w_q(t)dt
\\=&\lambda\int_{\mathbb{R}^q}
\eta(t)\overline{G(x,y;t)}w_q(t)dt.
\end{align*}
The above equation holds provided that
\begin{align}\label{eq:phi}
\lambda\eta(t)=\int_{\mathbb{R}^q}
\eta(t')\E[G(\X,\Y;t)\overline{G(\X,\Y;t')}]w_q(t')dt':=\varphi(\eta)(t),
\end{align}
where $\varphi(\cdot)$ is the corresponding operator defined in (\ref{eq:phi}).
Thus $\eta$ is the eigenfunction of $\varphi$ and $\lambda$ is the
corresponding eigenvalue. When $\X$ and $(\Y-\mu_{\Y})^2$ are
independent, we have
\begin{align*}
\E[G(\X,\Y;t)\overline{G(\X,\Y;t')}]=\var(\Y)(f_{\X}(t-t')-f_{\X}(t)f_{\X}(-t')).
\end{align*}
Suppose that
$$f_{\X}(t-t')-f_{\X}(t)f_{\X}(-t')=\sum^{+\infty}_{j=1}\lambda_{\X,j}\beta_{\X,j}(t)\overline{\beta_{\X,j}(t')},$$
where $\int_{\mathbb{R}^q}
\beta_{\X,i}(t)\overline{\beta_{\X,j}(t)}w_q(t)dt=\mathbf{1}\{i=j\}.$ In this
case, Condition (13) reduces to
\begin{align}\label{con1}
\frac{(\sum^{+\infty}_{i=1}\lambda_{\X,i}^t)^{2/t}}{\sum^{+\infty}_{i=1}\lambda_{\X,i}^2}\rightarrow
0,
\end{align}
for some $t>2.$
%The key message from the above derivation is that
%for different kernel $K(\cdot,\cdot)$, the corresponding quadratic
%form captures a different aspect of the underlying distribution. And
%the conditions we need to impose rely on both the kernel and the
%underlying distribution.

Recall that
$H(Z_i,Z_k)=V(Y_i,Y_k)\widetilde{U}(X_i,X_k)=V(Y_i,Y_k)\sum^{p}_{j=1}U_j(x_{ij},x_{kj})$.
Below we focus on the case where $K(\cdot,\cdot)=H(\cdot,\cdot)$ and
$\Z=Z$. Define $G_j(u_j,y;t)=(f_{x_j}(t)-e^{\I tu_j})(\mu-y)$ for $u_j,y,t\in\mathbb{R}$.
Then it can be shown that
\begin{align*}
H(z,z')=\int_{\mathbb{R}}\sum^{p}_{j=1}G_j(u_j,y;t)\overline{G_j(u_j',y';t)}w_1(t)dt,
\end{align*}
where $z=(u_1,\dots,u_p,y)^T$ and $z'=(u_1',\dots,u_p',y')^T$. The
eigenfunction of $A_H$ satisfies that
\begin{align*}
\E[H(Z,z)h(Z)]=\int_{\mathbb{R}}
\sum^{p}_{j=1}\E[G_j(x_j,Y;t)h(Z)]\overline{G_j(u_j,y;t)}w_1(t)dt=\lambda
h(z).
\end{align*}
Thus we must have $h(z)=\int_{\mathbb{R}}
\sum^{p}_{j=1}\eta_j(t)\overline{G_j(u_j,y;t)}w_1(t)dt.$ It implies
that
\begin{equation}\label{eq:egen2}
\begin{split}
&\int_{\mathbb{R}}\sum^{p}_{j=1}\left(\int_{\mathbb{R}}
\sum^{p}_{i=1}\eta_i(t')\E[G_j(x_j,Y;t)\overline{G_i(x_i,Y;t')}]w_1(t')dt'\right)\overline{G_j(u_j,y;t)}w_1(t)dt
\\=&\lambda\int_{\mathbb{R}}
\sum^{p}_{j=1}\eta_j(t)\overline{G_j(u_j,y;t)}w_1(t)dt.
\end{split}
\end{equation}
Let $\varphi_{ji}(g)(t)=\int_{\mathbb{R}}
g(t')\E[G_j(x_j,Y;t)\overline{G_i(x_i,Y;t')}]w_1(t')dt'$. When
$X=(x_1,\dots,x_p)^T$ and $(Y-\mu)^2$ are independent, we have
\begin{equation*}
\E[G_j(x_j,Y;t)\overline{G_i(x_i,Y;t')}]=\var(Y)(f_{x_j,x_i}(t,-t')-f_{x_j}(t)f_{x_i}(-t')),
\end{equation*}
in which case $\varphi_{ji}(g)(t)=\var(Y)\int_{\mathbb{R}}
g(t')(f_{x_j,x_i}(t,-t')-f_{x_j}(t)f_{x_i}(-t'))w_1(t')dt'.$ Notice
that (\ref{eq:egen2}) holds provided that
\begin{equation}
\sum^{p}_{i=1}\varphi_{ji}(\eta_i)=\lambda\eta_j,
\end{equation}
for all $1\leq j\leq p$. Let $L^2(w_1)$ be the space of functions
such that $\int_{\mathbb{R}} f(t)^2w_1(t)dt<\infty$ for $f\in
L^2(w_1)$. For $g=(g_1,\dots,g_p)^T$ with $g_j\in L^2(w_1)$, define
the operator
$$\Phi(g)=\left(\sum^{p}_{j=1}\varphi_{1j}(g_j),\dots,\sum^{p}_{j=1}\varphi_{pj}(g_j)\right).$$
Therefore, $\lambda$ is the eigenvalue associated with $\Phi.$ Define $Tr(\cdot)$ as the nuclear norm for a Hermitian
operator. Then with $t=4$, the condition needed to ensure the asymptotic normality [Hall (1984)] is
\begin{align*}
\frac{Tr^{1/2}(\Phi^4)}{Tr(\Phi^2)}\rightarrow 0.
\end{align*}

\begin{remark}
{\rm
We remark that the above analysis can be extended to a more general
class of dependence measures such as the distance covariance in
Sz\'{e}kely et al. (2007).
}
\end{remark}

\subsection{Conditions in Theorem 2.1}
We study the following conditions imposed in Theorem 2.1 [see (8), (10) and (11) therein],
\begin{align}
&\frac{\E[G(Z,Z')^2]}{\mathcal{S}^4}\rightarrow 0,
\label{eq:con1}\\
&\frac{\E [H(Z,Z')^4]}{n\mathcal{S}^4}\rightarrow 0, \label{eq:con2}
\\ &\frac{\E
[\widetilde{U}(X,X'')^2V(Y,Y')^2]}{\mathcal{S}^2}=o(n), \label{eq-r1}
\\&\frac{\E [\widetilde{U}(X,X')^2]\E [V(Y,Y')^2]}{\mathcal{S}^2}=o(n^2). \label{eq-r2}
\end{align}
Note that (2) is slightly stronger than the second part of (8) in
the paper. To proceed, we summarize the general conditions we
shall study in this subsection:
\begin{align}
&0<c\leq \text{var}(Y|X)=E[(Y-E[Y|X])^2|X]\leq
\E[(Y-\E[Y|X])^4|X]^{1/2}\leq C<+\infty, \label{eq-1}\\
&\frac{\E[\widetilde{U}(X,X')\widetilde{U}(X',X'')\widetilde{U}(X'',X''')\widetilde{U}(X''',X)]}{(\sum_{j,j'=1}^pdcov(x_j,x_{j'})^2)^2}\rightarrow
0,
\label{eq:con1-2}\\
&\frac{\E
[\widetilde{U}(X,X')^4]}{n(\sum_{j,j'=1}^pdcov(x_j,x_{j'})^2)^2}\rightarrow
0, \label{eq:con2-2}
\end{align}
where (\ref{eq-1}) holds almost surely for some constants $c$ and
$C.$

%\begin{remark}
%{\rm Suppose $Y=g(X)+\sigma(X)\epsilon$, where $E[\epsilon|X]=0$,
%$\E[\epsilon^2|X]\geq c_1$ and $\E[\epsilon^4|X]<C_1$. Note that
%$\E[Y|X]=g(X)$. Suppose $0<c<\sigma(X)<C<\infty,$ then it is
%straightforward to see that (\ref{eq-1}) is fulfilled.}
%\end{remark}

Under (\ref{eq-1}) and the null hypothesis (i.e. $E[Y|X]=\mu$),
\begin{align*}
\mathcal{S}^2=\E[H(Z,Z')^2]=&\E[\E[(Y-\mu)^2(Y'-\mu)^2\widetilde{U}(X,X')^2|X,X']]
\\ \geq&
c^2 \E[\widetilde{U}(X,X')^2]=c^2\sum_{j,j'=1}^pdcov(x_j,x_{j'})^2.
\end{align*}
Also note that
\begin{align*}
\E [\widetilde{U}(X,X'')^2V(Y,Y')^2]\leq& \E[\E
[\widetilde{U}(X,X'')^2V(Y,Y')^2|X,X']]
\\ \leq& C^2\E[\widetilde{U}(X,X''
)^2]\leq C^2\sum_{j,j'=1}^pdcov(x_j,x_{j'})^2.
\end{align*}
Thus (\ref{eq-r1}) and (\ref{eq-r2}) hold under (\ref{eq-1}) and the
null hypothesis. Using the conditioning argument, we have
\begin{align*}
&\E[G(Z,Z')^2]\leq C^4
\E[\widetilde{U}(X,X'')\widetilde{U}(X',X'')\widetilde{U}(X,X''')\widetilde{U}(X',X''')],\\
&\E [H(Z,Z')^4]\leq C^4 \E[\widetilde{U}(X,X')^4].
\end{align*}
Thus (\ref{eq:con1-2})-(\ref{eq:con2-2}) imply
(\ref{eq:con1})-(\ref{eq:con2}) under (\ref{eq-1}). To sum up,
(\ref{eq-1})-(\ref{eq:con2-2}) together imply
(\ref{eq:con1})-(\ref{eq-r2}).

Using the fact that
$$U_j(u_j,u_j')=\int_{\mathbb{R}}(f_j(t)-e^{\I t
u_j})(f_j(-t)-e^{-\I t u_j'})w(t)dt,\quad \I=\sqrt{-1},$$ we have
\begin{align*}
&\E[\widetilde{U}(X,X'')\widetilde{U}(X',X'')\widetilde{U}(X,X''')\widetilde{U}(X',X''')]
\\
=&\sum_{j,k,l,m=1}^p\int\{f_{jl}(t_1,t_2)-f_j(t_1)f_l(t_2)\}\{f_{lm}(-t_2,-t_4)-f_l(-t_2)f_m(-t_4)\}
\\&\{f_{mk}(t_4,t_3)-f_m(t_4)f_k(t_3)\}\{f_{kj}(-t_3,-t_1)-f_k(-t_3)f_j(-t_1)\}
\\&w(t_1)w(t_2)w(t_3)w(t_4)dt_1dt_2dt_3dt_4,
\end{align*}
where $f_{jk}$ and $f_j$ denote the (joint) characteristic functions for
$(X_j,X_k)$ and $X_j$ respectively. Similarly, we have
\begin{align*}
&\E[\widetilde{U}(X,X')^4]
\\=&
\sum_{j,k,l,m=1}^p\int \left|\E[\{f_j(t_1)-e^{\I t_1
X_j}\}\{f_k(t_2)-e^{\I t_2 X_k}\}\{f_l(t_3)-e^{\I t_3
X_l}\}\{f_m(t_4)-e^{\I t_4 X_m}\}]\right|^2
\\&w(t_1)w(t_2)w(t_3)w(t_4)dt_1dt_2dt_3dt_4,
\end{align*}
where
\begin{equation}\label{eq:ch}
\begin{split}
&\E[\{f_j(t_1)-e^{\I t_1 X_j}\}\{f_k(t_2)-e^{\I t_2
X_k}\}\{f_l(t_3)-e^{\I t_3 X_l}\}\{f_m(t_4)-e^{\I t_4 X_m}\}]
\\=&\E[\{f_j(t_1)f_k(t_2)-e^{\I t_1 X_j}f_k(t_2)-f_j(t_1)e^{\I t_2
X_k}+e^{\I t_1X_j+\I t_2 X_k}\}
\\&\{f_l(t_3)f_m(t_4)-e^{\I t_3
X_l}f_m(t_4)-f_l(t_3)e^{\I t_4 X_m}+e^{\I t_3X_l+\I t_4 X_m}\}]
\\=&-3f_j(t_1)f_k(t_2)f_l(t_3)f_m(t_4)+f_{jl}(t_1,t_3)f_k(t_2)f_m(t_4)+f_{jk}(t_1,t_2)f_l(t_3)f_m(t_4)
\\&+f_{jm}(t_1,t_4)f_k(t_2)f_l(t_3)+f_{km}(t_2,t_4)f_j(t_1)f_l(t_3)+f_{kl}(t_2,t_3)f_j(t_1)f_m(t_4)
\\&+f_{lm}(t_3,t_4)f_j(t_1)f_k(t_2)-f_{jlm}(t_1,t_3,t_4)f_k(t_2)-f_{klm}(t_2,t_3,t_4)f_j(t_1)-f_{jkl}(t_1,t_2,t_3)f_m(t_4)
\\&-f_{jkm}(t_1,t_2,t_4)f_l(t_3)+f_{jklm}(t_1,t_2,t_3,t_4).
\end{split}
\end{equation}
Hence (\ref{eq:con1-2})-(\ref{eq:con2-2}) can be translated into
conditions on the joint characteristic functions.

\textbf{Banded dependence structure:} For the ease of notation, we
assume that $\pi(i)=i.$ The result can be extended to the case where
$\pi(\cdot)$ is an arbitrary permutation of $\{1,2\dots,p\}$. In
this case, we know that $X_i$ and $X_{j}$ are independent provided
that $|i-j|>L$. By Lemma \ref{lemma1} below, we have
\begin{align*}
&\E[\widetilde{U}(X,X')\widetilde{U}(X',X'')\widetilde{U}(X'',X''')\widetilde{U}(X''',X)]
\\ = & \sum_{j=1}^{p}\sum_{k=j-L}^{j+L}\sum_{l=k-L}^{k+L}\sum_{m=j-L}^{j+L}\E[U_j(x_j,x'_j)U_k(x'_k,x''_k)U_l(x_l'',x'''_l)U_m(x'''_m,x_m)]
\\ \leq &
C\sum_{j=1}^{p}\sum_{k=j-L}^{j+L}\sum_{l=k-L}^{k+L}\sum_{m=j-L}^{j+L}
\E[|x_j-\mu_j|]\E[|x_k-\mu_k|]\E[|x_l-\mu_l|]\E[|x_m-\mu_m|]
\\ \leq & C'p(L+1)^3 \max_{1\leq j\leq p}(E[|x_j-\mu_j|])^4.
%\\ \leq & C\bigg\{pL^3 \max_{j\neq k,k\neq l,l\neq m,m\neq
%j}\E[U_j(x_j,x'_j)U_k(x'_k,x''_k)U_l(x_l'',x'''_l)U_m(x'''_m,x_m)]
%\\&+pL^2 \max_{j\neq l,l\neq m}\E[U_j(x_j,x'_j)U_j(x'_j,x''_j)U_l(x_l'',x'''_l)U_m(x'''_m,x_m)]
%\\&+pL \max_{j\neq m}\E[U_j(x_j,x'_j)U_j(x'_j,x''_j)U_j(x_j'',x'''_j)U_m(x'''_m,x_m)]
%+p \max_{1\leq j\leq p}\E[U_j(x_j,x_j')^4]\bigg\}
%\\ \leq & C'p(L+1)^3\max_{1\leq j\leq p}(\E[x_j^2])^2,
\end{align*}
By the H\"{o}lder's inequality and Lemma \ref{lemma1}, we obtain
\begin{align*}
\E[\widetilde{U}(X,X')^4]\leq
&C\sum_{j=1}^{p}\sum_{k,l,m=j-3L}^{j+3L}|\E[U_j(x_j,x_j')U_k(x_k,x_k')U_l(x_l,x_l')U_m(x_m,x_m')]|
\\&+C'\left\{\sum_{|j-k|\leq
L}dcov(x_j,x_k)^2\right\}^2
\\ \leq
&Cp(6L+1)^3\max_{1\leq j\leq p}\E[U_j(x_j,x_j')^4]
\\&+C'\left\{\sum_{|j-k|\leq L}dcov(x_j,x_k)^2\right\}^2
\\ \leq & C''\left[p(L+1)^3\max_{1\leq j\leq p}\text{var}(x_j)^2+\left\{\sum_{|j-k|\leq
L}dcov(x_j,x_k)^2\right\}^2\right].
\end{align*}
%Also note that
%$$(\sum_{j,j'=1}^pdcov(x_j,x_{j'})^2)^2\geq
%C\left\{p^2L^2\min_{0<|i-j|\leq L}dcov(x_i,x_j)^4+p^2\min_{1\leq
%j\leq p}dcov(x_j)^4\right\}.$$
Therefore, (\ref{eq:con1-2}) and (\ref{eq:con2-2}) are implied by
%\begin{align*}
%&\frac{pL^3 \max_{j\neq k,k\neq l,l\neq m,m\neq
%j}\E[\widetilde{U}_j(x_j,x'_j)\widetilde{U}_k(x'_k,x''_k)\widetilde{U}_l(x_l'',x'''_l)\widetilde{U}_m(x'''_m,x_m)]}{(\sum_{|j-k|\leq
%L}dcov(x_j,x_{k})^2)^2}\rightarrow 0,
%\\&\frac{pL^2 \max_{j\neq l,l\neq m}\E[\widetilde{U}_j(x_j,x'_j)\widetilde{U}_j(x'_j,x''_j)\widetilde{U}_l(x_l'',x'''_l)\widetilde{U}_m(x'''_m,x_m)]}{(\sum_{|j-k|\leq
%L}dcov(x_j,x_{k})^2)^2}\rightarrow 0,
%\\ &\frac{pL \max_{j\neq m}\E[\widetilde{U}_j(x_j,x'_j)\widetilde{U}_j(x'_j,x''_j)\widetilde{U}_j(x_j'',x'''_j)\widetilde{U}_m(x'''_m,x_m)]}{(\sum_{|j-k|\leq
%L}dcov(x_j,x_{k})^2)^2}\rightarrow 0,
%\\&\frac{p\max_{1\leq j\leq p}\E[x_j^4]}{(\sum_{|j-k|\leq
%L}dcov(x_j,x_{k})^2)^2}\rightarrow 0,
%\\ &\frac{pL^3\max_{1\leq j\leq
%p}\E[x_j^4]}{n(\sum_{|j-k|\leq L}dcov(x_j,x_{k})^2)^2}\rightarrow 0.
%\end{align*}
%These conditions are of course implied by the stronger assumption
%that
\begin{align}\label{eq-band}
\frac{p(L+1)^3\max\left\{n^{-1}\max_{1\leq j\leq
p}\text{var}(x_j)^2,\max_{1\leq j\leq p}(\E[|x_j-\mu_j|])^4 \right\}
}{(\sum_{|j-k|\leq L}dcov(x_j,x_{k})^2)^2}\rightarrow 0.
\end{align}
In particular, if $L=o(p^{1/3})$ and $\max_{1\leq j\leq
p}\text{var}(x_j)/\min_{1\leq j\leq p}dcov(x_j,x_j)^2$ is bounded from above,
(\ref{eq-band}) is satisfied.

\begin{lemma}\label{lemma1}
We have
\begin{align*}
&\E[U_j(x_j,x'_j)^4]\leq C\text{var}(x_j)^2,
\\&
|\E[U_j(x_j,x'_j)U_k(x'_k,x''_k)U_l(x_l'',x'''_l)U_m(x'''_m,x_m)]|\leq
C'\E[|x_j-\mu_j|]\E[|x_k-\mu_k|]\E[|x_l-\mu_l|]\E[|x_m-\mu_m|],
\end{align*}
for some constant $C,C'>0$.
\end{lemma}
\begin{proof}
By the triangle inequality, we have
\begin{align*}
\left|\E[K(x_j,x')]-K(x,x')\right|\leq \E[K(x,x_j')]
\end{align*}
for $x,x'\in\mathbb{R}$. Thus we have
$$|U_j(x,x')|\leq \max\{|\E[K(x_j,x_j')]-2\E[K(x,x_j')]|,\E[K(x_j,x_j')]\}:=a_j(x),$$
which implies that
$$\E[U_j(x_j,x'_j)^4]\leq \E[a_j(x_j)^2]\E[a_j(x_j')^2]$$
as $a_j(x_j)$ and $a_j(x_j')$ are independent. Some simple algebra
shows that
\begin{align*}
\E[a_j(x_j)^2]\leq C\text{var}(x_j)^2.
\end{align*}
Therefore, we have $\E[U_j(x_j,x'_j)^4]\leq C^2\text{var}(x_j)^2$.
Similarly, we get
\begin{align*}
&|\E[U_j(x_j,x'_j)U_k(x'_k,x''_k)U_l(x_l'',x'''_l)U_m(x'''_m,x_m)]|
\\ \leq&
\E[a_j(x_j)]\E[a_k(x_k')]\E[a_l(x_l'')]\E[a_m(x_m''')]
\\ \leq & C\E[|x_j-\mu_j|]\E[|x_k-\mu_k|]\E[|x_l-\mu_l|]\E[|x_m-\mu_m|].
\end{align*}

\end{proof}

%When $X\sim N(0,\Sigma)$ with $\Sigma=(\sigma_{jk})_{j,k=1}^p$ and
%$\sigma_{ii}=1$, the analysis in the next section shows that
%\begin{align*}
%\E[\widetilde{U}_j(x_j,x'_j)\widetilde{U}_k(x'_k,x''_k)\widetilde{U}_l(x_l'',x'''_l)\widetilde{U}_m(x'''_m,x_m)]\leq
% C|\sigma_{jk}\sigma_{kl}\sigma_{lm}\sigma_{mj}|
%\end{align*}
%and $dcov(x_j,x_k)^2\geq \pi^{-1}\sigma_{jk}^2.$ Therefore the
%conditions become
%\begin{align}
%\frac{\text{Tr}(\breve{\Sigma}^4)}{\text{Tr}^2(\Sigma^2)}\rightarrow
%0,\quad \frac{pL^3}{n\text{Tr}^2(\Sigma^2)}\rightarrow 0,
%\end{align}
%with $\breve{\Sigma}=(|\sigma_{ij}|)_{i,j=1}^{p}$.

%When $x_1,\dots,x_p$ are mutually independent,
%(\ref{eq:con1-2})-(\ref{eq:con2-2}) become
%\begin{align}
%&\frac{\sum_{j=1}^{p}\E[\widetilde{U}_j(x_j,x'_j)\widetilde{U}(x'_j,x''_j)\widetilde{U}(x''_j,x'''_j)\widetilde{U}(x'''_j,x_j)]}{(\sum_{j=1}^pdcov(x_j)^2)^2}\rightarrow 0, \label{eq:m1}\\
%&\frac{\sum_{j=1}^{p}\E[\widetilde{U}_j(x_j,x'_j)^4]}{n(\sum_{j=1}^pdcov(x_j)^2)^2}\rightarrow 0.\label{eq:m2}
%\end{align}
%Because $\E[\widetilde{U}_j(x_j,x'_j)^4]\leq C_1\E[x_j^4]$, (\ref{eq:m1})-(\ref{eq:m2}) are implied by
%\begin{align*}
%&\frac{\sum_{j=1}^{p}\E[x_j^4]}{(\sum_{j=1}^pdcov(x_j)^2)^2}\rightarrow 0.
%\end{align*}

\textbf{Conditions under Gaussianity:} Under the Gaussian assumption, we verify (\ref{eq:con1-2}) and (\ref{eq:con2-2}) in the following two steps.
\\
\textbf{Step 1:} We have
\begin{align*}
&V(\sigma_{jl},\sigma_{lm},\sigma_{mk},\sigma_{kj})
\\:=&\E[U_j(x_j,x''_j)U_k(x'_k,x''_k)U_l(x_l,x'''_l)U_m(x'_m,x'''_m)]
\\=&\int\{f_{jl}(t_1,t_2)-f_j(t_1)f_l(t_2)\}\{f_{lm}(-t_2,-t_4)-f_l(-t_2)f_m(-t_4)\}
\\&\{f_{mk}(t_4,t_3)-f_m(t_4)f_k(t_3)\}\{f_{kj}(-t_3,-t_1)-f_k(-t_3)f_j(-t_1)\}
\\&w(t_1)w(t_2)w(t_3)w(t_4)dt_1dt_2dt_3dt_4
\\=&\int e^{-(t_1^2+t_2^2+t_3^2+t_4^2)}(1-e^{-\sigma_{jl}t_1t_2})(1-e^{-\sigma_{lm}t_2t_4})(1-e^{-\sigma_{mk}t_3t_4})(1-e^{-\sigma_{kj}t_1t_3})
\\&w(t_1)w(t_2)w(t_3)w(t_4)dt_1dt_2dt_3dt_4.
\end{align*}
Using the power series expansion, we have
\begin{align*}
&(1-e^{-\sigma_{jl}t_1t_2})(1-e^{-\sigma_{lm}t_2t_4})(1-e^{-\sigma_{mk}t_3t_4})(1-e^{-\sigma_{kj}t_1t_3})
\\=&\sum_{i_1,i_2,i_3,i_4=1}^{\infty}\frac{(-1)^{i_1+i_2+i_3+i_4}}{i_1!i_2!i_3!i_4!}\sigma_{jl}^{i_1}\sigma_{lm}^{i_2}\sigma_{mk}^{i_3}\sigma_{kj}^{i_4}(t_1t_2)^{i_1}(t_2t_4)^{i_2}(t_3t_4)^{i_3}(t_1t_3)^{i_4},
\end{align*}
which implies that
\begin{align*}
&V(\sigma_{jl},\sigma_{lm},\sigma_{mk},\sigma_{kj})
\\=&C_1\int\sum_{i_1,i_2,i_3,i_4=1}^{\infty}\frac{(-1)^{i_1+i_2+i_3+i_4}}{i_1!i_2!i_3!i_4!}\sigma_{jl}^{i_1}\sigma_{lm}^{i_2}\sigma_{mk}^{i_3}\sigma_{kj}^{i_4}
e^{-(t_1^2+t_2^2+t_3^2+t_4^2)}
\\&(t_1t_2)^{i_1-1}(t_2t_4)^{i_2-1}(t_3t_4)^{i_3-1}(t_1t_3)^{i_4-1}dt_1dt_2dt_3dt_4
\\=&C_1\sigma_{jl}\sigma_{lm}\sigma_{mk}\sigma_{kj}\sum_{i_1,i_2,i_3,i_4=0}^{\infty}\frac{(-1)^{i_1+i_2+i_3+i_4}}{(i_1+1)!(i_2+1)!(i_3+1)!(i_4+1)!}\sigma_{jl}^{i_1}\sigma_{lm}^{i_2}\sigma_{mk}^{i_3}\sigma_{kj}^{i_4}
\int e^{-(t_1^2+t_2^2+t_3^2+t_4^2)}
\\&(t_1t_2)^{i_1}(t_2t_4)^{i_2}(t_3t_4)^{i_3}(t_1t_3)^{i_4}dt_1dt_2dt_3dt_4.
\end{align*}
Note that $\int
e^{-(t_1^2+t_2^2+t_3^2+t_4^2)}(t_1t_2)^{i_1}(t_2t_4)^{i_2}(t_3t_4)^{i_3}(t_1t_3)^{i_4}dt_1dt_2dt_3dt_4$
is nonzero if and only if $i_1,i_2,i_3,i_4$ are all even or odd
numbers simultaneously. Thus using the fact that $|\sigma_{kl}|\leq
1$, we have
\begin{align*}
&V(\sigma_{jl},\sigma_{lm},\sigma_{mk},\sigma_{kj})
\\ = &C_1\sigma_{jl}\sigma_{lm}\sigma_{mk}\sigma_{kj}\sum_{i_1,i_2,i_3,i_4 \text{ are all even or odd}}\frac{1}{(i_1+1)!(i_2+1)!(i_3+1)!(i_4+1)!}
\sigma_{jl}^{i_1}\sigma_{lm}^{i_2}\sigma_{mk}^{i_3}\sigma_{kj}^{i_4}
\\&\int
e^{-(t_1^2+t_2^2+t_3^2+t_4^2)}(t_1t_2)^{i_1}(t_2t_4)^{i_2}(t_3t_4)^{i_3}(t_1t_3)^{i_4}dt_1dt_2dt_3dt_4
\\ \leq & C_1|\sigma_{jl}\sigma_{lm}\sigma_{mk}\sigma_{kj}|\sum_{i_1,i_2,i_3,i_4 \text{ are all even or odd}}\frac{1}{(i_1+1)!(i_2+1)!(i_3+1)!(i_4+1)!}\int e^{-(t_1^2+t_2^2+t_3^2+t_4^2)}
\\&(t_1t_2)^{i_1}(t_2t_4)^{i_2}(t_3t_4)^{i_3}(t_1t_3)^{i_4}dt_1dt_2dt_3dt_4
\\ = & |\sigma_{jl}\sigma_{lm}\sigma_{mk}\sigma_{kj}|V(1,1,1,1),
\end{align*}
where
$V(1,1,1,1)=\E[U_j(x_j,x''_j)U_k(x'_j,x''_j)U_l(x_j,x'''_j)U_m(x'_j,x'''_j)]\leq
\E[U_j(x_j,x_j')^4]<\infty$ with
$x_j,x_j',x_j'',x_j'''\sim^{i.i.d} N(0,1)$. On the other hand, in
view of the proof of Theorem 7 of Sz\'{e}kely et al. (2007), we have
$$dcov(x_j,x_{j'})^2=\frac{\sigma_{j,j'}^2}{\pi^2}\sum_{k=1}^{+\infty}\frac{2^{2k}-2}{(2k)!}\sigma_{j,j'}^{2(k-1)}\int_{\mathbb{R}^2}e^{-t^2-s^2}(ts)^{2(k-1)}dtds\geq \frac{\sigma_{j,j'}^2}{\pi},$$
which implies that $\sum_{j,j'=1}^pdcov(x_j,x_{j'})^2\geq
\pi^{-1}\sum_{j,j'=1}^{p}\sigma_{j,j'}^2=\pi^{-1}\text{Tr}(\Sigma^2).$
Thus Condition (\ref{eq:con1-2}) is implied by the assumption $\frac{\sum_{j,k,l,m=1}^{p}|\sigma_{jk}\sigma_{kl}\sigma_{lm}\sigma_{mj}|}{\text{Tr}^2(\Sigma^2)}\rightarrow
0$.

%$\text{Tr}(\breve{\Sigma}^4)/\text{Tr}^2(\Sigma^2)\rightarrow 0.$

%\begin{remark}
%{\rm Consider the case of equal correlation, i.e., $\sigma_{ij}=\rho$ for
%$i\neq j$ and $\sigma_{ii}=1$. We have
%\begin{align*}
%U_{jklm}(t_1,t_2,t_3,t_4)=&\sum^{+\infty}_{s=2}\frac{(-1)^s\rho^s}{s!}\sum_{*}\frac{s!}{s_1!s_2!s_3!s_4!s_5!s_6!}t_1^{s_1+s_2+s_3}t_2^{s_1+s_4+s_5}t_3^{s_2+s_4+s_6}t_4^{s_3+s_5+s_6}.
%\end{align*} }
%\end{remark}

\textbf{Step 2:} %Assume that $\max_{i\neq j}|\sigma_{ij}|\leq c<1$. %Let $a_1=t_jt_k\sigma_{jk}$, $a_2=t_jt_l\sigma_{jl}$, $a_3=t_jt_m\sigma_{jm}$, $a_4=t_kt_m\sigma_{km}$, $a_5=$
Using (\ref{eq:ch}) and the power series expansion, we deduce after laborious calculations that
\begin{align*}
%&\E [\widetilde{U}(X,X')^4]
%\\ \leq & C^4\sum_{j,k,l,m=1}^p\int
%\bigg|-3f_j(t_1)f_k(t_2)f_l(t_3)f_m(t_4)+f_{jl}(t_1,t_3)f_k(t_2)f_m(t_4)+f_{jk}(t_1,t_2)f_l(t_3)f_m(t_4)
%\\&+f_{jm}(t_1,t_4)f_l(t_2)f_k(t_3)+f_{km}(t_2,t_4)f_j(t_1)f_l(t_3)+f_{kl}(t_2,t_3)f_j(t_1)f_m(t_4)
%\\&+f_{lm}(t_3,t_4)f_j(t_1)f_k(t_2)-f_{jlm}(t_1,t_3,t_4)f_k(t_2)-f_{klm}(t_2,t_3,t_4)f_j(t_1)-f_{jkl}(t_1,t_2,t_3)f_m(t_4)
%\\&-f_{jkm}(t_1,t_2,t_4)f_l(t_3)+f_{jklm}(t_1,t_2,t_3,t_4)\bigg|^2
%\\&w(t_1)w(t_2)w(t_3)w(t_4)dt_1dt_2dt_3dt_4
& \E[U_j(x_j,x_j')U_k(x_k,x_k')U_l(x_l,x_l')U_m(x_m,x_m')]
\\ =& \int e^{-(t_1^2+t_2^2+t_3^2+t_4^2)}
U_{jklm}(t_1,t_2,t_3,t_4)^2w(t_1)w(t_2)w(t_3)w(t_4)dt_1dt_2dt_3dt_4,
\end{align*}
where
\begin{align*}
U_{jklm}(t_1,t_2,t_3,t_4)=&\sum^{+\infty}_{u=2}\sum_{\mathbf{s}}\frac{(-1)^u}{\mathbf{s}!}t_1^{s_1+s_2+s_3}t_2^{s_1+s_4+s_5}t_3^{s_2+s_4+s_6}t_4^{s_3+s_5+s_6}\sigma_{jk}^{s_1}\sigma_{jl}^{s_2}\sigma_{jm}^{s_3}\sigma_{kl}^{s_4}\sigma_{km}^{s_5}\sigma_{lm}^{s_6}
\end{align*}
with $\sum_{\mathbf{s}}$ denoting the summation over all
$\mathbf{s}=(s_1,s_2,s_3,s_4,s_5,s_6)$ such that
$s_1+s_2+s_3+s_4+s_5+s_6=u,s_1+s_2+s_3\geq 1,s_1+s_4+s_5\geq
1,s_2+s_4+s_6\geq 1$ and $s_3+s_5+s_6\geq 1,$ and $\mathbf{s}!=\prod^{6}_{i=1} s_i!$. Therefore
\begin{align*}
 &\E[U_j(x_j,x_j')U_k(x_k,x_k')U_l(x_l,x_l')U_m(x_m,x_m')]
 \\ = & \sum^{+\infty}_{u,v=2}\sum_{\mathbf{s}}\sum_{\mathbf{r}}\frac{(-1)^{u+v}}{\mathbf{s}!\mathbf{r}!}\sigma_{jk}^{s_1+r_1}\sigma_{jl}^{s_2+r_2}\sigma_{jm}^{s_3+r_3}\sigma_{kl}^{s_4+r_4}\sigma_{km}^{s_5+r_5}\sigma_{lm}^{s_6+r_6}\int
e^{-(t_1^2+t_2^2+t_3^2+t_4^2)}
\\&t_1^{s_1+s_2+s_3+r_1+r_2+r_3}t_2^{s_1+s_4+s_5+r_1+r_4+r_5}t_3^{s_2+s_4+s_6+r_2+r_4+r_6}t_4^{s_3+s_5+s_6+r_3+r_5+r_6}
\\&w(t_1)w(t_2)w(t_3)w(t_4)dt_1dt_2dt_3dt_4
\\=& \sum_{h=4}^{+\infty}(-1)^h\sum_{u+v=h;u,v\geq 2}\sum_{\mathbf{s}}\sum_{\mathbf{v}}\frac{1}{\mathbf{s}!\mathbf{r}!}\sigma_{jk}^{s_1+r_1}\sigma_{jl}^{s_2+r_2}\sigma_{jm}^{s_3+r_3}\sigma_{kl}^{s_4+r_4}\sigma_{km}^{s_5+r_5}\sigma_{lm}^{s_6+r_6}a_{\mathbf{s},\mathbf{r}},
\end{align*}
where
\begin{align*}
a_{\mathbf{s},\mathbf{r}}=&\int e^{-(t_1^2+t_2^2+t_3^2+t_4^2)}t_1^{s_1+s_2+s_3+r_1+r_2+r_3}t_2^{s_1+s_4+s_5+r_1+r_4+r_5}t_3^{s_2+s_4+s_6+r_2+r_4+r_6}t_4^{s_3+s_5+s_6+r_3+r_5+r_6}
\\&w(t_1)w(t_2)w(t_3)w(t_4)dt_1dt_2dt_3dt_4,
\end{align*}
and $s_1+s_2+s_3+r_1+r_2+r_3$, $s_1+s_4+s_5+r_1+r_4+r_5$,
$s_2+s_4+s_6+r_2+r_4+r_6$, $s_3+s_5+s_6+r_3+r_5+r_6$ are all even
numbers. Under the above constraint, the integration inside the curly bracket is always nonnegative.
For $h\geq 4$ and $\mathbf{s},\mathbf{r}$ satisfying the above constraint, the term $|\sigma_{jk}^{s_1+r_1}\sigma_{jl}^{s_2+r_2}\sigma_{jm}^{s_3+r_3}\sigma_{kl}^{s_4+r_4}\sigma_{km}^{s_5+r_5}|$
is bounded by one of the following terms
\begin{align*}
&\sigma_{jk}^2\sigma_{lm}^2,\sigma_{jl}^2\sigma_{km}^2,\sigma_{jm}^2\sigma_{kl}^2,\sigma_{jk}^2\sigma_{jl}^2\sigma_{jm}^2,\sigma_{kj}^2\sigma_{kl}^2\sigma_{km}^2,\sigma_{lj}^2\sigma_{lk}^2\sigma_{lm}^2,\sigma_{mj}^2\sigma_{mk}^2\sigma_{ml}^2,
\\&|\sigma_{jk}\sigma_{kl}\sigma_{lm}\sigma_{mj}|,|\sigma_{jk}\sigma_{km}\sigma_{ml}\sigma_{lj}|,|\sigma_{jl}\sigma_{lk}\sigma_{km}\sigma_{mj}|.
\end{align*}
%This provides a way to classify the terms in the summation.
On the other hand, as $\E[U_j(x_j,x_j')^4]=\sum_{h=4}^{+\infty}(-1)^h\sum_{u+v=h;u,v\geq 2}\sum_{\mathbf{s}}\sum_{\mathbf{v}}a_{\mathbf{s},\mathbf{r}}/(\mathbf{s}!\mathbf{r}!)<\infty$, we must have
$\sum_{u+v=h;u,v\geq 2}\sum_{\mathbf{s}}\sum_{\mathbf{v}}a_{\mathbf{s},\mathbf{r}}/(\mathbf{s}!\mathbf{r}!)\rightarrow 0$
as $k\rightarrow +\infty.$
%\textbf{To get this result, we have to change the order of summation and integration. To be rigorous, we should argue why it is valid to do so.}
Thus we get
\begin{align*}
&\E[U_j(x_j,x_j')U_k(x_k,x_k')U_l(x_l,x_l')U_m(x_m,x_m')]
\\ \leq & (\sigma_{jk}^2\sigma_{lm}^2+\sigma_{jl}^2\sigma_{km}^2+\sigma_{jm}^2\sigma_{kl}^2+|\sigma_{jk}\sigma_{kl}\sigma_{lm}\sigma_{mj}|+|\sigma_{jk}\sigma_{km}\sigma_{ml}\sigma_{lj}|+|\sigma_{jl}\sigma_{lk}\sigma_{km}\sigma_{mj}|)
\\&\times\left(\sum_{h=4}^{+\infty}c^{h-4}\sum_{u+v=h;u,v\geq 2}\sum_{\mathbf{s}}\sum_{\mathbf{v}}\frac{1}{\mathbf{s}!\mathbf{r}!}a_{\mathbf{s},\mathbf{r}}\right)
\\&+(\sigma_{jk}^2\sigma_{jl}^2\sigma_{jm}^2+\sigma_{kj}^2\sigma_{kl}^2\sigma_{km}^2+\sigma_{lj}^2\sigma_{lk}^2\sigma_{lm}^2+\sigma_{mj}^2\sigma_{mk}^2\sigma_{ml}^2)
\\&\times \left(\sum_{h=6}^{+\infty}c^{h-6}\sum_{u+v=h;u,v\geq 2}\sum_{\mathbf{s}}\sum_{\mathbf{v}}\frac{1}{\mathbf{s}!\mathbf{r}!}a_{\mathbf{s},\mathbf{r}}\right)
\\ \leq & C_1(\sigma_{jk}^2\sigma_{lm}^2+\sigma_{jl}^2\sigma_{km}^2+\sigma_{jm}^2\sigma_{kl}^2+|\sigma_{jk}\sigma_{kl}\sigma_{lm}\sigma_{mj}|+|\sigma_{jk}\sigma_{km}\sigma_{ml}\sigma_{lj}|+|\sigma_{jl}\sigma_{lk}\sigma_{km}\sigma_{mj}|
\\&+\sigma_{jk}^2\sigma_{jl}^2\sigma_{jm}^2+\sigma_{kj}^2\sigma_{kl}^2\sigma_{km}^2+\sigma_{lj}^2\sigma_{lk}^2\sigma_{lm}^2+\sigma_{mj}^2\sigma_{mk}^2\sigma_{ml}^2).
\end{align*}
It follows that
\begin{align*}
&\E[\widetilde{U}(X,X')^4]
\\ \leq& C_1\sum_{j,k,l,m=1}^{p}(\sigma_{jk}^2\sigma_{lm}^2+\sigma_{jl}^2\sigma_{km}^2+\sigma_{jm}^2\sigma_{kl}^2+|\sigma_{jk}\sigma_{kl}\sigma_{lm}\sigma_{mj}|+|\sigma_{jk}\sigma_{km}\sigma_{ml}\sigma_{lj}|+|\sigma_{jl}\sigma_{lk}\sigma_{km}\sigma_{mj}|
\\&+\sigma_{jk}^2\sigma_{jl}^2\sigma_{jm}^2+\sigma_{kj}^2\sigma_{kl}^2\sigma_{km}^2+\sigma_{lj}^2\sigma_{lk}^2\sigma_{lm}^2+\sigma_{mj}^2\sigma_{mk}^2\sigma_{ml}^2)
\\ \leq& C_2\left(\text{Tr}^2(\Sigma^2)+\text{Tr}(\breve{\Sigma}^4)+\sum_{j=1}^{p}(\sum_{k=1}^{p}\sigma_{jk}^2)^3\right),
\end{align*}
where $\breve{\Sigma}=(|\sigma_{jk}|)$. Therefore, (\ref{eq:con2-2}) is satisfied under (15) in the paper.
%\textbf{The coefficients of
%$\sigma_{jk}^{s_1+r_1}\sigma_{jl}^{s_2+r_2}\sigma_{jm}^{s_3+r_3}\sigma_{kl}^{s_4+r_4}\sigma_{km}^{s_5+r_5}$
%can be positive or negative. This makes it hard to bound
%$$\E[\widetilde{U}_j(x_j,x_j')\widetilde{U}_k(x_k,x_k')\widetilde{U}_l(x_l,x_l')\widetilde{U}_m(x_m,x_m')].$$ }

\subsection{Consistency of wild Bootstrap and variance estimator}
In this subsection, we prove Theorem 2.2 (variance ratio consistency) and
Theorem 2.3 (the bootstrap consistency) under Conditions (8), (10) and (11) in the main paper.
In particular, the variance ratio consistency is an intermediate step in the proof of bootstrap consistency, see the details in Step 6 below.

Define
$d_{kl}=\sum^{p}_{j=1}\widetilde{A}_{kl}(j)\widetilde{B}_{kl}$
and note that $MDD^{*}_n(Y|x_j)^2=\frac{2}{n(n-1)}\sum_{k<l}d_{kl}
e_{k}e_{l}.$ Denote by $\E^*$ and $\text{cov}^*/\text{var}^*$
the expectation and covariance/variance in the bootstrap world. Then we have $\E^*[MDD^{*}_n(Y|x_j)^2]=0$ for any $1\leq j\leq p,$ and
$$\text{var}^*\left(\sqrt{\binom{n}{2}}\sum_{j=1}^pMDD^{*}_n(Y|x_j)^2\right)=\frac{2}{n(n-1)}\sum_{k<l}d_{kl}^2=\hat{\mathcal{S}}^2.$$
Define the event
\begin{equation}\label{eq1}
\mathcal{A}=\left\{\frac{\sum_{j\geq j'}\left(\sum^{j'-1}_{i=1}d_{ij}d_{ij'}\right)^2}{(\sum_{k<l}d_{kl}^2)^2}\rightarrow 0\right\}.
\end{equation}
Our proof involves the following four steps:
\begin{enumerate}
\item Conditional on $\mathcal{A}$, we show
$$\sqrt{\binom{n}{2}}\frac{\sum^{p}_{j=1}MDD^{*}_n(Y|x_j)^2}{\hat{\mathcal{S}}}\rightarrow^{d^*} N(0,1),$$
where $\rightarrow^{d^*}$ denotes convergence in distribution with
respect to $\{e_i\}$.
\item Conditional on $\mathcal{A}$, we show
$$\frac{\hat S^{*}}{\hat{\mathcal{S}}}\rightarrow^{p^*} 1,$$
where $\rightarrow^{p^*}$ denotes convergence in probability with
respect to $\{e_i\}$.
\item Under (\ref{eq-1})-(\ref{eq:con2-2}), we prove that
$$\frac{\sum_{j\geq j'}\left(\sum^{j'-1}_{i=1}d_{ij}d_{ij'}\right)^2}{(\sum_{k<l}d_{kl}^2)^2}\rightarrow^p 0.$$
\item Finally, combining the above results, we obtain the convergence in distribution in probability, i.e.,
$$\sqrt{\binom{n}{2}}\frac{\sum^{p}_{j=1}MDD^{*}_n(Y|x_j)^2}{\hat{S}^*}\rightarrow^{\mathcal{D}^*}
N(0,1).$$
\end{enumerate}

For clarity, we present the proofs in the following 6 steps.

\textbf{Step 1:} The basic
idea is to apply the martingale CLT to the quadratic form
$\sum^{p}_{j=1}MDD^{*}_n(Y|x_j)^2.$ The arguments are similar to those in Section
\ref{sec:normal}. Define
$S_r^*=\sum^{r}_{j=2}\sum_{i=1}^{j-1}d_{ij}e_ie_j.$ It is
straightforward to see that $S_r^*$ is a mean-zero martingale with
respect to the filtration
$\mathcal{F}_r^*=\sigma\{e_1,e_2,\dots,e_r\}$. We establish the
asymptotic normality by Corollary 3.1 of Hall and Heyde (1980).
Define $W_j^*=\sum_{i=1}^{j-1}d_{ij}e_ie_j$ and note that
$\E^*[W_j^{*2}|\mathcal{F}_{j-1}^*]=\sum_{i,k=1}^{j-1}d_{ij}d_{kj}e_ie_k,$
and $2\sum_{j=1}^{n}\E^*[W_j^{*2}]/\{n(n-1)\}=\hat{\mathcal{S}}^2$.
Direct calculation yields that for $j\geq j'$,
\begin{align*}
&\text{cov}^*(\E[W_j^{*2}|\mathcal{F}_{j-1}^*],\E[W_{j'}^{*2}|\mathcal{F}_{j'-1}^*])
\\=&\sum_{i,k=1}^{j-1}\sum_{i',k'=1}^{j'-1}d_{ij}d_{kj}d_{i'j'}d_{k'j'}\text{cov}(e_ie_k,e_{i'}e_{k'})
\\=&\sum_{i,k,i',k'=1}^{j'-1}d_{ij}d_{kj}d_{i'j'}d_{k'j'}\text{cov}(e_ie_k,e_{i'}e_{k'})
\\=&2\sum_{i=1}^{j'-1}d_{ij}^2d_{ij'}^2+2\sum_{1\leq i\neq k\leq j'-1}d_{ij}d_{ij'}d_{kj}d_{kj'}=2\left(\sum^{j'-1}_{i=1}d_{ij}d_{ij'}\right)^2.
\end{align*}
It implies that on $\mathcal{A}$
\begin{align*}
&\frac{4}{n^2(n-1)^2}\sum_{j,j'=2}^{n}\text{cov}^*(\E[W_j^{*2}|\mathcal{F}_{j-1}^*],\E[W_{j'}^{*2}|\mathcal{F}_{j'-1}^*])
\\=&\frac{8}{n^2(n-1)^2}\sum_{j,j'=2}^{n}\left(\sum^{\min\{j,j'\}-1}_{i=1}d_{ij}d_{ij'}\right)^2=o(\hat{\mathcal{S}}^4).
\end{align*}
We also note that on $\mathcal{A}$
\begin{align*}
\sum_{j=2}^{n}\E^*[|W_j^*|^4]=&3\sum_{j=2}^{n}\sum_{i_1,i_2,i_3,i_4=1}^{j-1}d_{i_1j}d_{i_2j}d_{i_3j}d_{i_4j}\E[e_{i_1}e_{i_2}e_{i_3}e_{i_4}]
\\=&9\sum_{j=2}^{n}\sum_{i=1}^{j-1}d_{ij}^4+9\sum_{j=2}^{n}\sum_{1\leq i_1\neq i_2\leq j-1}d_{i_1j}^2d_{i_2j}^2=9\sum_{j=2}^{n}\left(\sum_{i=1}^{j-1}d_{ij}^2\right)^2
\\=&o\left(\left(\sum_{k<l}d_{kl}^2\right)^2\right).
\end{align*}
Therefore by the martingale CLT (see Hall and Heyde (1980)), we have on $\mathcal{A}$,
\begin{align*}
\sqrt{\binom{n}{2}}\frac{\sum^{p}_{j=1}MDD^{*}_n(Y|x_j)^2}{\hat{\mathcal{S}}}\rightarrow^{d^*} N(0,1).
\end{align*}

\textbf{Step 2:} To show the ratio consistency for the bootstrap
variance estimator, we note that conditional on the sample, $\hat
S^{*2}$ is an unbiased estimator for $\hat{\mathcal{S}}^2$.
Therefore, it suffices to show that conditional on $\mathcal{A}$,
$\frac{\text{var}^*(\hat S^{*2})}{\hat{\mathcal{S}}^4}\rightarrow
0.$ Simple algebra shows that
\begin{align*}
\frac{\text{var}^*(\hat{S}^{*2})}{\hat{\mathcal{S}}^4}=&\frac{1}{\hat{\mathcal{S}}^4}\text{var}^*\left(\frac{1}{\binom{n}{2}}
\sum_{1\leq k<l\leq n}d_{kl}^2 e_{k} ^2
e_{l}^2\right)
\\=&\frac{4}{n^2(n-1)^2\hat{\mathcal{S}}^4}\bigg\{\sum_{k<l}d_{kl}^4\text{cov}(e_k^2e_l^2,e_{k}^2e_{l}^2)+\sum_{k<l\neq l'}d_{kl}^2d_{kl'}^2\text{cov}(e_k^2e_l^2,e_{k}^2e_{l'}^2)
\\&+\sum_{k\neq k'<l}d_{kl}^2d_{k'l}^2\text{cov}(e_k^2e_l^2,e_{k'}^2e_{l}^2)+2\sum_{k<l<l'}d_{kl}^2d_{ll'}^2\text{cov}(e_k^2e_l^2,e_{l}^2e_{l'}^2)\bigg\}
\\ \leq & \frac{C}{(\sum_{k<l}d_{kl}^2)^2}\bigg\{\sum_{k<l}d_{kl}^4+\sum_{k\neq l\neq l'}d_{kl}^2d_{ll'}^2\bigg\},
\end{align*}
for some constant $C>0$. On $\mathcal{A}$, we have
\begin{align*}
\frac{\sum_{k<l}d_{kl}^4}{(\sum_{k<l}d_{kl}^2)^2}\rightarrow
0,\quad \frac{\sum_{k\neq l\neq
l'}d_{kl}^2d_{kl'}^2}{(\sum_{k<l}d_{kl}^2)^2}\rightarrow 0,
\end{align*}
which indicates that
$$\frac{\hat S^{*}}{\hat{\mathcal{S}}}\rightarrow^{p^*} 1.$$

\textbf{Step 3:} To deal with the $\mathcal{U}$-centered version of $A_{ij}$ (here we omit the superscript), we write $\widetilde{A}_{ij}$ in the following form,
\begin{align*}
\widetilde{A}_{ij}=&\frac{n-3}{n-1}A_{ij}-\frac{n-3}{(n-1)(n-2)}\sum_{l\notin \{i,j\}}A_{il}-\frac{n-3}{(n-1)(n-2)}\sum_{k\notin \{i,j\} }A_{kj}
\\&+\frac{1}{(n-1)(n-2)}\sum_{k,l\notin \{i,j\},k\neq l}A_{kl}
\\=&\frac{n-3}{n-1}\bar{A}_{ij}-\frac{n-3}{(n-1)(n-2)}\sum_{l\notin \{i,j\}}\bar{A}_{il}-\frac{n-3}{(n-1)(n-2)}\sum_{k\notin \{i,j\}}\bar{A}_{kj}
\\&+\frac{2}{(n-1)(n-2)}\sum_{k,l\notin \{i,j\},k<l}\bar{A}_{kl}
\\:=&I_{n,1}-I_{n,2}-I_{n,3}+I_{n,4}
\end{align*}
where $\bar{A}_{ij}=A_{ij}-\E[A_{il}|X_i]-\E[A_{kj}|X_j]+\E[A_{kl}]$
is the double centered version of $A_{ij}$. Recall that $\bar{A}_{ij}^{(k)}=-U_k(x_{ik},x_{jk}).$ A useful
property is that the four terms $I_{n,k},1\leq k\leq 4$ are
uncorrelated with each other.
\begin{remark}
{\rm This expression provides an alternative way of justifying the unbiasedness of the $\mathcal{U}$-centered estimator. For example, $\E[\widetilde{A}_{ij}^2]$ is simply equal to the sum of the variances of the above four terms (since they are uncorrelated by construction), which is equal to
\begin{align*}
\left\{\frac{(n-3)^2}{(n-1)^2}+\frac{2(n-3)^2}{(n-1)^2(n-2)}+\frac{2(n-3)}{(n-1)^2(n-2)}\right\}\E[\bar{A}_{ij}^2]=\frac{n-3}{n-1}\E[\bar{A}_{ij}^2].
\end{align*}
Hence $\sum_{i\neq j}\widetilde{A}_{ij}^2/\{n(n-3)\}$ is an unbiased estimator for $\E[\bar{A}_{ij}^2]$.
}
\end{remark}

\textbf{Step 4:} In this step, we prove that
\begin{align}\label{eq:st4}
\frac{\sum_{j\geq j'}\left(\sum^{j'-1}_{i=1}d_{ij}d_{ij'}\right)^2}{(\sum_{k<l}d_{kl}^2)^2}\rightarrow^p 0.
\end{align}
The above condition can be viewed as the sample version of Condition (8) in the paper. In Step 6 below, we show that $\hat{\mathcal{S}}^2$ is ratio-consistent under the null, i.e.,
\begin{align*}
\frac{\hat{\mathcal{S}}^2}{\mathcal{S}^2}\rightarrow^{p} 1.
\end{align*}
To simplify the calculation, we shall assume (\ref{eq-1}).
However, we emphasize that it is not essential and can be relaxed.
Using the expression in Step 3, we can show $\E[\widetilde{B}_{ij}^4|X_i,X_j]\leq C_1^4$ for some constant $C_1>0.$
Using this fact, (\ref{eq-1}) and again the expression in Step 3 (see more details in Step 5 below), we get
\begin{align*}
\E [d_{ij}^4]=&\sum_{j_1,j_2,j_3,j_4=1}^{p}\E[\widetilde{A}_{ij}^{(j_1)}\widetilde{A}_{ij}^{(j_2)}\widetilde{A}_{ij}^{(j_3)}\widetilde{A}_{ij}^{(j_4)}\widetilde{B}_{ij}^4]
\\ \leq & C_1^4\sum_{j_1,j_2,j_3,j_4=1}^{p}\E[\widetilde{A}_{ij}^{(j_1)}\widetilde{A}_{ij}^{(j_2)}\widetilde{A}_{ij}^{(j_3)}\widetilde{A}_{ij}^{(j_4)}]=C_1'\E[\widetilde{U}(X,X')^4].
\end{align*}
Similarly, for $i<j'<j$, we have
\begin{align*}
\E [d_{ij}^2d_{ij'}^2]\leq& C_1^4\sum_{j_1,j_2,j_3,j_4=1}^{p}\E[\widetilde{A}_{ij}^{(j_1)}\widetilde{A}_{ij}^{(j_2)}\widetilde{A}_{ij'}^{(j_3)}\widetilde{A}_{ij'}^{(j_4)}]
\\=& C_2'\E[\widetilde{U}(X,X')^2\widetilde{U}(X,X'')^2]+C_2''n^{-1}\E[\widetilde{U}(X,X')^4],
\end{align*}
and for $i\neq i'<j'<j,$
\begin{align*}
\E [d_{ij}d_{ij'}d_{i'j}d_{i'j'}]\leq & C_1^4 \sum_{j_1,j_2,j_3,j_4=1}^{p}\E[\widetilde{A}_{ij}^{(j_1)}\widetilde{A}_{ij'}^{(j_2)}\widetilde{A}_{i'j}^{(j_3)}\widetilde{A}_{i'j'}^{(j_4)}]
\\=& C_3' \E[\widetilde{U}(X,X')\widetilde{U}(X',X'')\widetilde{U}(X'',X''')\widetilde{U}(X''',X)]+C_3''n^{-1}\E[\widetilde{U}(X,X')^4].
\end{align*}
Therefore, by (\ref{eq:con1-2})-(\ref{eq:con2-2}), we have,
\begin{align*}
&\frac{\E\sum_{j>j'}\left(\sum^{j'-1}_{i=1}d_{ij}d_{ij'}\right)^2}{n^4\mathcal{S}^4}
\\ \leq & C_1\left(\frac{\E[\widetilde{U}(X,X')\widetilde{U}(X',X'')\widetilde{U}(X'',X''')\widetilde{U}(X''',X)]}{\mathcal{S}^4}+\frac{\E[\widetilde{U}(X,X')^4]}{n\mathcal{S}^4}\right)\rightarrow 0,
\end{align*}
which implies (\ref{eq:st4}) by the Markov inequality.

\textbf{Step 5:} We provide some details on the calculation in Step 4 above. The key idea is to use the alternative expression given in Step 3. Let
$$c_n=\frac{(n-3)^4}{(n-1)^4}+\frac{2(n-3)^4}{(n-1)^4(n-2)^3}+\frac{2(n-3)}{(n-1)^4(n-2)^3},$$
such that $c_n\rightarrow 1.$ We begin by considering
\begin{align*}
\E[d_{ij}^4]\leq & C_1^4\sum_{j_1,j_2,j_3,j_4=1}^p\E[\widetilde{A}_{ij}^{(j_1)}\widetilde{A}_{ij}^{(j_2)}\widetilde{A}_{ij}^{(j_3)}\widetilde{A}_{ij}^{(j_4)}]
\\=&c_nC_1^4\sum_{j_1,j_2,j_3,j_4=1}^p\E[\bar{A}_{ij}^{(j_1)}\bar{A}_{ij}^{(j_2)}\bar{A}_{ij}^{(j_3)}\bar{A}_{ij}^{(j_4)}]+\mathcal{R}_n
\\=&c_nC_1^4\E[\widetilde{U}(X,X')^4]+\mathcal{R}_n,
\end{align*}
where $\mathcal{R}_n=O(n^{-1}\E[\widetilde{U}(X,X')^4])$ is the smaller order term. Some typical terms in $\mathcal{R}_n$ include
\begin{align*}
&O(n^{-4})\sum_{j_1,j_2,j_3,j_4=1}^p\E[\bar{A}^{(j_1)}_{ij}\bar{A}^{(j_2)}_{jl}\bar{A}^{(j_3)}_{lm}\bar{A}^{(j_4)}_{mi}]=O\left(n^{-4}\E[\widetilde{U}(X,X')^4]\right),
\\&O(n^{-4})\sum_{j_1,j_2,j_3,j_4=1}^p\E[\bar{A}^{(j_1)}_{ij}\bar{A}^{(j_2)}_{ij}\bar{A}^{(j_3)}_{kl}\bar{A}^{(j_4)}_{kl}]=O\left(n^{-4}\E[\widetilde{U}(X,X')^4]\right),\quad \{i,j\}\cap\{k,l\}=\emptyset,
\\&O(n^{-2})\sum_{j_1,j_2,j_3,j_4=1}^p\E[\bar{A}^{(j_1)}_{ij}\bar{A}^{(j_2)}_{ij}\bar{A}^{(j_3)}_{im}\bar{A}^{(j_4)}_{im}]=O\left(n^{-2}\E[\widetilde{U}(X,X')^4]\right),\quad m\in\{i,j\},
\end{align*}
where we have used the H\"{o}lder's inequality.
The calculation for $\E [d_{ij}^2d_{ij'}^2]$ and $\E [d_{ij}d_{ij'}d_{i'j}d_{i'j'}]$ is similar and hence is skipped.

\textbf{Step 6:} To finish the proof, we show the ratio consistency
$\frac{\hat{\mathcal{S}}^2}{\mathcal{S}^2}\rightarrow^{p} 1.$ By
the Markov inequality, we only need to show that
\begin{align*}
\frac{\text{var}(\sum_{k<l}d_{kl}^2)}{n^4\mathcal{S}^4}+\frac{|\E[d_{kl}^2]-\mathcal{S}^2|^2}{\mathcal{S}^4}\rightarrow 0.
\end{align*}
Using the expression in Step 3, we have,
\begin{align*}
&\E[d_{kl}^2]=\sum_{j,j'=1}^{p}\E[\widetilde{A}_{kl}^{(j)}\widetilde{A}_{kl}^{(j')}\widetilde{B}_{kl}^2]=(1+O(n^{-1}))\mathcal{S}^2+C_1n^{-1}\E[\widetilde{U}(X,X')^4]^{1/2},
\end{align*}
which implies $\frac{|\E[d_{kl}^2]-\mathcal{S}^2|^2}{\mathcal{S}^4}\rightarrow 0$. On the other hand, we note that
\begin{align}
&\frac{\text{var}(\sum_{k<l}d_{kl}^2)}{n^4\mathcal{S}^4} \nonumber
\\=&
\frac{\sum_{k<l}\sum_{k'<l'}\text{cov}(d_{kl}^2,d_{k'l'}^2)}{n^4\mathcal{S}^4}
\nonumber
\\=&\frac{2\sum_{k<l<l'}\text{cov}(d_{kl}^2,d_{l'l}^2)+\sum_{k<l}\text{var}(d_{kl}^2)}{n^4\mathcal{S}^4}+\frac{\sum_{k<l,k'<l',\{k,l\}\cap\{k',l'\}=\emptyset}\text{cov}(d_{kl}^2,d_{k'l'}^2)}{n^4\mathcal{S}^4} \nonumber
\\ \leq
&\frac{2\sum_{k<l<l'}\E[d_{kl}^2d_{l'l}^2]+\sum_{k<l}\E[d_{kl}^4]}{n^4\mathcal{S}^4}+\frac{\sum_{k<l,k'<l',\{k,l\}\cap\{k',l'\}=\emptyset}\text{cov}(d_{kl}^2,d_{k'l'}^2)}{n^4\mathcal{S}^4}.
\label{var-eq2}
\end{align}
Following similar arguments in Step 4, we can show that the first term in
(\ref{var-eq2}) converges to zero. To deal with the second term in
(\ref{var-eq2}), we again use the expression in Step 3. Using the
fact that
$\text{cov}(\bar{A}_{kl}^{(j_1)}\bar{A}_{kl}^{(j_2)}\bar{B}_{kl}^2,\bar{A}_{k'l'}^{(j_3)}\bar{A}_{k'l'}^{(j_4)}\bar{B}_{k'l'}^2)=0$
for any $1\leq j_1,j_2,j_3,j_4\leq p$ and
$\{k,l\}\cap\{k',l'\}=\emptyset$, and the H\"{o}lder's inequality,
we know
$\sum_{k<l,k'<l',\{k,l\}\cap\{k',l'\}=\emptyset}\text{cov}(d_{kl}^2,d_{k'l'}^2)$
is bounded by $C_1n^{3}\E[\widetilde{U}(X,X')^4]$ for some
sufficiently large constant $C_1$. Therefore by (\ref{eq:con2-2}),
we have
\begin{align*}
\frac{\sum_{k<l,k'<l',\{k,l\}\cap\{k',l'\}=\emptyset}\text{cov}(d_{kl}^2,d_{k'l'}^2)}{n^4\mathcal{S}^4}\leq
\frac{C_2\E[\widetilde{U}(X,X')^4]}{n(\sum_{j,j'=1}^pdcov(x_j,x_{j'})^2)^2}\rightarrow 0.
\end{align*}
This completes the proof.

\subsection{Characterization of local alternative models}\label{sec:local}
We provide some discussions on the local alternative models. Let
$\I=\sqrt{-1}$ and $w(t)=1/(\pi t^2)$. Define
$G_j(u_j,y;t)=(f_{x_j}(t)-e^{\I tu_j})(\mu-y)$ for $u_j, y, t \in
\mathbb{R}$. By Lemma 1 of Sz\'{e}kely et al. (2007), it can be
shown that
\begin{align}\label{ch}
H(z,z')=\int_{\mathbb{R}}\sum^{p}_{j=1}G_j(u_j,y;t)\overline{G_j(u_j',y';t)}w(t)dt,
\end{align}
where $z=(u_1,\dots,u_p,y)^T$ and $z'=(u_1',\dots,u_p',y')^T$.
Recall that $\widetilde{\mathcal{L}}(x,y)=\E [\widetilde{U}(x,\X)V(y,\Y)]$. To have a better understanding of our local alternative model, we shall show that the two conditions
\begin{align}
&\var(\widetilde{\mathcal{L}}(X,Y))=o(n^{-1}\mathcal{S}^2), \label{c1} \\
&\var(\widetilde{\mathcal{L}}(X,Y'))=o(\mathcal{S}^2),  \label{c2}
\end{align}
have a similar interpretation as equation (4.2) in Zhong and Chen (2011). %As we have discussed under the null, $\mathcal{S}^2=O(\E H^2(Z,Z'))=O(\text{Tr}^2(\Phi))$.
Using representation (\ref{ch}), we have
\begin{align*}
\widetilde{\mathcal{L}}(x,y)=\int_{\mathbb{R}}\sum^{p}_{j=1}G_j(x_j,y;t)\overline{\E [G_j(x_j',Y';t)]}w(t)dt.
\end{align*}
Thus (\ref{c1}) can be re-expressed as
\begin{align*}
\var(\widetilde{\mathcal{L}}(X,Y))=&\int_{\mathbb{R}}\int_{\mathbb{R}}\sum^{p}_{j,k=1}\E\{(G_j(x_j,Y;t_1)-\E G_j(x_j,Y;t_1))\overline{(G_k(x_k,Y;t_2)-\E G_k(x_k,Y;t_2))}\}
\\&\E [G_k(x_k',Y';t_2)]\overline{\E [G_j(x_j',Y';t_1)]} w(t_1)w(t_2)dt_1dt_2=O(n^{-1}\mathcal{S}^2),
\end{align*}
and (\ref{c2}) can be written as
\begin{align*}
\var(\widetilde{\mathcal{L}}(X,Y'))=&\int_{\mathbb{R}}\int_{\mathbb{R}}\sum^{p}_{j,k=1}\E\{G_j(x_j,Y';t_1)\overline{G_k(x_k,Y';t_2)}\}
\\&\E [G_k(x_k',Y';t_2)]\overline{\E [G_j(x_j',Y';t_1)]} w(t_1)w(t_2)dt_1dt_2
\\=&\int_{\mathbb{R}}\int_{\mathbb{R}}\sum^{p}_{j,k=1}\E\{(f_{x_j}(t_1)-e^{\I t_1 u_j})(f_{x_k}(-t_2)-e^{-\I t_2 u_k})\}\var(Y)
\\&\E [G_k(x_k',Y';t_2)]\overline{\E [G_j(x_j',Y';t_1)]} w(t_1)w(t_2)dt_1dt_2.
\end{align*}
By the definition $H^*(Z,Z')=H(Z,Z')-\widetilde{\mathcal{L}}(X,Y)-\widetilde{\mathcal{L}}(X',Y')+\E[\widetilde{U}(X,X')V(Y,Y')]$ and $\mathcal{S}^2=O(\var(H^*(Z,Z')))$.
Also note that $H^*$ can be re-expressed as
\begin{align}\label{hstar}
H^*(Z,Z')=\int_{\mathbb{R}}\sum^{p}_{j=1}\{G_j(x_j,Y;t)-\E G_j(x_j,Y;t)\}\overline{\{G_j(x_j',Y';t)-\E [G_j(x_j',Y';t)]\}}w(t)dt.
\end{align}
We can define an operator $\Phi^*$ in a similar fashion as $\Phi$ by
replacing $G_j(x_j,Y;t)$ with its demeaned version $G_j(x_j,Y;t)-\E
[G_j(x_j,Y;t)]$. Denote by $\Phi_0$ with
$\E[G_j(x_j,Y;t)\overline{G_i(x_i,Y;t')}]$ being replaced by
$\var(Y)(f_{x_j,x_i}(t,-t')-f_{x_j}(t)f_{x_i}(-t'))$ in $\Phi$. Note
that $\Phi=\Phi_0$ when $X$ and $Y$ are independent. Let
$h(t)=(h_1(t),\dots,h_p(t))$ with $h_k(t)=E [G_k(x_k,Y;t)]$. For
$g(t)=(g_1(t),\dots,g_p(t))$ and
$\widetilde{g}(t)=(\widetilde{g}_1(t),\dots,\widetilde{g}_p(t))$,
define the inner product,
$$<g,\widetilde{g}>_w=\sum^{p}_{j=1}\int_{\mathbb{R}}w(t)g_j(t)\overline{\widetilde{g}_j(t)}dt.$$
Thus Conditions (\ref{c1})-(\ref{c2}) become
\begin{align}
&<\Phi^*(h),h>_w=O(n^{-1}Tr((\Phi^*)^2)), \label{cc1}
\\&<\Phi_0(h),h>_w=O(Tr((\Phi^*)^2)). \label{cc2}
\end{align}
Notice that $h(t)=\mathbf{0}_{p\times 1}$ for all $t$ under the null. Conditions (\ref{cc1})-(\ref{cc2})
quantify the distance between the alternatives and the null
hypothesis. The characterization of the local alternative model is somewhat abstract but this is sensible because MDD targets at very broad alternatives, i.e., arbitrary type of conditional mean dependence.

\subsection{Asymptotic analysis for conditional quantile dependence testing}
In this subsection, we prove Proposition 3.1.
We first state the following lemma whose proof is given in Shao and Zhang
(2014).

%\begin{lemma}
%Suppose $F_{Y,n}(\hat{Q}_{\tau})=\tau$, where $F_{Y,n}$ is the
%empirical distribution function based on $\{Y_i\}^{n}_{i=1}$. Then
%$\hat{B}^*_{Q,kl}=$
%\end{lemma}

\begin{lemma}\label{lemma-quan1}
Suppose the distribution function of $Y$ satisfies Assumption 3.1, then there exist $\epsilon_0>0$ and $c>0$ such that for
$\epsilon\in (0,\epsilon_0)$,
$$P\left(\frac{1}{n}\sum^{n}_{k=1}\left|\W_k-W_k\right|> \epsilon\right)\leq 3\exp(-2nc\epsilon^2).$$
\end{lemma}
Let $B_{Q,kl}=|W_{k}-W_{l}|^2/2$ and $B^*_{Q,kl}$ be
its $\mathcal{U}$-centered version. Note that
\begin{align*}
MDD_n(\W|x_j)^2-MDD_n(W|x_j)^2=\frac{1}{n(n-3)}\sum_{k\neq
l}\widetilde{A}_{kl}(j)(\hat{B}_{Q,kl}^*-B_{Q,kl}^*).
\end{align*}
Thus we need to show that
\begin{align*}
\frac{1}{\mathcal{S}_Q}\sum_{k\neq
l}(\hat{B}_{Q,kl}^*-B_{Q,kl}^*)\sum_{j=1}^{p}\widetilde{A}_{kl}(j)=o_p(n).
\end{align*}
Direct calculation yields that
\begin{align*}
\hat{B}^*_{Q,kl}=&\frac{-1}{n-2}(\W_k^2+\W_l^2)-\W_k\W_l+\frac{1}{(n-1)(n-2)}\sum^{n}_{j=1} \W_j^2
\\&+\frac{1}{n-2}(\W_k+\W_l)\sum_{j=1}^{n}\W_j-\frac{1}{(n-1)(n-2)}\left(\sum_{j=1}^{n}\W_j\right)^2,
\end{align*}
and a similar decomposition for $B_{Q,kl}^*.$ Hence we have
\begin{align*}
&\frac{1}{\mathcal{S}_Q}\sum_{k\neq
l}(\hat{B}_{Q,kl}^*-B_{Q,kl}^*)\sum_{j=1}^{p}\widetilde{A}_{kl}(j)
\\=&\frac{-2}{\mathcal{S}_Q(n-2)}\sum_{k\neq
l}(\W_k^2-W_k^2)\sum_{j=1}^{p}\widetilde{A}_{kl}(j)+\frac{1}{\mathcal{S}_Q(n-1)(n-2)}\sum^{n}_{i=1}(\W_i^2-W_i^2)\sum_{k\neq
l}\sum_{j=1}^{p}\widetilde{A}_{kl}(j)
\\&-\frac{1}{\mathcal{S}_Q}\sum_{k\neq
l}(\W_k\W_l-W_kW_l)\sum_{j=1}^{p}\widetilde{A}_{kl}(j)+\frac{2}{(n-2)S_Q}\sum_{k\neq l}\left(\W_k\sum_{i=1}^{n}\W_i-W_k\sum_{i=1}^{n}W_i\right)\sum_{j=1}^{p}\widetilde{A}_{kl}(j)
\\&-\frac{1}{(n-1)(n-2)S_Q}\left\{\left(\sum_{i=1}^{n}\W_i\right)^2-\left(\sum_{i=1}^{n}W_i\right)^2\right\}\sum_{k\neq l}\sum_{j=1}^{p}\widetilde{A}_{kl}(j)
\\=&J_{1,n}+J_{2,n}+J_{3,n}+J_{4,n}+J_{5,n},
\end{align*}
where $J_{i,n}$ with $1\leq i\leq 5$ are defined implicity. Notice that
\begin{align*}
|J_{1,n}|^2\leq &\frac{4}{\mathcal{S}_Q^2(n-2)^2}\sum_{k\neq
l}(\W_k^2-W_k^2)^2\sum_{k\neq l}\left(\sum_{j=1}^{p}\widetilde{A}_{kl}(j)\right)^2\\
\leq &\frac{C(n-1)}{\mathcal{S}_Q^2(n-2)^2}\sum_{k=1}^n|\W_k-W_k|\sum_{k\neq l}\left(\sum_{j=1}^{p}\widetilde{A}_{kl}(j)\right)^2,
\end{align*}
and
\begin{align*}
|J_{2,n}|\leq &\frac{C'}{\mathcal{S}_Q(n-1)(n-2)}\sum^{n}_{i=1}|\W_i-W_i|\left|\sum_{k\neq
l}\sum_{j=1}^{p}\widetilde{A}_{kl}(j)\right| \\
\leq & \frac{C'n}{\mathcal{S}_Q(n-1)(n-2)}\sum^{n}_{i=1}|\W_i-W_i|\left\{\sum_{k\neq
l}\left(\sum_{j=1}^{p}\widetilde{A}_{kl}(j)\right)^2\right\}^{1/2},
\end{align*}
for some constants $C,C'>0.$ By Lemma \ref{lemma-quan1}, we need to show that
$$\frac{1}{n^2\mathcal{S}_Q^2}\E\sum_{k\neq l}\left(\sum_{j=1}^{p}\widetilde{A}_{kl}(j)\right)^2=O(1).$$
By Sz\'{e}kely and Rizzoz (2014) and (11) in the paper, we have
\begin{align*}
\E\left[\sum_{k\neq
l}\left(\sum_{j=1}^{p}\widetilde{A}_{kl}(j)\right)^2\right]=&\sum^{p}_{j,j'=1}\E\left(\sum_{k\neq
l}\widetilde{A}_{kl}(j)\widetilde{A}_{kl}(j')\right)
=n(n-3)\sum^{p}_{j,j'=1}dcov(x_j,x_{j'})^2=O(n^2\mathcal{S}_Q^2).
\end{align*}
Then we have $J_{1,n}=o_p(n)$ and $J_{2,n}=o_p(n)$. Recall that $Q_{\tau}=Q_{\tau}(Y)$
is the $\tau$th quantile of $Y$. For $J_{3,n}$, we note that
\begin{align*}
J_{3,n}=&-\frac{1}{\mathcal{S}_Q}\sum_{k\neq
l}(\W_k\W_l-W_kW_l)\sum_{j=1}^{p}\widetilde{A}_{kl}(j)
\\=&-\frac{1}{\mathcal{S}_Q}\sum_{k\neq
l}(\W_k-W_k)\W_l\sum_{j=1}^{p}\widetilde{A}_{kl}(j)-\frac{1}{\mathcal{S}_Q}\sum_{k\neq
l}(\W_l-W_l)W_k\sum_{j=1}^{p}\widetilde{A}_{kl}(j)
\\=&-\frac{1}{\mathcal{S}_Q}\sum_{k\neq
l}\mathbf{1}\{Q_{\tau}<Y_k\leq \hat{Q}_{\tau}\}\W_l\sum_{j=1}^{p}\widetilde{A}_{kl}(j)+\frac{1}{\mathcal{S}_Q}\sum_{k\neq
l}\mathbf{1}\{\hat{Q}_{\tau}<Y_k\leq Q_{\tau}\}\W_l\sum_{j=1}^{p}\widetilde{A}_{kl}(j)
\\&-\frac{1}{\mathcal{S}_Q}\sum_{k\neq
l}\mathbf{1}\{Q_{\tau}<Y_l\leq \hat{Q}_{\tau}\}W_k\sum_{j=1}^{p}\widetilde{A}_{kl}(j)+\frac{1}{\mathcal{S}_Q}\sum_{k\neq
l}\mathbf{1}\{\hat{Q}_{\tau}<Y_l\leq Q_{\tau}\}W_k\sum_{j=1}^{p}\widetilde{A}_{kl}(j)
\\=&I_{1,n}+I_{2,n}+I_{3,n}+I_{4,n}.
\end{align*}
Let
$\mathcal{A}_{\delta}=\{\hat{Q}_{\tau}-Q_{\tau}>\delta\}\cup \{\hat{Q}_{\tau}-Q_{\tau}\leq -\delta\}$.
For any $\epsilon>0,$ choose $0<\delta<\delta_0$ such that
$P(\mathcal{A}_{\delta})\leq \epsilon/2$ for large enough $n.$
We claim that that for small enough $\delta$, and $C_i=Q_{\tau}$ or $\hat{Q}_{\tau}$,
\begin{equation}\label{eq:ass1}
\begin{split}
&\E[\mathbf{1}\{Y_{j_1}\leq C_1,Y_{j_2}\leq C_2,\dots,Y_{j_k}\leq C_k, -\delta<\hat{Q}_{\tau}-Q_{\tau}\leq \delta\}|X_1,\dots,X_n]
\\=&\E[\mathbf{1}\{Y_{j_1}\leq C_1,Y_{j_2}\leq C_2,\dots,Y_{j_k}\leq C_k,-\delta<\hat{Q}_{\tau}-Q_{\tau}\leq \delta\}],
\end{split}
\end{equation}
where $k=2,3,4$. The proof of (\ref{eq:ass1}) is given at the end. By (\ref{eq:ass1}), we have for any $\varepsilon>0,$
\begin{align*}
&P(|I_{1,n}/n|>\varepsilon)
\\ \leq &
P(|I_{1,n}/n|>\varepsilon,\mathcal{A}_{\delta}^c)+P(\mathcal{A}_{\delta})\leq
P(|I_{1,n}/n|>\varepsilon,\mathcal{A}_{\delta}^c)+\epsilon/2
\\ \leq & P\left(\left|\frac{1}{n\mathcal{S}_Q}\sum_{k\neq
l}\mathbf{1}\{Q_{\tau}<Y_k\leq \hat{Q}_{\tau},\mathcal{A}_{\delta}^c\}\W_l\sum_{j=1}^{p}\widetilde{A}_{kl}(j)\right|>\varepsilon\right)+\epsilon/2
\\ \leq & \frac{1}{\varepsilon^2n^2\mathcal{S}_Q^2}\E\left|\sum_{k\neq
l}\mathbf{1}\{Q_{\tau}<Y_k\leq \hat{Q}_{\tau},\mathcal{A}_{\delta}^c\}\W_l\sum_{j=1}^{p}\widetilde{A}_{kl}(j)\right|^2+\epsilon/2
\\ \leq & \frac{1}{\varepsilon^2n^2\mathcal{S}_Q^2}\sum_{\{k,l\}\neq \{k',l'\},k\neq l,k'\neq l'}
\E[\mathbf{1}\{Q_{\tau}<Y_k\leq \hat{Q}_{\tau},Q_{\tau}<Y_{k'}\leq \hat{Q}_{\tau},\mathcal{A}_{\delta}^c\}]\\
&\times\W_l\W_{l'}\sum_{j,j'=1}^{p}\E[\widetilde{A}_{kl}(j)\widetilde{A}_{k'l'}(j')]+\frac{C \delta}{\varepsilon^2n^2\mathcal{S}_Q^2}\E\sum_{k\neq
l}\left(\sum_{j=1}^{p}\widetilde{A}_{kl}(j)\right)^2+\epsilon/2,
\\ \leq &\frac{C\delta}{\varepsilon^2n^2\mathcal{S}_Q^2}\sum_{\{k,l\}\neq \{k',l'\},k\neq l,k'\neq l'}\left|\sum_{j,j'=1}^{p}\E\widetilde{A}_{kl}(j)\widetilde{A}_{k'l'}(j')\right|
+\frac{C \delta}{\varepsilon^2n^2\mathcal{S}_Q^2}\E\sum_{k\neq
l}\left(\sum_{j=1}^{p}\widetilde{A}_{kl}(j)\right)^2+\epsilon/2,
\end{align*}
for some $C>0$, where we have used the fact that
\begin{align*}
&\E[\mathbf{1}\{Q_{\tau}<Y_k\leq \hat{Q}_{\tau},Q_{\tau}<Y_{k'}\leq \hat{Q}_{\tau},\mathcal{A}_{\delta}^c\}\W_l\W_{l'}]
\\ \leq& C\E [\mathbf{1}\{Q_{\tau}<Y_k\leq Q_{\tau}+\delta\}]\leq CG_2(\delta_0)\delta.
\end{align*}
Define the following quantities
\begin{align*}
&d_{1}(j)=\E[|x_{j}-x_j'|],\quad d_{2}(j,j')=\E[|x_j-x_j'||x_{j'}-x_{j'}''|],\quad d_{3}(j,j')=\E[|x_j-x_{j}'||x_{j'}-x_{j'}'|].
\end{align*}
We have for $\{k,l\}\cap \{k',l'\}=\emptyset$ and any $1\leq j,j'\leq p$,
\begin{align*}
&\E[\widetilde{A}_{kl}(j)\widetilde{A}_{k'l'}(j')]
\\=&d_1(j)d_1(j')-\frac{4}{n-2}\left\{2d_2(j,j')+(n-3)d_1(j)d_1(j')\right\}
\\&-\frac{n}{(n-1)(n-2)^2}\left\{(n-2)(n-3)d_1(j)d_1(j')+2d_3(j,j')+4(n-2)d_2(j,j')\right\}
\\&+\frac{4}{(n-2)^2}\left\{(n-2)(n-3)d_1(j)d_1(j')+3(n-2)d_2(j,j')+d_3(j,j')\right\}
\\=&\frac{2}{(n-1)(n-2)}(d_1(j)d_1(j')-2d_2(j,j')+d_3(j,j'))=\frac{2}{(n-1)(n-2)}dcov(x_j,x_{j'})^2.
\end{align*}
Notice that in the summation over $\{k,l\}\neq \{k',l'\},k\neq l,k'\neq l'$, we have $O(n^4)$ such terms.
Using similar calculation, we have for $l\neq l'$, $l\neq k$ and $l'\neq k$,
\begin{align*}
&\E[\widetilde{A}_{kl}(j)\widetilde{A}_{kl'}(j')]
\\=&d_2(j,j')-\frac{n-3}{(n-2)^2}\left\{(n-2)d_2(j,j')+d_3(j,j')\right\}-\frac{1}{n-2}\left\{2d_2(j,j')+(n-3)d_1(j)d_1(j')\right\}
\\&+\frac{3}{(n-2)^2}\left\{(n-2)(n-3)d_1(j)d_1(j')+3(n-2)d_2(j,j')+d_3(j,j')\right\}
\\&-\frac{n}{(n-1)(n-2)^2}\left\{(n-2)(n-3)d_1(j)d_1(j')+2d_3(j,j')+4(n-2)d_2(j,j')\right\}
\\&-\frac{1}{n-2}\left\{(n-3)d_1(j)d_1(j')+2d_2(j,j')\right\}
\\=&-\frac{(n-3)}{(n-1)(n-2)}(d_1(j)d_1(j')-2d_2(j,j')+d_3(j,j'))=-\frac{(n-3)}{(n-1)(n-2)}dcov(x_j,x_{j'})^2.
\end{align*}
And we have $O(n^3)$ such terms in the summation over $\{k,l\}\neq \{k',l'\},k\neq l,k'\neq l'$.
Thus we deduce that
\begin{align*}
\sum_{\{k,l\}\neq \{k',l'\},k\neq l,k'\neq l'}\sum_{j,j'=1}^{p}\left|\E\widetilde{A}_{kl}(j)\widetilde{A}_{k'l'}(j')\right|\leq Cn^2\sum^{p}_{j,j'=1}dcov(x_j,x_{j'})^2.
\end{align*}
By (11) in the paper, $W$ is independent of $X$ and thus
$$\frac{\sum^{p}_{j,j'=1}dcov(x_j,x_{j'})^2}{S_Q^2}=\tau^{-2}(1-\tau)^{-2}.$$
Hence we can choose a small enough $\delta$ such that
$$P(|I_{1,n}/n|>\varepsilon)\leq \epsilon,$$
which suggests that $I_{1,n}=o_p(n)$. Similarly we have $I_{j,n}=o_p(n)$ for $2\leq j\leq 4,$ which implies that $J_{3,n}=o_p(n).$
Using similar arguments, we can show that $J_{4,n}=o_p(n)$. Finally note that by Lemma \ref{lemma-quan1},
$$J_{5,n}/n\leq o_p(1)\frac{1}{nS_Q}\left|\sum_{k\neq l}\sum_{j=1}^{p}\widetilde{A}_{kl}(j)\right|.$$
Because
$$\E\frac{1}{n^2S_Q^2}\left|\sum_{k\neq l}\sum_{j=1}^{p}\widetilde{A}_{kl}(j)\right|^2=O(1),$$
we get $J_{5,n}=o_p(n)$.

Finally we prove (\ref{eq:ass1}). Note that $\hat{Q}_{\tau}=\inf\{y: F_n(y)\geq \tau\}$, where $F_n$ is the empirical distribution function based on $\{Y_i\}^{n}_{i=1}$.
Thus
$-\delta<\hat{Q}_{\tau}-Q_{\tau}\leq \delta$ is equivalent to
\begin{align*}
&F_n(Q_{\tau}+\delta)=\frac{1}{n}\sum^n_{i=1}\mathbf{1}\{Y_i\leq Q_{\tau}+\delta\}\geq \tau,\\
&F_n(Q_{\tau}-\delta)=\frac{1}{n}\sum^n_{i=1}\mathbf{1}\{Y_i\leq Q_{\tau}-\delta\}< \tau.
\end{align*}
We see that the event $\mathbf{1}\{Y_{j_1}\leq C_1,Y_{j_2}\leq C_2,\dots,Y_{j_k}\leq C_k, -\delta<\hat{Q}_{\tau}-Q_{\tau}\leq \delta\}$
only depends on $\mathbf{1}\{Y_{j_i}\leq C_i\}$ for $1\leq i\leq k$, and $\mathbf{1}\{Y_i\leq Q_{\tau}\pm \delta\}$. Moreover, given that $-\delta<\hat{Q}_{\tau}-Q_{\tau}\leq \delta$, the value of $\hat{Q}_{\tau}$
is determined by $\mathbf{1}\{Y_i\leq Q_{\tau}+a\}$ for $-\delta\leq a\leq \delta$ and $1\leq i\leq n$.
Hence the event $\mathbf{1}\{Y_{j_1}\leq C_1,Y_{j_2}\leq C_2,\dots,Y_{j_k}\leq C_k, -\delta<\hat{Q}_{\tau}-Q_{\tau}\leq \delta\}$ is determined by
$\mathbf{1}\{Y_i\leq Q_{\tau}+a\}$ for $-\delta\leq a\leq \delta$ and $1\leq i\leq n$.
Therefore, (\ref{eq:ass1}) holds under Assumption 3.2. The proof for Proposition 3.1 is thus completed.

\section{Extension to factorial designs}\label{sec:factor}
As motivated by the study in Zhong and Chen (2011), we propose an
extension of our MDD-based test to the situation where the
observation $(X_i, Y_i)$ is not a simple random sample but has a
factorial design structure, as often the case in microarray study.
Following Zhong and Chen (2011), we shall focus on the two way
factorial designs with two factors $A$ and $B$, where $A$ has $I$
levels and $B$ has $J$ levels. In the latter paper, they assumed the
observations $(X_{ijk},Y_{ijk})$ in the $i$th level of $A$ and $j$th
level of $B$ satisfy a linear model:
\[
\mathbb{E}(Y_{ijk}|X_{ijk})=\mu_{ij}+X^T_{ijk}\beta,\quad i=1,\dots,I;\quad j=1,\dots,J,
\]
for $k=1,2,\dots,n_{ij}$, where $X_{ijk}=(x_{ijk,1},\dots,x_{ijk,p})^T$, $n_{ij}$ denotes the number
of observations in cell $(i,j)$, $p$ is the dimension of the covariates, and $\mu_{ij}$ denotes the fixed or
random effect corresponding to the cell $(i,j)$. In Zhong and Chen
(2011), they tested the hypothesis $\tilde{H}_0: \beta=\beta_0$
versus $\tilde{H}_a:\beta \ne \beta_0$. In particular, when
$\beta_0=0$, the null corresponds to the conditional mean
independence of $Y$ given $X$ in each of the cell $(i, j)$ for the
factorial design regardless of the nuisance parameters $\mu_{ij}$.
In this section, we generalize the test in Section 2.2 of the main paper
to the factorial design case.

We introduce some notation first. Let  $(X_{ij},Y_{ij})
\stackrel{D}= ( X_{ijk}, Y_{ijk})$ for $k=1,2,...,n_{ij}$, where $
X_{ij}=(x_{ij,1},...,x_{ij,p})^T$. One way to formulate the test is to consider
the following hypothesis

\[
H''_0: {\mathbb E}[ Y_{ij}| x_{ij,h}] ={\mathbb E}[ Y_{ij}]\quad
\text{almost surely},
\]
for all $1 \le h  \le p$,  $1 \le i \le I$, and $1 \le j \le J$ versus
the alternative that
\[
H''_a: P\left({\mathbb E}[  Y_{ij} |  x_{ij,h}] \ne {\mathbb E}[ Y_{ij}]\right)>0
\]
for some $1 \le h  \le p$,  $1\le i \le I$, and $1 \le j \le J$.

Throughout the discussions, we assume independence across different
cells, i.e.,
\begin{equation}
   \text{$(Y_{ij},X_{ij})$ is independent of $(Y_{i'j'},X_{i'j'})$},
    \label{eq:con0_fac}
\end{equation}
for $ (i, j) \ne (i', j')$.
%This condition is weaker than the one made in Zhong and Chen (2011), where they assumed independence in a factor model structure.
Let $MDD_{n_{ij}}({Y}_{ij}|{x}_{ij,h})^2$ be the unbiased estimator
for $MDD({Y}_{ij}|{x}_{ij,h})^2$ based on the sample
$(Y_{ijk},X_{ijk})_{k=1}^{n_{ij}}$. By the Hoeffding decomposition,
we have that under the $H''_0$,
\begin{align*}
&MDD_{n_{ij}}({Y}_{ij}|{x}_{ij,h})^2
\\&=\frac{1}{\dbinom{n_{ij}}{2}} \sum_{1\le k < l \le n_{ij}}
U_{ij,h}(x_{ijk,h},x_{ijl,h})
V_{ij}(Y_{ijk},Y_{ijl})+R_{ij,h,n_{ij}}
\end{align*}
where $R_{ij,h,n_{ij}}$ is an asymptotically negligible remainder
term. Define
\[
S_{ij}^2=\sum_{h,h'=1}^p{\sigma}_{h,h'}^{ij},\quad
{\sigma}_{h,h'}^{ij}=E [U_{ij,h}(x_{ij,h},x'_{ij,h})
U_{ij,h'}(x_{ij,h'},x'_{ij,h'}) V_{ij}^2(Y_{ij},Y'_{ij})],
\]
where $(Y_{ij}',{X}_{ij}')$ is an independent copy of
$(Y_{ij},X_{ij})$. A natural estimator for ${\sigma}_{h,h'}^{ij}$ is
given by
\[
\hat \sigma_{h,h'}^{ij}=\frac{1}{\dbinom{n_{ij}}{2}} \sum_{1\le k <
l \le n_{ij}} {\widetilde  A}_{kl}^{ij}(h){\widetilde
A}_{kl}^{ij}(h')({\widetilde  B}_{kl}^{ij})^2
\]
with  ${\widetilde  A}_{kl}^{ij}(h)$ and $ {\widetilde  B}_{kl}^{ij} $
being the ${\cal U}$-centered versions of
$A_{kl}^{ij}(h)=|x_{ijk,h}-x_{ijl,h}|$ and ${\widetilde
B}_{kl}^{ij}=|Y_{ijk}-Y_{ijl}|^2/2$ respectively. Denote $\hat \Sigma_{ij}=(\hat
\sigma_{h,h'}^{ij})_{h,h'=1}^p$ and $\hat {\cal S}^2_{ij}={\mathbf
1}'_p \hat \Sigma_{ij} {\mathbf 1}_p$. The test statistics we
considered for the case of factorial designs are
\[
\breve{T}_{F,n}= \frac{\sum_{i=1}^I \sum_{j=1}^J
\sqrt{\dbinom{n_{ij}}{2}} \sum_{h=1}^p
MDD_{n_{ij}}({Y}_{ij}|{x}_{ij,h})^2 }{\sqrt{\sum_{i=1}^I
\sum_{j=1}^J {\cal S}^2_{ij}}},
\]
and
\[
T_{F,n}= \frac{\sum_{i=1}^I \sum_{j=1}^J \sqrt{\dbinom{n_{ij}}{2}}
\sum_{h=1}^p MDD_{n_{ij}}({Y}_{ij}|{x}_{ij,h})^2
}{\sqrt{\sum_{i=1}^I \sum_{j=1}^J \hat {\cal S}^2_{ij}}}.
\]

Similarly, we can define $Z'_{ij}=({X}_{ij}',Y'_{ij})$,
$Z''_{ij}=({X}_{ij}'',Y''_{ij})$ and
$Z'''_{ij}=({X}_{ij}''',Y'''_{ij})$ to be independent copies of
$Z_{ij}=({X}_{ij},Y_{ij})$. Recall the definitions $
\widetilde{U}_{ij}({X}_{ij},{X}_{ij}')=\sum_{h=1}^p
\widetilde{U}_{ij,h}(x_{ij,h},x'_{ij,h}) $, $ H_{ij}(Z_{ij},Z'_{ij})
= \widetilde{U}_{ij}({X}_{ij},{X}_{ij}') V_{ij}(Y_{ij},Y_{ij}') $
and $ G_{ij}(Z_{ij},Z'_{ij})= \E[ H_{ij}(Z_{ij},Z''_{ij})
H_{ij}(Z'_{ij},Z''_{ij}) | (Z_{ij},Z'_{ij})] $. For each $ (i,j) \in
\{1,...,I\} \times \{1,...,J\} $, we impose the following
conditions:
\begin{equation}
\begin{split}
    &\frac{\E [G_{ij}(Z_{ij},Z'_{ij})^2]}{\{\E[
        H_{ij}(Z_{ij},Z'_{ij})^2]\}^2}\rightarrow 0,
    \\&\frac{\E [H_{ij}(Z_{ij},Z'_{ij})^4]/n+\E [H_{ij}(Z_{ij},Z''_{ij})^2H_{ij}(Z'_{ij},Z''_{ij})^2]}{n\{\E
        [H_{ij}(Z_{ij},Z'_{ij})^2]\}^2}\rightarrow 0, \label{eq:con1_fac}
\end{split}
\end{equation}
and
\begin{eqnarray}
    && \frac{\E [\widetilde{U}_{ij}({X}_{ij},{X}_{ij}'')^2V_{ij}(Y_{ij},Y'_{ij})^2]}{\mathcal{S}_{ij}^2}=O(1), \label{eq-r1-fac} \\
%    && \frac{\E \widetilde{U}_{ij}({X}_{ij},{X}_{ij}'')\widetilde{U}_{ij}(X'_{ij},X''_{ij})V_{ij}(Y_{ij},Y'_{ij})^2}{\mathcal{S}_{ij}^2}=O(1), \label{eq-r2-fac}\\
    &&\frac{\E [\widetilde{U}_{ij}({X}_{ij},{X}_{ij}')^2]\E
        [V_{ij}(Y_{ij},Y'_{ij})^2]}{\mathcal{S}_{ij}^2}=O(1).
    \label{eq-r3-fac}
\end{eqnarray}

\begin{theorem}\label{thm:factor}
    Under the assumption (\ref{eq:con0_fac})-(\ref{eq-r3-fac}), and the null hypothesis $H''_0$, we have
    \[
 \breve{T}_{F,n} \rightarrow^d N(0,1).
    \]
\end{theorem}

Theorem \ref{thm:factor} is readily attained by slightly modifying
the proof of Theorem 2.1 of the main paper. We impose conditions on each cell,
which generalize those in Theorem 2.1 in the paper to the case of
factorial designs.

%\begin{remark}
% { \rm  As in Assumption (\ref{eq:con0_fac}), we assume independence across different cells. When this assumption is violated, one should
%    consider a more sophisticated test statistic which takes the dependence across cells into account. This is beyond the scope of the current paper and it is left for future research.}

    %In this case, suitable conditions should be imposed on
    %   $\sum_{i=1}^I\sum^{J}_{j=1} H_{ij}(Z_{ij},Z_{ij}) $ instead of each $H_{ij}$ as the case in Assumption (\ref{eq:con1_fac}).
%\end{remark}

\section{Additional Simulation Results}
We conduct additional simulations to assess the
finite sample performance of the proposed MDD-based tests.

% % % % % % % % % % % % % % % % % % % % % % % % % %
% % % % % % % % % % % % % % % % % % % % % % % % % %
% % % % % % % % % % % EXAMPLE 1 % % % % % % % % % %
% % % % % % % % % % % % % % % % % % % % % % % % % %
% % % % % % % % % % % % % % % % % % % % % % % % % %

\subsection{Conditional mean independence}

\begin{example}\label{ex21}
{\rm
We consider the simple linear model:
\[
Y_{i}= X_i^T\beta+\epsilon_{i},\quad i=1,2,\dots,n,
\]
where $X_i=(x_{i1},x_{i2},\dots,x_{ip})^T$ is a $p$-dimensional
vector of covariates, $\beta=(\beta_1,\dots,\beta_p)^T$ is the
regression coefficient and $\epsilon_i$ is the error that is
independent of $X_i$. The covariates are generated from the
following model:
\begin{equation}\label{eq:cov}
x_{ij}=(\varsigma_{ij}+\varsigma_{i0})/\sqrt{2},\quad j=1,2,...,p,
\end{equation}
where
$(\varsigma_{i0},\varsigma_{i1},\dots,\varsigma_{ip})^T\overset{i.i.d.}{\sim}
N(0,I_{p+1})$.  Hence the covariates are strongly correlated with
the pairwise correlation equal to 0.5.
} \end{example}

Under the null, $\beta=\mathbf{0}_{p\times 1}$; under both the
sparse and non-sparse alternative, we fix $|\beta|_p=0.06$ as in
Example 4.1 of the paper. We set $n=100$, $p=50,100,200$, and
also consider three configurations for the error $\epsilon$:
$N(0,1)$, $t_3$, and $\chi^2_1-1$. Table \ref{table:eg21} presents
the empirical sizes and powers of the proposed test and the ZC test
for the cases of $p<n$, $p=n$ and $p>n$. The empirical sizes of both
tests are reasonably close to the 10\% nominal level for three
different error distributions. At the 5\% significance level,
however, we see both tests have slightly inflated rejection
probabilities under the null, which is similar to what we observe in
Example 4.1. For the empirical powers, our test is highly
comparable to ZC test under the simple linear model.

% % % % % % % % % % % % % % % % % % % % % % % % % %
% % % % % % % % % % % % % % % % % % % % % % % % % %
% % % % % % % % %EXAMPLE 2 % %% % % % % % % % % % %
% % % % % % % % % % % % % % % % % % % % % % % % % %
% % % % % % % % % % % % % % % % % % % % % % % % % %

\begin{example}\label{ex22}
{\rm
This example specifies another non-linear relationship between $Y$ and $X$. The
model is given by
\begin{equation}\label{eq:ex3}
Y_{i}=\frac{p}{\sqrt{\sum^{p}_{j=1}\beta_j
x_{ij}^2}}+\epsilon_{i},\quad i=1,2,\dots,n,
\end{equation}
where the covariates  $X_i=(x_{i1}, x_{i2},...,x_{ip})^T$ are
generated according to (\ref{eq:cov}). Again we consider three
configurations for the error, namely $N(0,1)$, $t_3$, and $\chi^2-1$
respectively.
}
\end{example}

Under
the non-sparse alternative, we set $\beta_j=1$ for $j=1,2,...,n/2$. For sparse alternative,
$\beta_j=1$ for $j=1,2,...,5$. The configurations for $p$ and $n$
are the same as Example \ref{ex21}.

Table \ref{table:eg22} summarizes the empirical powers for Example
\ref{ex22}. Similar to  Example 4.2 of the paper,  ZC test,
which is designed for linear model, exhibits very little power under
model mis-specification. The powers for ZC test are even below the
nominal level under sparse alternatives. We notice that our proposed
test has gradually decreasing powers as dimension increases under
both sparse and non-sparse alternatives, which might be explained by
the fact that the covariates enter into the denominator of the model
and the signals are weakened by increasing the denominator in
(\ref{eq:ex3}) as the dimension increases. Overall, the power
performance is reasonably good for our proposed model-free test.

% % % % % % % % % % % % % % % % % % % % % % % % % %
% % % % % % % % % % % % % % % % % % % % % % % % % %
% % % % % % % % %EXAMPLE 3 % %% % % % % % % % % % %
% % % % % % % % % % % % % % % % % % % % % % % % % %
% % % % % % % % % % % % % % % % % % % % % % % % % %

\subsection{Conditional quantile independence}

 \begin{example}\label{ex23}
{\rm
 This example considers a simple linear model with
heteroscedesticity:
\[Y_{i}=  X_i^T \beta + \left(1+ X_i^T \beta \right)^2\epsilon_{i}, \quad i=1,2,\dots,n,\]
where $X_i=(x_{i1}, x_{i2},...,x_{ip})^T$ is a $p$ dimensional
vector and $\epsilon_i$ is the error independent of $X_i$. All the
other configurations are the same as Example \ref{ex21}.
}
\end{example}

Table \ref{table:eg23} shows the empirical sizes and powers for
different configurations in Example \ref{ex23}. The sizes are generally
precise at 10\% level and slightly inflated at 5\% level, and they
do not seem to depend on the error distribution much. The empirical
powers apparently depend on the error distribution. For $N(0,1)$,
$t_3$ and Cauchy(0,1) error, they have higher powers at $\tau=0.5$,
$0.75$ than that at $\tau=0.25$. The powers are almost 1 in the
non-sparse alternatives and suffer big reduction under sparse
alternatives. In addition, the powers at $\tau=0.75$ under sparse
case are higher than that at $\tau=0.5$. In contrast, the powers for
$\chi^2_1-1$ error configuration perform the best at $\tau=0.25$ and
$0.75$, while suffer a great power loss at $\tau=0.5$. Moreover, the
powers at $\tau=0.25$ under sparse alternatives are higher than that
at $\tau=0.75$. These phenomenon should be related to the error
distribution. Notice that normal, student's t and Cauchy
distributions are all symmetric; while $\chi^2_1-1$ is skewed. The performance of MDD-based test for conditional mean is very similar to that for Example 4.4.

% % % % % % % % % % % % % % % % % % % % % % % % % %
% % % % % % % % % % % % % % % % % % % % % % % % % %
% % % % % % % % %EXAMPLE 4 % %% % % % % % % % % % %
% % % % % % % % % % % % % % % % % % % % % % % % % %
% % % % % % % % % % % % % % % % % % % % % % % % % %

As we observe some size distortion in Table 4 for
Example 4.4 of the paper, We further apply the wild bootstrap
to approximate the finite sample distribution of $T_{Q,n}$ in the
same way as the studentized bootstrap statistic $T_n^*$ used in the
conditional mean dependence testing. According to
Table~\ref{table:eg54boot}, the bootstrap helps to reduce the size
distortion substantially and the size corresponding to the bootstrap
approximation is fairly accurate. Similar results are observed for
$10\%$ level and are not presented. We also tried the wild bootstrap
for the case of $\tau=0.5$, and nonzero $\beta$s, which falls under
our null hypothesis. There is some reduction in size distortion in
this case (results not shown) but the size of the bootstrap based
test is still quite inflated for large $p$ (say, $p=100,200$). This
again points to the importance of our local quantile independence
assumption (Assumption 3.2), which is violated in this case
and seems crucial to the validity of our normal and bootstrap
approximation.

\subsection{Conditional mean independence under factorial designs}
In this subsection, we evaluate the finite sample performance of the
proposed test for conditional mean independence under factorial
designs in comparison with the ZC test.

\begin{example}\label{ex55}
{\rm
Consider the factorial design with
non-linear models in each cell,
\[
Y_{ijk}=\mu_{ij}+\sqrt{\sum^{p}_{h=1}\beta_h
x_{ijk,h}^2}+\epsilon_{ijk},\quad k=1,2,\dots,n_{ij},
\]
where $(\mu_{11}, \mu_{12}, \mu_{21}, \mu_{22})=(1,3,3,4)$, $n_{11}=n_{12}=n_{21}=n_{22}=n$ and $\epsilon_{ijk}\sim^{i.i.d} N(0,4)$. We
consider two distributions for the covariates $ X_{ijk}$ in each
cell:
\\
Case 1: $ X_{ijk} $ is generated independently from the moving
average model as in Example 4.1 of the paper.
\begin{equation*}
x_{ijk,h}=\alpha_{ij1} z_{ijk,h} +\alpha_{ij2}
z_{ijk,(h+1)}+\cdots+\alpha_{ijT_{ij}} z_{ijk,(h+T_{ij}-1)}+\mu_{ij,h},
\end{equation*}
for $i=1,2$, $j=1,2$, and $k=1,\dots,n_{ij}$. By choosing $(T_{11},
T_{12}, T_{21}, T_{22})=(10,15,20,25)$, the dependence structure in
each cell is different from the others.
\\
Case 2: $x_{ijk,h}=(\varsigma_{ijk,h}+\varsigma_{ijk,0})/\sqrt{2}$, where
$(\varsigma_{ijk,0},\varsigma_{ijk,1},\dots,\varsigma_{ijk,p})^T\overset{i.i.d}{\sim}
N(0,I_{p+1})$ for $i=1,2$, $j=1,2$, $k=1,\dots,n_{ij}$, and
$h=1,\dots,p$.

We consider $n_{ij} = 30, 50, 70$, and $p = 100, 150, 200$. For the sparse and non-sparse alternatives, we keep $|\beta|_p = 0.06$.

}
\end{example}

Table \ref{table:eg55} summarizes the results for both cases. There is slight
size distortion at 5\% significance level. Both tests have satisfactory power in Case 1. The ZC test,
however, suffers from a significant power loss in Case 2 in part due
to model mis-specification. In contrast, we observe that our proposed
test retains the high power under the non-sparse alternative and the
power increases rapidly as sample size increases. Although our test
exhibits less power under sparse alternatives, it still outperforms
ZC test in power by a noticeable amount.

 \begin{example}\label{ex24}
{\rm This example is based on the two factor balanced
design with two levels for each factor which has been considered in
Zhong and Chen (2011),
\[
Y_{ijk}=\mu_{ij}+X_{ijk}^T\beta+\epsilon_{ijk},\quad
k=1,2,\dots,n_{ij},
\]
where $(\mu_{11},\mu_{12},\mu_{21},\mu_{22})=(1,3,3,4)$, and
$n_{11}=n_{12}=n_{21}=n_{22}$. Within each cell, $ X_{ijk} $ is
generated independently from the following moving average model,
\begin{equation*}
x_{ijk,h}=\alpha_{ij1} z_{ijk,h} +\alpha_{ij2}
z_{ijk,(h+1)}+...+\alpha_{ijT_{ij}} z_{ijk,(h+T_{ij}-1)}+\mu_{ij,h},
\end{equation*}
for $i=1,2$, $j=1,2$, and $k=1,\dots,n_{ij}$. By choosing $(T_{11},
T_{12},T_{21},T_{22})=(10,15,20,25)$, the dependence structure in
each cell is different from each other. We also consider two error
distributions for $\epsilon_i$: $N(0,4)$ and centralized gamma
distribution with shape parameter 1 and scale parameter 0.5.
Consider $n_{ij}=30,50,70$, and $p=100,150,200$. For both the sparse and
non-sparse alternatives, we keep $|\beta|_p=0.06$.
}
\end{example}

The simulation results are presented in Table \ref{table:eg24}. The
size phenomenon is similar to what we observed in other examples. The
powers for both tests in the non-sparse and sparse alternatives are
quite high in most cases. We observe that the two tests are very
much comparable in terms of size and power under the simple linear model with factorial design.

\subsection{Further comparison with existing methods}

In this subsection, we further compare our  test for conditional mean independence with some recent ones developed by McKeague and Qian (2015) and the discussants for the latter paper. McKeague and Qian (2015) proposed an adaptive re-sampling test (MQ test, hereafter)  for detecting the existence of significant predictors in a linear regression model.
 They test $H_{10}: \cov(Y,x_j)=0,~j=1,\cdots,p$ versus $H_{11}: \cov(Y,x_j)\not=0,~\mbox{for at least one}~j=1,\cdots,p$. Under finite second moment assumptions for $Y$ and $X$,
 our null hypothesis $H_0'$ implies their $H_{10}$ but there are situations where their null hypothesis holds but $Y$ is marginally conditionally mean dependent on $x_j$ for some $j$, i.e., $H_a'$ holds; see below for an example. MQ's test is based on the estimated marginal regression coefficient of the selected predictor, which has the strongest marginal correlation with the response. It was shown in MQ that their test's limiting distribution over the parameter space has a discontinuity at the origin, which causes the inconsistency of naive bootstrap. MQ proposed a modified bootstrap method that is adaptive to the nonregular behavior of their test statistic by introducing a thresholding parameter, the choice of which is based on another level of bootstrap. While their method is well suited to identify the most significant predictor, it is computationally expensive due to the use of double bootstrap, and its applicability to large $p$ case is theoretically unknown and computationally prohibitive.

 In the discussions of McKeague and Qian (2015), Zhang and Laber (2015) and Chatterjee and Lahiri (2015) have proposed alternative tests that overcome some limitation of MQ test
 either computationally or methodologically. In particular, Zhang and Laber (2015) proposed to use the largest  marginal t-statistic in magnitude to achieve scale-invariance, and use a parametric bootstrap procedure to obtain the critical values. Their test can be implemented much faster than MQ test.   They also proposed an adaptive parametric bootstrap procedure that can be adaptive to unknown level of sparsity. The second test by Zhang and Laber (2015) is, however, a bit more complex to implement and computationally more expensive and it does not seem to lead to a substantial power gain as seen from Table 2 in their paper, so we decided not to include it into our comparison. It is worth mentioning that the simpler test procedure proposed in Zhang
 and Laber (2015) is based on direction estimation of the asymptotic variance in Theorem 1 of McKeague and Qian (2015), which requires stronger independence (rather than uncorrelatedness)
 conditions of the regression error with covariates.
 Chatterjee and Lahiri (2015) proposed a $L_2$ type test statistic via aggregating the marginal $t$-statistic using $L_2$ norm, and the limiting null distribution can be well approximated by the naive bootstrap. Hence their test does not suffer the non-continuity issue, does not require the selection of a tuning parameter, and is computationally simple.

 Following the suggestion of associate editor and a referee, we shall compare the performance of our test (MDD), with MQ test, ZL test (Zhang and Laber's parametric bootstrap test corresponding to $\widehat{\xi}_n$ in their paper), CL test (Chatterjee and Lahiri 2015). We do not aim to provide a comprehensive simulation study to compare these four tests as the latter three tests have been partially compared in Zhang and Lahiri (2015), and Chatterjee and Lahiri (2015). Also our MDD-based test targets a different null hypothesis so in some situations our test is not directly comparable to the other three. This point will be highlighted via a numerical example below.

\begin{example}
\label{eg:ar}
Consider the following linear model,
$$Y_i = X_i^T \beta + \epsilon_i, ~~~ i=1,...,n$$
where $X_i = (x_{i1},x_{i2},...,x_{ip})^T$ is generated from multivariate Gaussian with zero means and covariance  $\Sigma=(\sigma_{ij})_{i,j=1}^p$, where $\sigma_{ij}=\rho^{|i-j|} $; we consider four cases $\rho = 0.1, 0.5, 0.8, -0.5$ respectively.
Here $\beta=(\beta_1,...,\beta_p)^T$ is the regression coefficient and $\epsilon_i$ is the error that is independent of $X_i$ and also generated from i.i.d $N(0,1)$.
\end{example}

\begin{example}
\label{eg:indep}
Consider the following two models,
\begin{itemize}
\item [i)] $Y_i = X_i^T \beta + \epsilon_i$, $X_i = (x_{i1},x_{i2},...,x_{ip})^T$ is generated from  i.i.d  N(0,1);
\item [ii)] $Y_i = g(X_i)^T \beta + \epsilon_i$, where $g(x)=(x_1^2,...,x_p^2)^T$, $X_i$ is the same as in case i);
\end{itemize}
for $ i = 1,2,...,n$. Here $\beta=(\beta_1,...,\beta_p)^T$ is the regression coefficient and $\epsilon_i$ is the error that is independent of $X_i$ and also generated from i.i.d $N(0,1)$.
\end{example}

Following the setting in Chatterjee and Lahiri (2015), we consider three scenarios: under $H_0$,  $\beta =\boldsymbol{0}_{p \times 1}$; under non-sparse $H_a$,
$\beta =  p^{-1/2} \cdot  c \cdot  \boldsymbol{1}_{p \times 1}$; under sparse $H_a$, $\beta = (c,0,0,...,0)$. We choose $c=0.2$ for Example \ref{eg:ar} and $c=0.5$ for Example \ref{eg:indep}; we fix sample size $n=200$ and $p=10,50,150$. For MQ test, we use fixed thresholds $\lambda_n=(2,4,4.5,5,10)$  as the use of double bootstrap  to select thresholding parameter is too expensive for $p=50, 150$ in our simulations. It is worth noting that in Example~\ref{eg:indep}(ii),  $cov(Y, x_i) = 0$ for $i =1,...,p$, but $\E(Y|x_i) \ne \E(Y) $. Therefore, this model falls under the alternative hypothesis for our proposed test, but under the null for ZL, CL and MQ tests.

From Table \ref{table:ar}, we observe that all tests (mdd, ZL, CL and MQ with $\lambda_n=5, 10$) have quite reasonable sizes  for different $\rho$s; mdd and CL tests have higher power under the dense alternative since they are both $L_2$ type statistics, while the other two $L_\infty$ tests: ZL and MQ test perform better under sparse alternatives. The power generally increases as the correlation $\rho$ increases and the positive dependence within the covariates seems to enhance the power under both  dense and sparse alternatives. When $\rho=-0.5$, none of the tests have good power under the dense alternative.  Notice that MQ's performance is fairly sensitive to the threshold value used, thus  a double bootstrap approach will hopefully provide more stable size/power at the expense of heavy computation. Overall, the performance of mdd is on a par with that of CL, and the performance of ZL and MQ with thresholding parameter $\lambda_n=5,10$ are comparable.
% provide a more reasonable choices, but decreases the usefulness in practice due to the computational cost.

Table \ref{table:indep} further compares the four tests  when the covariates are independent. It appears that the size accuracy is similar to that reported for Example \ref{eg:ar}.
For Case i), the powers of mdd and CL tests are highest under dense alternatives as expected, whereas the powers of ZL and MQ are higher under sparse alternatives.  It can be seen that  the power of mdd test under dense and sparse alternatives are roughly the same for both i) and ii), which can be explained by the fact that for independent covariates, the power of our MDD-based test is proportional to the signal strength $|\beta|$, which is fixed at $|\beta|=0.5$ here. Under the non-linear model in case ii), the power of mdd test is satisfactory for small $p$ and it decreases as the dimension increases. By contrast, the ZL and CL tests exhibit certain size inflation when $p=10$.  Again MQ's size is sensitive to the  thresholding value used and seems accurate for $\lambda_n=5,10$. This example highlights the fact that our mdd test targets marginal conditional mean dependence, whereas
MQ, ZL and CL are developed for detecting marginal uncorrelatedness. Overall, the results demonstrate that our proposed test is highly competitive in the case when mdd can be compared to MQ test and variations (see Example~\ref{eg:ar}),  and when they are not directly comparable (see Example~\ref{eg:indep} ii), the mdd test is able to detect marginal conditional mean dependence and perform reasonably well.

%   to the existing methods under linear model and also remains powerful in detecting non-linear conditional mean dependence.

\section{Data illustration}\label{sec:data}
We apply the proposed test to the data set described in Lkhagvadorj
et al. (2009) and Zhong and Chen (2011). The data is from clinical
outcomes in a randomized factorial design experiment, where 24
six-month-old Yorkshire gilts from a line selected for high feed
efficiency were used. All the gilts were genotyped for the
melanocortin-4 receptor gene (MC4R) variant at position 298; and 12
gilt homozygous for N298 and 12 gilt homozygous for D298. Two feed
treatments are randomly assigned to each group of the MC4R genotype.
One is ad libitum a crude protein standard swine diet, the other one
is fasting diet, which leads to decreased body weight, backfat, and
serum urea concentration and increased serum non-esterified fatty
acid. The genotype and feed treatments are the two factors in the
randomized complete factorial design. A total of six pigs were used
for each combination of genotype and feed treatment across the four
blocks. The goal of this study is to identify conditional mean
dependence of triiodothyronine ($T_3$) measurements on gene sets, as
described in Zhong and Chen (2011), where $T_3$ is a vital thyroid
hormone that increases the metabolic rate, protein synthesis, and
stimulates breakdown of cholesterol.

The gene sets mentioned above are defined by the Gene Ontology (GO
term). The dataset includes the gene expression values for 24,123
genes in the gilts' liver tissues. These genes are then classified
into different gene sets (GO term) according to their biological
functions among three categories: cellular component, molecular
function and biological process. The dataset included 6176 GO terms
in total, where each of them contains some of the genes from the
24,123 genes collected. Our response is the $T_3$ measurements in
the blood. We aim to find the gene sets (GO terms) that have an
impact on the $T_3$ measurements in terms of conditional mean after
accounting for the two factors we have in the design.

We use $i$, $j$, $k$ to denote the indices for feed treatment, MC4R
genotype and observations, respectively. Our response is $Y_{ijk}$,
the $T_3$ measurement for the $k$th pig in $i$th feed treatment and
$j$th genotype. We want to test whether the $l$th GO term, denoted
as $X_{ijk}^l$, contributes to the conditional mean of the $T_3$
measurement or not. The following four different factorial designs
have been considered in Zhong and Chen (2011),
\begin{eqnarray*}
    \text{Design I:} &&  \E(Y_k)=\alpha+(X_k^{l})^T\beta, ~~~~ k=1,...,24;\\
    \text{Design II:} &&  \E(Y_{ik})=\alpha+\mu_i+(X_{ik}^{l})^T\beta, ~~~~ k=1,...,12;\\
    \text{Design III:} &&  \E(Y_{jk})=\alpha+\gamma_j+(X_{jk}^{l})^T\beta, ~~~~ k=1,...,12;\\
    \text{Design IV:} &&  \E(Y_{ijk})=\alpha+\mu_i+\gamma_j+(\mu\gamma)_{ij}+(X_{ijk}^{l})^T\beta, ~~~~ k=1,...,6;
\end{eqnarray*}
for $i=1,2$, $j=1,2$, and $l=1,..,6176$, which correspond to 6176
total gene sets (GO terms). To avoid the linear model assumption, we
test instead the following hypothesis for each of the designs above
\begin{eqnarray*}
    \text{Design I:} && \E(Y_k|X_k^{l})= \E(Y_k), ~~~~ k=1,...,24;\\
    \text{Design II:} && \E(Y_{ik}|X_{ik}^{l})= \E(Y_{ik}), ~~~~ k=1,...,12;\\
    \text{Design III:} && \E(Y_{jk}|X_{jk}^{l})= \E(Y_{jk}), ~~~~ k=1,...,12;\\
    \text{Design IV:} && \E(Y_{ijk}|X_{ijk}^{l}) =\E(Y_{ijk}), ~~~~ k=1,...,6;
\end{eqnarray*}
for $i=1,2$, $j=1,2$, and $l=1,..,6176$. Note that the dimensions of
the GO terms $X_{ijk}^l$ range from 1 to 5158. The MDD-based test is
implemented for the GO terms with dimension $p_l \ge 5$. The
remaining GO terms are tested using a simple $F$-test.

\begin{figure}[!ht]
    \centering
   \includegraphics[width=4in]{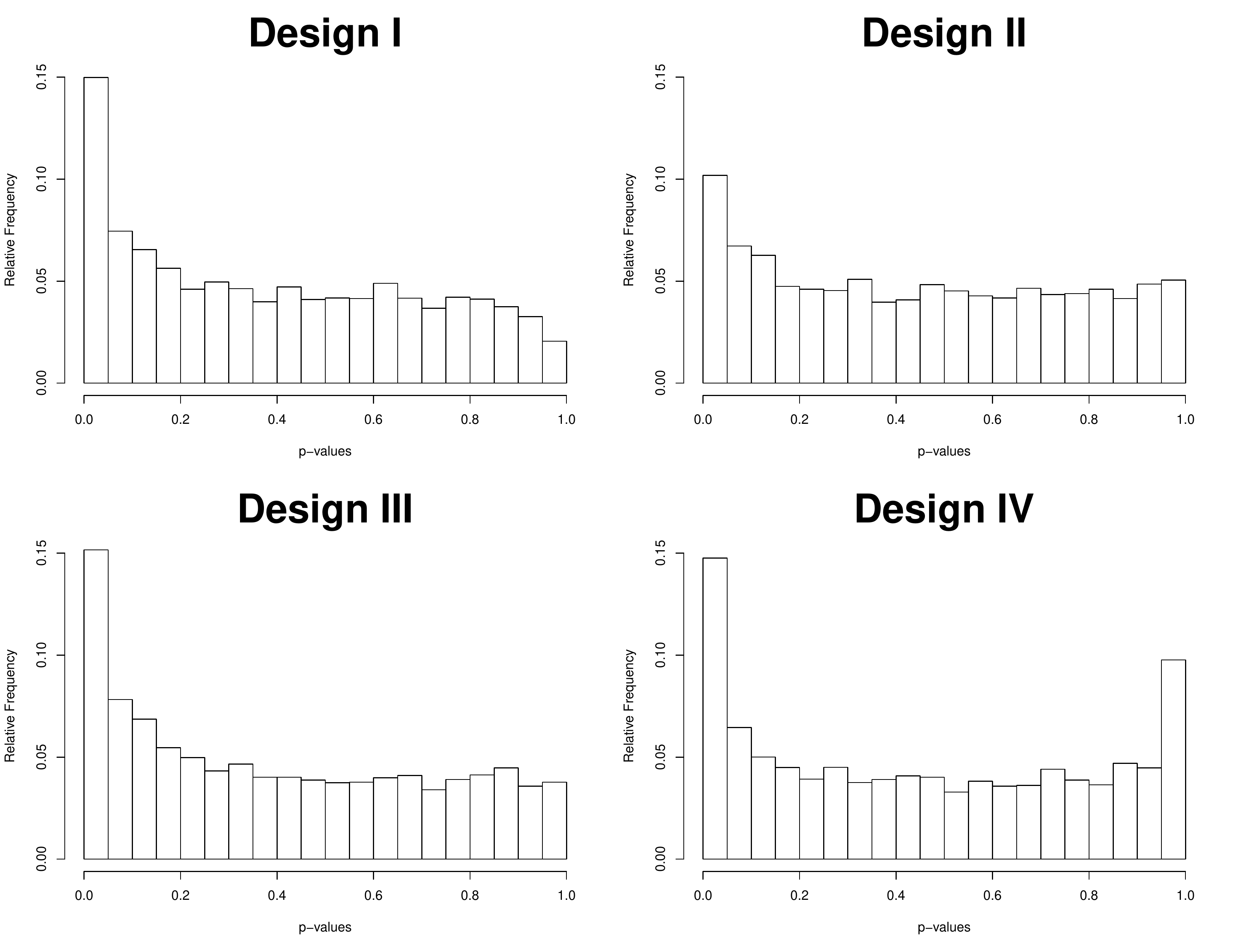}
    \caption{Histograms of the p-values on all gene sets}
    \label{fig:gene}
\end{figure}

We draw the histograms of p-values for all the GO terms as shown in
Figure \ref{fig:gene}. From the histograms we can observe that the
p-values from Design I and Design III are similar; while the results
from Design II has lower portion of small p-values than the other
three designs; Design IV has more large p-values relatively to
others.

\begin{figure}[ht!]
      \centering
      \includegraphics[width=4.5in]{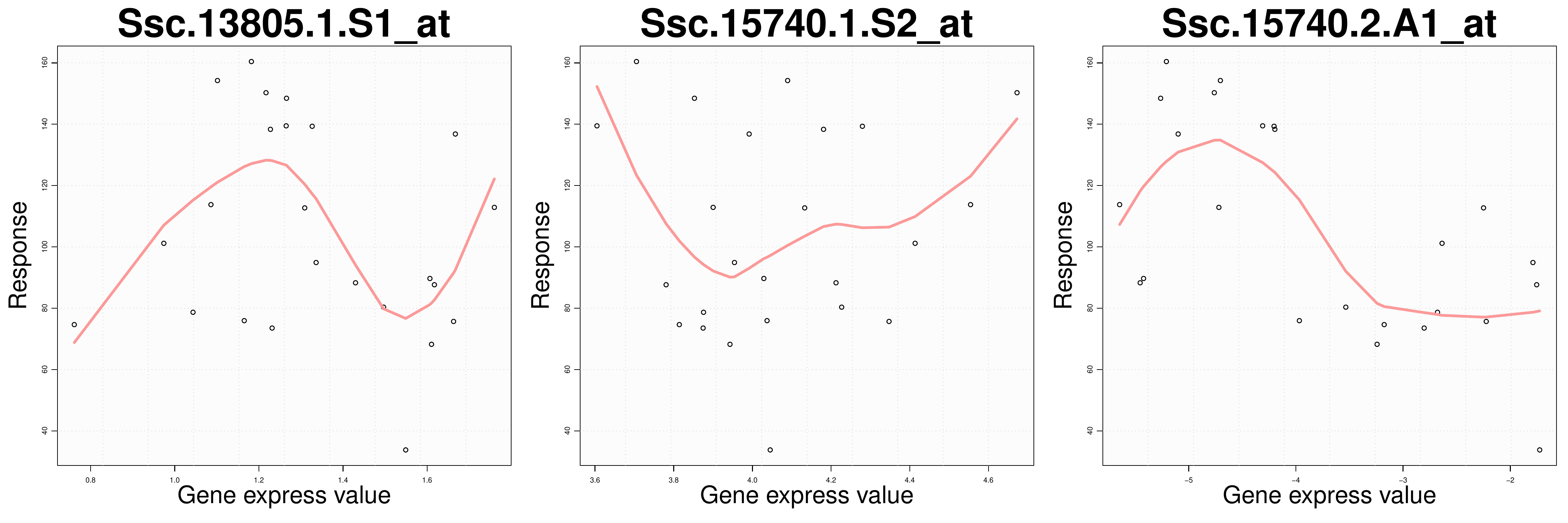}\\
      \includegraphics[width=3.2in]{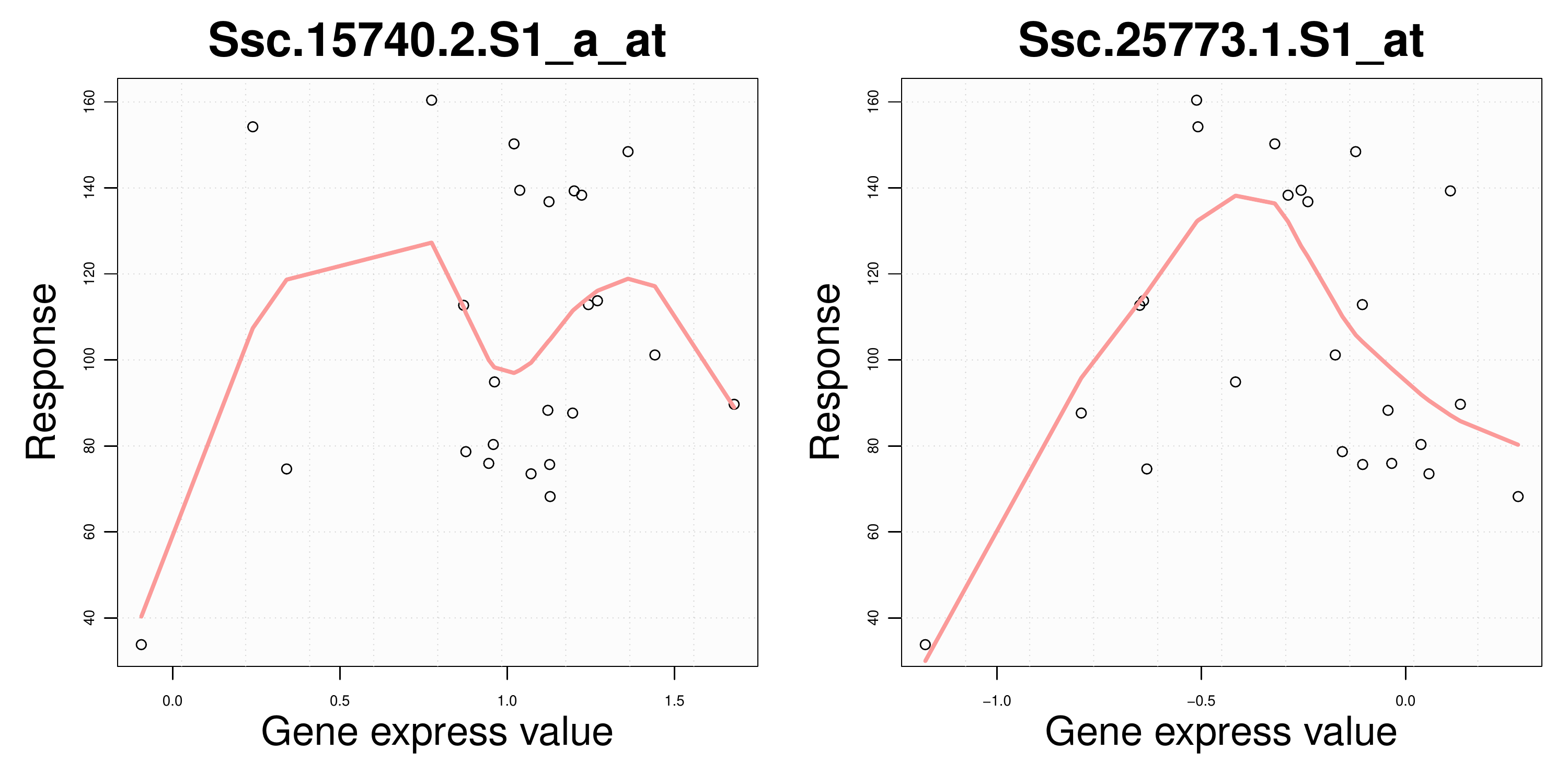}
      \caption{Scatter plots of the $T_3$ measurements versus five gene express values in the GO term 0043536; the red lines correspond to the LOESS fitting.}
      \label{fig:nonlinear}
      \end{figure}

We then use the Benjamini--Hochberg step-up procedure to control the false discovery rate at level $\alpha=0.05$. Namely, for $m$ hypothesis, we find the largest integer $k$ such that $P_{(k)} \leq \frac{k}{m} \alpha$, where $P_{(k)}$ is the ordered p-values from the total $m$ hypothesis tests. Then we reject the null hypothesis for all $H_{(i)}$ for $i = 1, \ldots, k$, where $H_{(i)}$ is the null hypothesis corresponds to the ordered p-value $P_{(i)}$. After controlling the false discovery rate for the p-values, we
find 5 gene sets that are signified as significant under all four
designs. They are GO:0032012, GO:0005086, GO:0043536, GO:0005161 and
GO:0045095. Compared with Zhong and Chen (2011), gene sets
GO:0032012, GO:0005086 are among the three significant gene sets
they found using their method. To demonstrate the non-linear dependence of the GO terms we found, we choose GO:0043536, which contains 5 genes, and present the scatter plot of the response versus the five gene expression values and also fit a local polynomial regression line (LOESS) with degree 2 in Figure \ref{fig:nonlinear}. From the plot, we do observe some non-linear dependence between the response and the five genes.

Our null considers the conditional
mean independence of the $T_3$ measurement with the gene sets, which
includes the null hypothesis in Zhong and Chen (2011) under the
linear model assumption. Thus the above finding seems reasonable
since we expect to detect more significant gene sets as our test is
more powerful than ZC test when the gene set contributes to the mean
of $T_3$ in a non-linear fashion as shown in our simulations.

% % % % % % % % % % % % % % % % % % % % % % % % % %
% % % % % % % % % % % % % % % % % % % % % % % % % %
% % % % % % % % % % % EXAMPLE 1 % % % % % % % % % %
% % % % % % % % % % % % % % % % % % % % % % % % % %
% % % % % % % % % % % % % % % % % % % % % % % % % %

\newpage

\begin{table}[h!]\footnotesize
    \centering
    \caption{Empirical sizes and powers of the MDD-based test for conditional mean independence and the ZC test at significance levels 5\% and 10\% for Example \ref{ex21}.}
    \label{table:eg21}
    \begin{tabular}{ccccccc}
        \toprule
        &&&\multicolumn{2}{c}{mdd}&\multicolumn{2}{c}{ZC}\\
        error & case & $p$ & 5\% & 10\% & 5\% & 10\% \\
        \hline

        \multirow{9}{*}{$N(0,1)$} & \multirow{3}{*}{$H_0$ } &   50 &  0.078 & 0.112 & 0.075 & 0.098 \\
        &  &  100 & 0.075&  0.111 & 0.078 & 0.109 \\
        &  &  200 & 0.065 & 0.091 & 0.062 & 0.089 \\
        \cline{3-7}
        & \multirow{3}{*}{\parbox{1.2cm}{\centering non-sparse $H_a$}} &  50 & 0.575 & 0.635 & 0.595 & 0.640 \\
        &   &  100 & 0.842 & 0.879 & 0.862 & 0.887 \\
        &   &  200 & 0.982 & 0.989 & 0.986 & 0.992 \\
        \cline{3-7}
        &\multirow{3}{*}{\parbox{1.2cm}{\centering sparse $H_a$}} &    50 & 0.189 & 0.243 & 0.199 & 0.248 \\
        &  &  100 & 0.183 & 0.230 & 0.191 & 0.238 \\
        &  &  200 & 0.186 & 0.224 & 0.187 & 0.224 \\
        \hline
        \multirow{9}{*}{$t_3$} & \multirow{3}{*}{$H_0$ } &  50 & 0.084 & 0.110 & 0.080 & 0.104 \\
        & &  100 & 0.075 & 0.106 & 0.072 & 0.097 \\
        & &  200 & 0.082 & 0.112& 0.063 & 0.091 \\
        \cline{3-7}
        & \multirow{3}{*}{\parbox{1.2cm}{\centering non-sparse $H_a$}} &  50 &0.619 & 0.656 & 0.620 & 0.669 \\
        & &  100 & 0.848 & 0.877 & 0.829 & 0.865 \\
        & &  200 & 0.949 & 0.956 & 0.954 & 0.957 \\
        \cline{3-7}
        &\multirow{3}{*}{\parbox{1.2cm}{\centering sparse $H_a$}} & 50 & 0.234 & 0.278 & 0.226 & 0.275 \\
        & &  100 & 0.218 & 0.254 & 0.212 & 0.253 \\
        & &  200 & 0.226 & 0.269 & 0.215 & 0.262 \\
        \hline
        \multirow{9}{*}{$\chi^2_1-1$} & \multirow{3}{*}{$H_0$ }&
        50 & 0.083 & 0.114 & 0.083 & 0.110 \\
        & &  100 & 0.061 & 0.096 & 0.058 & 0.086 \\
        & &  200 & 0.058 & 0.084 & 0.047 & 0.078 \\
        \cline{3-7}
        & \multirow{3}{*}{\parbox{1.2cm}{\centering non-sparse $H_a$}} & 50 & 0.746 & 0.791 & 0.744 & 0.788 \\
        & &  100 & 0.908 & 0.930 & 0.916 & 0.932 \\
        & &  200 & 0.978 & 0.988 & 0.978 & 0.987 \\
        \cline{3-7}
        &\multirow{3}{*}{\parbox{1.2cm}{\centering sparse $H_a$}} &   50 & 0.272 & 0.319 & 0.264 & 0.317 \\
        & &  100 & 0.265 & 0.314 & 0.256 & 0.307 \\
        & &  200 & 0.250 & 0.313 & 0.246 & 0.297 \\
        \hline
    \end{tabular}
\end{table}

% % % % % % % % % % % % % % % % % % % % % % % % % %
% % % % % % % % % % % % % % % % % % % % % % % % % %
% % % % % % % % % % % EXAMPLE 2 % % % % % % % %
% % % % % % % % % % % % % % % % % % % % % % % % % %
% % % % % % % % % % % % % % % % % % % % % % % % % %
\newpage

\begin{table}[h!]
\footnotesize
    \centering
    \caption{\small Empirical powers of the MDD-based test for conditional mean independence and the ZC test at significance levels 5\% and 10\% for Examples \ref{ex22}.}
    \label{table:eg22}
    \begin{tabular}{ccccccc}
        \toprule
        &&&\multicolumn{2}{c}{mdd}&\multicolumn{2}{c}{ZC}\\
        error & case & $p$ &  5\% & 10\% & 5\% & 10\%\\
        \hline
        \multirow{6}{*}{\scriptsize  $N(0,1)$} &  \multirow{3}{*}{\parbox{1.2cm}{\centering non-sparse $H_a$}} &   50 & 0.997 & 0.999 & 0.077 & 0.111 \\
        &  &  100   & 0.986 & 1.000 & 0.092 & 0.126 \\
        &  &  200  & 0.973 & 0.997 & 0.120 & 0.155 \\
        \cline{3-7}
        &\multirow{3}{*}{\parbox{1.2cm}{\centering sparse $H_a$}} &   50 & 0.933 & 0.968 & 0.020 & 0.036 \\
        &  &  100  & 0.874 & 0.933 & 0.020 & 0.040 \\
        &  &  200  & 0.825 & 0.911 & 0.019 & 0.036 \\
        \hline
        \multirow{6}{*}{\scriptsize  $t_3$} & \multirow{3}{*}{\parbox{1.2cm}{\centering non-sparse $H_a$}} &  50 &  0.789 & 0.888 & 0.074 & 0.109 \\
        & &  100  & 0.685 & 0.821 & 0.073 & 0.106 \\
        & &  200  & 0.593 & 0.750 & 0.095 & 0.119 \\
        \cline{3-7}
        &\multirow{3}{*}{\parbox{1.2cm}{\centering sparse $H_a$}} & 50  & 0.930 & 0.968 & 0.019 & 0.037 \\
        & &  100  & 0.876 & 0.932 & 0.020 & 0.036 \\
        & &  200  & 0.823 & 0.910 & 0.019 & 0.035 \\
        \hline
        \multirow{6}{*}{\scriptsize $\chi^2_1-1$} & \multirow{3}{*}{\parbox{1.2cm}{\centering non-sparse $H_a$}} &  50  & 0.854 & 0.912 & 0.095 & 0.116 \\
        & &  100  & 0.777 & 0.874 & 0.095 & 0.122 \\
        & &  200  & 0.727 & 0.833 & 0.079 & 0.108 \\
        \cline{3-7}
        &\multirow{3}{*}{\parbox{1.2cm}{\centering sparse $H_a$}} &   50  & 0.933 & 0.967 & 0.019 & 0.033 \\
        & &  100 & 0.874 & 0.929 & 0.021 & 0.036 \\
        & &  200 & 0.822 & 0.912 & 0.017 & 0.036 \\
        \hline
    \end{tabular}
\end{table}

% % % % % % % % % % % % % % % % % % % % % % % % % %
% % % % % % % % % % % % % % % % % % % % % % % % % %
% % % % % % % % % % % EXAMPLE 3   % % % % % % % %
% % % % % % % % % % % % % % % % % % % % % % % % % %
% % % % % % % % % % % % % % % % % % % % % % % % % %
\newpage

\begin{table}[htp!]
\footnotesize
    \centering
    \caption{Empirical sizes and powers of the MDD-based test for conditional quantile independence at significance levels 5\% and 10\% for Example \ref{ex23}.}
    \label{table:eg23}
    \begin{tabular}{ccccccccccc}
        \toprule
        &&&\multicolumn{2}{c}{$N(0,1)$}&\multicolumn{2}{c}{$t_3$}&\multicolumn{2}{c}{$\chi^2_1-1$}&\multicolumn{2}{c}{Cauchy(0,1)}\\
        $\tau$ & case & $p$ & 5\% & 10\% & 5\% & 10\% & 5\% & 10\% & 5\% & 10\%  \\
        \hline

        \multirow{9}{*}{0.25} & \multirow{3}{*}{$H_0$ }  &   50 & 0.084 & 0.112 & 0.058 & 0.077 & 0.074 & 0.099 & 0.073 & 0.094 \\
   &  &  100 & 0.073 & 0.104 & 0.067 & 0.094 & 0.071 & 0.099 & 0.069 & 0.101 \\
   &  &  200 & 0.064 & 0.085 & 0.058 & 0.095 & 0.069 & 0.094 & 0.062 & 0.082 \\
        \cline{3-11}
 & \multirow{3}{*}{\parbox{1.2cm}{ \centering non-sparse $H_a$}} &    50 & 0.134 & 0.173 & 0.133 & 0.167 & 0.934 & 0.948 & 0.144 & 0.193 \\
   &  &  100 & 0.148 & 0.192 & 0.164 & 0.199 & 0.952 & 0.965 & 0.225 & 0.271 \\
   &  &  200 & 0.168 & 0.218 & 0.185 & 0.211 & 0.874 & 0.901 & 0.285 & 0.331 \\
        \cline{3-11}
        &\multirow{3}{*}{\parbox{1.2cm}{\centering sparse $H_a$}} &   50 & 0.095 & 0.123 & 0.081 & 0.107 & 0.615 & 0.663 & 0.082 & 0.106 \\
   &  &  100 & 0.091 & 0.126 & 0.081 & 0.101 & 0.612 & 0.670 & 0.082 & 0.120 \\
   &  &  200 & 0.082 & 0.109 & 0.083 & 0.114 & 0.590 & 0.646 & 0.076 & 0.109 \\
        \hline
        \multirow{9}{*}{0.5} & \multirow{3}{*}{$H_0$ }  &   50 & 0.083 & 0.106 & 0.076 & 0.097 & 0.064 & 0.089 & 0.068 & 0.099 \\
   &  &  100 & 0.072 & 0.100 & 0.063 & 0.087 & 0.062 & 0.091 & 0.057 & 0.092 \\
   &  &  200 & 0.073 & 0.100 & 0.053 & 0.073 & 0.063 & 0.088 & 0.071 & 0.095 \\
        \cline{3-11}
        & \multirow{3}{*}{\parbox{1.2cm}{ \centering non-sparse $H_a$}}    &   50 & 0.519 & 0.575 & 0.464 & 0.530 & 0.107 & 0.128 & 0.333 & 0.397 \\
   &  &  100 & 0.764 & 0.819 & 0.722 & 0.774 & 0.122 & 0.154 & 0.568 & 0.636 \\
   &  &  200 & 0.930 & 0.943 & 0.895 & 0.920 & 0.208 & 0.251 & 0.802 & 0.844 \\
        \cline{3-11}
        &\multirow{3}{*}{\parbox{1.2cm}{\centering sparse $H_a$}}  &   50 & 0.163 & 0.209 & 0.153 & 0.190 & 0.077 & 0.102 & 0.120 & 0.160 \\
   &  &  100 & 0.169 & 0.206 & 0.141 & 0.180 & 0.069 & 0.092 & 0.124 & 0.162 \\
   &  &  200 & 0.175 & 0.212 & 0.139 & 0.160 & 0.071 & 0.091 & 0.121 & 0.160 \\
        \hline
        \multirow{9}{*}{0.75} & \multirow{3}{*}{$H_0$ }  &   50 & 0.067 & 0.098 & 0.080 & 0.108 & 0.065 & 0.085 & 0.062 & 0.084 \\
   &  &  100 & 0.057 & 0.085 & 0.078 & 0.098 & 0.073 & 0.095 & 0.076 & 0.111 \\
   &  &  200 & 0.064 & 0.095 & 0.063 & 0.085 & 0.065 & 0.088 & 0.069 & 0.095 \\
        \cline{3-11}
        & \multirow{3}{*}{\parbox{1.2cm}{ \centering non-sparse $H_a$}}    &   50 & 0.929 & 0.953 & 0.839 & 0.877 & 0.352 & 0.409 & 0.606 & 0.674 \\
   &  &  100 & 0.995 & 0.998 & 0.967 & 0.977 & 0.596 & 0.649 & 0.850 & 0.886 \\
   &  &  200 & 1.000 & 1.000 & 1.000 & 1.000 & 0.833 & 0.875 & 0.968 & 0.976 \\
        \cline{3-11}
        &\multirow{3}{*}{\parbox{1.2cm}{\centering sparse $H_a$}}  &   50 & 0.415 & 0.482 & 0.347 & 0.394 & 0.118 & 0.153 & 0.217 & 0.260 \\
   &  &  100 & 0.418 & 0.473 & 0.342 & 0.392 & 0.131 & 0.162 & 0.213 & 0.255 \\
   &  &  200 & 0.390 & 0.447 & 0.324 & 0.370 & 0.123 & 0.157 & 0.210 & 0.262 \\
        \hline
        \multirow{9}{*}{mean} & \multirow{3}{*}{$H_0$ }   &   50 & 0.071 & 0.103 & 0.069 & 0.101 & 0.067 & 0.098 & 0.052 & 0.097 \\
   &  &  100 & 0.074 & 0.101 & 0.075 & 0.097 & 0.071 & 0.101 & 0.060 & 0.103 \\
   &  &  200 & 0.067 & 0.094 & 0.078 & 0.105 & 0.066 & 0.088 & 0.054 & 0.096 \\
        \cline{3-11}
        & \multirow{3}{*}{\parbox{1.2cm}{ \centering non-sparse $H_a$}}   &   50 & 0.414 & 0.466 & 0.260 & 0.306 & 0.206 & 0.280 & 0.075 & 0.117 \\
   &  &  100 & 0.549 & 0.602 & 0.323 & 0.381 & 0.278 & 0.357 & 0.096 & 0.139 \\
   &  &  200 & 0.605 & 0.667 & 0.408 & 0.464 & 0.390 & 0.487 & 0.115 & 0.159 \\
        \cline{3-11}
        &\multirow{3}{*}{\parbox{1.2cm}{\centering sparse $H_a$}} &   50 & 0.181 & 0.233 & 0.136 & 0.180 & 0.106 & 0.140 & 0.056 & 0.106 \\
   &  &  100 & 0.171 & 0.207 & 0.119 & 0.158 & 0.105 & 0.142 & 0.064 & 0.109 \\
   &  &  200 & 0.178 & 0.211 & 0.138 & 0.171 & 0.100 & 0.141 & 0.059 & 0.101 \\
   \hline
\end{tabular}
\end{table}

\newpage

\begin{table}[h!]
\footnotesize
    \centering
    \caption{Size comparison for the proposed test using normal approximation (mdd) and the wild bootstrap approximation for Example 4.4 of the paper at the 5 \% nominal level.}
    \label{table:eg54boot}
    \begin{tabular}{cccccccccc}
        \toprule
        &&\multicolumn{2}{c}{$N(0,1)$}&\multicolumn{2}{c}{$t_3$}&\multicolumn{2}{c}{Cauchy(0,1)} &\multicolumn{2}{c}{$\chi^2_1-1$}\\
        $\tau$  & $p$ & mdd & Boot & mdd & Boot  & mdd & Boot  & mdd & Boot   \\
        \hline
       \multirow{3}{*}{0.25} &   50 & 0.084 & 0.061 & 0.058 & 0.043 & 0.073 & 0.054 & 0.074 & 0.055 \\
          &  100 & 0.073 & 0.053 & 0.067 & 0.048 & 0.069 & 0.051 & 0.071 & 0.053 \\
          &  200 & 0.064 & 0.043 & 0.058 & 0.042 & 0.062 & 0.035 & 0.069 & 0.045 \\
          \hline
          \multirow{3}{*}{0.50}    &   50 & 0.083 & 0.057 & 0.076 & 0.054 & 0.068 & 0.044 & 0.064 & 0.048 \\
          &  100 & 0.072 & 0.057 & 0.063 & 0.047 & 0.057 & 0.046 & 0.062 & 0.046 \\
          &  200 & 0.073 & 0.054 & 0.053 & 0.040 & 0.071 & 0.051 & 0.063 & 0.045 \\
               \hline
          \multirow{3}{*}{0.75}    &   50 & 0.067 & 0.051 & 0.080 & 0.059 & 0.062 & 0.050 & 0.065 & 0.041 \\
          &  100 & 0.057 & 0.044 & 0.078 & 0.052 & 0.076 & 0.052 & 0.073 & 0.056 \\
          &  200 & 0.064 & 0.050 & 0.063 & 0.042 & 0.069 & 0.051 & 0.065 & 0.047 \\
               \hline
    \end{tabular}
\end{table}

% % % % % % % % % % % % % % % % % % % % % % % % % %
% % % % % % % % % % % % % % % % % % % % % % % % % %
% % % % % % % %  factorial design  ex4 % %% % % % %
% % % % % % % % % % % % % % % % % % % % % % % % % %
% % % % % % % % % % % % % % % % % % % % % % % % % %
\newpage

\begin{table}[ht!]\footnotesize
    \centering
    \caption{Empirical sizes and powers of the MDD-based test for conditional mean independence and the ZC test at significance levels 5\% and 10\% for Example \ref{ex55}.}
    \label{table:eg55}
    \begin{tabular}{ccccccccccc}
        \hline
        &&&\multicolumn{4}{c}{Case 1}&\multicolumn{4}{c}{Case 2}\\
        &&&\multicolumn{2}{c}{mdd} &\multicolumn{2}{c}{ZC} &\multicolumn{2}{c}{mdd} &\multicolumn{2}{c}{ZC}\\
        & $n$ & $p$ & 5\% & 10\% & 5\% & 10\% & 5\% & 10\% & 5\% & 10\% \\
        \hline
        \multirow{9}{*}{$H_0$ } &    30 &  100 & 0.054 &  0.098  & 0.061 & 0.099 & 0.074 & 0.115 & 0.083 & 0.122\\
        &   30 &  150 & 0.056 & 0.105 & 0.061 & 0.117 & 0.060 & 0.102 & 0.067 & 0.114\\
        &   30 &  200 & 0.055 & 0.099 & 0.055 & 0.114 & 0.064 & 0.108 & 0.063 & 0.107\\
        &   50 &  100 & 0.064 & 0.104 & 0.057 & 0.103 & 0.062 & 0.098 & 0.067 & 0.102\\
        &   50 &  150 & 0.054 & 0.098 & 0.062 & 0.094 & 0.060 & 0.094 & 0.058 & 0.094\\
        &   50 &  200 & 0.067 & 0.110 & 0.062 & 0.112 & 0.071 & 0.108 & 0.070 & 0.105\\
        &   70 &  100 & 0.068 & 0.114 & 0.058 & 0.109 & 0.066 & 0.100 & 0.069 & 0.106\\
        &   70 &  150 & 0.075 & 0.109 & 0.073 & 0.112 & 0.065 & 0.102 & 0.067 &  0.106 \\
        &   70 &  200 & 0.061 & 0.098 & 0.063 & 0.099 & 0.058 & 0.104 & 0.070 &  0.107 \\
        \hline
        \multirow{9}{*}{\parbox{1.2cm}{ \centering non-sparse $H_a$}}  &  30 &  100 & 0.726 & 0.813 & 0.744 & 0.815 & 0.412 & 0.546 & 0.196 & 0.237\\
        &   30 &  150 & 0.566 & 0.674 & 0.577 & 0.682 & 0.459 & 0.596 & 0.226 & 0.266\\
        &   30 &  200 & 0.452 & 0.575 & 0.464 & 0.582 & 0.508 & 0.645 & 0.227 & 0.266\\
        &   50 &  100 & 0.972 & 0.986 & 0.975 & 0.984 & 0.788 & 0.908 & 0.194 & 0.236\\
        &   50 &  150 & 0.889 & 0.934 & 0.891 & 0.927 & 0.871 & 0.959 & 0.187 & 0.233\\
        &   50 &  200 & 0.786 & 0.855 & 0.783 & 0.849 & 0.903 & 0.975 & 0.251 & 0.281\\
        &   70 &  100 & 0.999 & 0.999 & 0.989 & 0.997 & 0.991 & 0.998 & 0.189 & 0.230\\
        &   70 &  150 & 0.978 & 0.986 & 0.963 & 0.977 & 0.998 & 0.999 & 0.200 & 0.236\\
        &   70 &  200 & 0.948 & 0.971 & 0.918 & 0.954 & 1.000 & 1.000 & 0.210 & 0.242\\
        \hline
        \multirow{9}{*}{\parbox{1.2cm}{\centering sparse $H_a$}} &  30 &  100 & 0.369 & 0.478 & 0.389 & 0.478 & 0.150 & 0.206 & 0.102 & 0.145\\
        &   30 &  150 & 0.283 & 0.396 & 0.301 & 0.395 & 0.128 & 0.185 & 0.091 & 0.127 \\
        &   30 &  200 & 0.229 & 0.326 & 0.239 & 0.339 & 0.159 & 0.233 & 0.104 & 0.160\\
        &   50 &  100 & 0.694 & 0.782 & 0.677 & 0.756 & 0.194 & 0.267 & 0.097 & 0.145\\
        &   50 &  150 & 0.523 & 0.622 & 0.503 & 0.617 & 0.172 & 0.240 & 0.097 & 0.133\\
        &   50 &  200 & 0.413 & 0.535 & 0.404 & 0.512 & 0.188 & 0.260 & 0.093 & 0.139\\
        &   70 &  100 & 0.875 & 0.918 & 0.838 & 0.882 & 0.251 & 0.356 & 0.113 & 0.145\\
        &   70 &  150 & 0.720 & 0.805 & 0.658 & 0.750 & 0.224 & 0.340 & 0.108 & 0.145\\
        &   70 &  200 & 0.611 & 0.699 & 0.545 & 0.638 & 0.233 & 0.312 & 0.111 & 0.161 \\
        \hline
    \end{tabular}
\end{table}

% % % % % % % % % % % % % % % % % % % % % % % % % %
% % % % % % % % % % % % % % % % % % % % % % % % % %
% % % % % % % % % % % EXAMPLE 4   % % % % % % % %
% % % % % % % % % % % % % % % % % % % % % % % % % %
% % % % % % % % % % % % % % % % % % % % % % % % % %

\newpage

\begin{table}[ht!]\footnotesize
    \centering
    \caption{Empirical sizes and powers of the MDD-based test for conditional mean independence and the ZC test at significance levels 5\% and 10\% for Example \ref{ex24}.}
    \label{table:eg24}
    \begin{tabular}{ccccccccccc}
        \hline
        &&&\multicolumn{4}{c}{Normal error}&\multicolumn{4}{c}{Gamma error}\\
        &&&\multicolumn{2}{c}{mdd} &\multicolumn{2}{c}{ZC} &\multicolumn{2}{c}{mdd} &\multicolumn{2}{c}{ZC}\\
        & $n$ & $p$ & 5\% & 10\% & 5\% & 10\% & 5\% & 10\% & 5\% & 10\% \\
        \hline
        \multirow{9}{*}{$H_0$ } &      30 &  100 & 0.054 & 0.098 & 0.061 & 0.099  & 0.055 & 0.113 & 0.055 & 0.100 \\
        &   30 &  150 & 0.056 & 0.105 & 0.061 & 0.117 & 0.055 & 0.100 & 0.047 & 0.093 \\
        &   30 &  200 & 0.055 & 0.099 & 0.055 & 0.114  & 0.052 & 0.096 & 0.046 & 0.090 \\
        &   50 &  100 & 0.064 & 0.104 & 0.057 & 0.103 & 0.064 & 0.108 & 0.053 & 0.105\\
        &   50 &  150 & 0.054 & 0.098 & 0.062 & 0.094 & 0.046 & 0.085 & 0.051 & 0.086\\
        &   50 &  200 & 0.067 & 0.110 & 0.062 & 0.112 & 0.057 & 0.098 & 0.052 & 0.092\\
        &   70 &  100 & 0.068 & 0.114 & 0.058 & 0.109 & 0.069 & 0.118 & 0.059 & 0.107\\
        &   70 &  150 & 0.075 & 0.109 & 0.073 & 0.112  & 0.044 & 0.090 & 0.041 & 0.094\\
        &   70 &  200 & 0.061 & 0.098 & 0.063 & 0.099 & 0.048 & 0.080 & 0.048 & 0.090\\
        \hline
        \multirow{9}{*}{\parbox{1.2cm}{ \centering non-sparse $H_a$}}    &      30 &  100 & 1.000 & 1.000 & 1.000 & 1.000  & 0.999 & 0.999 & 0.999 & 0.999\\
        &   30 &  150 & 0.996 & 1.000 & 0.998 & 1.000  & 0.989 & 0.992 & 0.993 & 0.995\\
        &   30 &  200 & 0.994 & 0.998 & 0.995 & 0.998  & 0.982 & 0.989 & 0.987 & 0.993 \\
        &   50 &  100 & 1.000 & 1.000 & 1.000 & 1.000  & 1.000 & 1.000 & 1.000 & 1.000\\
        &   50 &  150 & 1.000 & 1.000 & 1.000 & 1.000  & 1.000 & 1.000 & 1.000 & 1.000\\
        &   50 &  200 & 1.000 & 1.000 & 1.000 & 1.000  & 1.000 & 1.000 & 1.000 & 1.000\\
        &   70 &  100 & 1.000 & 1.000 & 1.000 & 1.000  & 1.000 & 1.000 & 1.000 & 1.000\\
        &   70 &  150 & 1.000 & 1.000 & 1.000 & 1.000  & 1.000 & 1.000 & 1.000 & 1.000 \\
        &   70 &  200 & 1.000 & 1.000 & 1.000 & 1.000  & 1.000 & 1.000 & 1.000 & 1.000 \\
        \hline
        \multirow{9}{*}{\parbox{1.2cm}{\centering sparse $H_a$}}   &      30 &  100 & 0.437 & 0.550 & 0.485 & 0.573 & 0.514 & 0.616 & 0.540 & 0.621 \\
        &   30 &  150 & 0.341 & 0.452 & 0.382 & 0.478 & 0.323 & 0.417 & 0.331 & 0.434  \\
        &   30 &  200 & 0.265 & 0.370 & 0.304 & 0.397  & 0.261 & 0.364 & 0.263 & 0.381  \\
        &   50 &  100 & 0.819 & 0.877 & 0.790 & 0.849  & 0.817 & 0.875 & 0.795 & 0.865\\
        &   50 &  150 & 0.644 & 0.754 & 0.621 & 0.729 & 0.652 & 0.729 & 0.621 & 0.712\\
        &   50 &  200 & 0.532 & 0.640 & 0.505 & 0.608 & 0.573 & 0.688 & 0.561 & 0.658 \\
        &   70 &  100 & 0.938 & 0.969 & 0.918 & 0.947 & 0.952 & 0.970 & 0.916 & 0.945\\
        &   70 &  150 & 0.847 & 0.905 & 0.792 & 0.850 & 0.864 & 0.907 & 0.805 & 0.865\\
        &   70 &  200 & 0.735 & 0.816 & 0.664 & 0.755 & 0.757 & 0.821 & 0.672 & 0.764\\
        \hline
    \end{tabular}
\end{table}

\newpage

\begin{table}[ht!]
\footnotesize
\caption{Empirical sizes and powers from four tests
for Example \ref{eg:ar}} \label{table:ar} \centering
\begin{tabular}{cccccccccc}
\hline
 \multirow{2}{*}{{\bf {\large  ($\rho = 0.1$)}}} &\multirow{2}{*}{p} & \multirow{2}{*}{mdd} &\multirow{2}{*}{ZL}    & \multirow{2}{*}{CL}  &  \multicolumn{5}{c}{MQ ($\lambda_n$)} \\
  \cline{6-10}
  &  &  &  &  & 2 & 4 & 4.5 & 5 & 10 \\
  \hline
 \multirow{3}{*}{ $H_0$ }     &   10 & 0.045 & 0.035 & 0.030 & 0.335 & 0.060 & 0.040 & 0.040 & 0.040 \\
    &   50 & 0.085 & 0.045 & 0.050 & 0.660 & 0.085 & 0.060 & 0.045 & 0.045 \\
    &  150 & 0.020 & 0.085 & 0.015 & 0.835 & 0.195 & 0.100 & 0.080 & 0.075 \\
    \hline
         \multirow{3}{*}{non-sparse $H_a$ }       &   10 & 0.565 & 0.345 & 0.570 & 0.695 & 0.430 & 0.405 & 0.380 & 0.350 \\
             &   50 & 0.250 & 0.105 & 0.220 & 0.665 & 0.200 & 0.150 & 0.120 & 0.095 \\
             &  150 & 0.145 & 0.095 & 0.135 & 0.845 & 0.315 & 0.165 & 0.120 & 0.080 \\
        \hline
    \multirow{3}{*}{sparse $H_a$ }                &   10 & 0.415 & 0.490 & 0.450 & 0.830 & 0.540 & 0.510 & 0.500 & 0.495 \\
                 &   50 & 0.185 & 0.355 & 0.155 & 0.895 & 0.415 & 0.360 & 0.340 & 0.320 \\
                 &  150 & 0.075 & 0.300 & 0.060 & 0.880 & 0.490 & 0.330 & 0.295 & 0.230 \\
   \hline
  \hline
 \multirow{2}{*}{{\bf {\large  ($\rho=0.5$ )}}}  &\multirow{2}{*}{p} & \multirow{2}{*}{mdd} &\multirow{2}{*}{ZL}    & \multirow{2}{*}{CL}  &  \multicolumn{5}{c}{MQ ($\lambda_n$)} \\
  \cline{6-10}
   &  &  &  &  & 2 & 4 & 4.5 & 5 & 10 \\
  \hline
\multirow{3}{*}{$H_0$}    &   10 & 0.040 & 0.020 & 0.035 & 0.285 & 0.045 & 0.035 & 0.030 & 0.030 \\
   &   50 & 0.065 & 0.035 & 0.045 & 0.610 & 0.075 & 0.050 & 0.035 & 0.020 \\
   &  150 & 0.025 & 0.025 & 0.020 & 0.805 & 0.145 & 0.070 & 0.050 & 0.040 \\
     \hline
        \multirow{3}{*}{non-sparse $H_a$ }
           &   10 & 0.975 & 0.900 & 0.980 & 0.985 & 0.960 & 0.945 & 0.935 & 0.900 \\
           &   50 & 0.905 & 0.405 & 0.905 & 0.915 & 0.770 & 0.615 & 0.510 & 0.425 \\
           &  150 & 0.695 & 0.195 & 0.740 & 0.910 & 0.665 & 0.420 & 0.300 & 0.205 \\
       \hline
    \multirow{3}{*}{sparse $H_a$ }              &   10 & 0.485 & 0.545 & 0.455 & 0.815 & 0.595 & 0.565 & 0.545 & 0.540 \\
               &   50 & 0.195 & 0.345 & 0.150 & 0.865 & 0.430 & 0.360 & 0.330 & 0.315 \\
               &  150 & 0.050 & 0.265 & 0.050 & 0.890 & 0.455 & 0.350 & 0.270 & 0.225 \\
   \hline
  \hline
    \multirow{2}{*}{{\bf {\large  ($\rho=0.8$ )}}} &\multirow{2}{*}{p} & \multirow{2}{*}{mdd} &\multirow{2}{*}{ZL}    & \multirow{2}{*}{CL}  &  \multicolumn{5}{c}{MQ ($\lambda_n$)} \\
     \cline{6-10}
     &  &  &  &  & 2 & 4 & 4.5 & 5 & 10 \\
     \hline
     \multirow{3}{*}{$H_0$}  &   10 & 0.075 & 0.030 & 0.035 & 0.175 & 0.085 & 0.060 & 0.045 & 0.045 \\
        &   50 & 0.055 & 0.035 & 0.040 & 0.500 & 0.100 & 0.080 & 0.065 & 0.055 \\
        &  150 & 0.045 & 0.020 & 0.025 & 0.705 & 0.070 & 0.040 & 0.030 & 0.030 \\
     \hline
    \multirow{3}{*}{non-sparse $H_a$ }           &   10 & 1.000 & 1.000 & 1.000 & 1.000 & 1.000 & 1.000 & 1.000 & 1.000 \\
            &   50 & 1.000 & 1.000 & 1.000 & 1.000 & 1.000 & 1.000 & 1.000 & 1.000 \\
            &  150 & 1.000 & 0.940 & 1.000 & 0.995 & 0.995 & 0.995 & 0.985 & 0.930 \\
  \hline
    \multirow{3}{*}{sparse $H_0$ }                  &   10 & 0.520 & 0.625 & 0.485 & 0.770 & 0.660 & 0.640 & 0.610 & 0.595 \\
                &   50 & 0.225 & 0.415 & 0.185 & 0.845 & 0.585 & 0.510 & 0.455 & 0.400 \\
                &  150 & 0.115 & 0.285 & 0.080 & 0.860 & 0.450 & 0.330 & 0.290 & 0.245 \\
        \hline
  \hline
  \multirow{2}{*}{{\bf {\large  ($\rho=-0.5$ )}}} &\multirow{2}{*}{p} & \multirow{2}{*}{mdd} &\multirow{2}{*}{ZL}    & \multirow{2}{*}{CL}  &  \multicolumn{5}{c}{MQ ($\lambda_n$)} \\
   \cline{6-10}
   &  &  &  &  & 2 & 4 & 4.5 & 5 & 10 \\
   \hline
    \multirow{3}{*}{$H_0$ }  &   10 & 0.060 & 0.035 & 0.045 & 0.370 & 0.055 & 0.055 & 0.055 & 0.055 \\
       &   50 & 0.075 & 0.085 & 0.060 & 0.610 & 0.105 & 0.070 & 0.070 & 0.060 \\
       &  150 & 0.025 & 0.080 & 0.030 & 0.845 & 0.205 & 0.145 & 0.115 & 0.095 \\
\hline
    \multirow{3}{*}{non-sparse $H_a$ }
       &   10 & 0.080 & 0.100 & 0.065 & 0.400 & 0.115 & 0.110 & 0.090 & 0.095 \\
       &   50 & 0.080 & 0.075 & 0.055 & 0.590 & 0.120 & 0.090 & 0.075 & 0.075 \\
       &  150 & 0.045 & 0.095 & 0.045 & 0.835 & 0.225 & 0.150 & 0.095 & 0.085 \\

\hline
    \multirow{3}{*}{sparse $H_0$ }         &   10 & 0.540 & 0.570 & 0.510 & 0.845 & 0.625 & 0.590 & 0.580 & 0.580 \\
       &   50 & 0.200 & 0.345 & 0.195 & 0.875 & 0.445 & 0.370 & 0.350 & 0.320 \\
       &  150 & 0.125 & 0.310 & 0.070 & 0.930 & 0.490 & 0.335 & 0.290 & 0.255 \\
      \hline
\end{tabular}
\end{table}

\newpage

\begin{table}[ht!]
\footnotesize
\caption{Empirical sizes and powers from four tests
for Example \ref{eg:indep}} \label{table:indep} \centering
\begin{tabular}{cccccccccc}
\hline
 \multirow{2}{*}{{\bf {\large  (i)}}}  &\multirow{2}{*}{p} & \multirow{2}{*}{mdd} &\multirow{2}{*}{ZL}    & \multirow{2}{*}{CL}  &  \multicolumn{5}{c}{MQ ($\lambda_n$)} \\
  \cline{6-10}
   &  &  &  &  & 2 & 4 & 4.5 & 5 & 10 \\
  \hline
\multirow{3}{*}{$H_0$}     &   10 & 0.055 & 0.080 & 0.065 & 0.360 & 0.105 & 0.085 & 0.070 & 0.070 \\
   &   50 & 0.080 & 0.060 & 0.065 & 0.680 & 0.095 & 0.070 & 0.070 & 0.050 \\
   &  150 & 0.060 & 0.075 & 0.035 & 0.865 & 0.245 & 0.130 & 0.090 & 0.070 \\
    \hline
        \multirow{3}{*}{non-sparse $H_a$ }   &   10 & 1.000 & 0.910 & 1.000 & 0.990 & 0.970 & 0.955 & 0.915 & 0.900 \\
       &   50 & 0.940 & 0.370 & 0.975 & 0.905 & 0.700 & 0.550 & 0.470 & 0.375 \\
       &  150 & 0.615 & 0.160 & 0.720 & 0.920 & 0.540 & 0.355 & 0.215 & 0.175 \\
       \hline
    \multirow{3}{*}{sparse $H_a$ }     &   10 & 1.000 & 1.000 & 1.000 & 1.000 & 1.000 & 1.000 & 1.000 & 1.000 \\
   &   50 & 0.910 & 1.000 & 0.945 & 1.000 & 1.000 & 1.000 & 1.000 & 1.000 \\
   &  150 & 0.620 & 1.000 & 0.715 & 1.000 & 1.000 & 1.000 & 1.000 & 1.000 \\
   \hline
  \hline
  \multirow{2}{*}{{\bf {\large  (ii)}}} &\multirow{2}{*}{p} & \multirow{2}{*}{mdd} &\multirow{2}{*}{ZL}    & \multirow{2}{*}{CL}  &  \multicolumn{5}{c}{MQ ($\lambda_n$)} \\
   \cline{6-10}
   &  &  &  &  & 2 & 4 & 4.5 & 5 & 10 \\
   \hline
    \multirow{3}{*}{non-sparse $H_a$ }    &   10 & 0.925 & 0.080 & 0.095 & 0.370 & 0.095 & 0.065 & 0.050 & 0.050 \\
         &   50 & 0.395 & 0.035 & 0.085 & 0.690 & 0.080 & 0.060 & 0.035 & 0.030 \\
         &  150 & 0.240 & 0.065 & 0.065 & 0.850 & 0.265 & 0.125 & 0.080 & 0.070 \\
\hline
  \multirow{3}{*}{sparse $H_0$ }
    &   10 & 0.890 & 0.095 & 0.100 & 0.275 & 0.095 & 0.080 & 0.065 & 0.065 \\
      &   50 & 0.385 & 0.060 & 0.055 & 0.595 & 0.085 & 0.050 & 0.040 & 0.020 \\
      &  150 & 0.215 & 0.090 & 0.045 & 0.790 & 0.225 & 0.120 & 0.090 & 0.060 \\

      \hline
\end{tabular}
\end{table}